%% file: Main.tex
\documentstyle{amsppt}
\input H+Wbook.macrosnew
\SetRokickiEPSFSpecial
\input imsmark.tex

\input option_keys
\magnification 1200

\voffset -.3in
\hoffset 0pt
\hsize=6.5 true in
\vsize=9 true in
\overfullrule 0pt

\def\peqno #1{\eqno(#1)}
\parskip 5 pt

\def \C{{\Bbb C}}
\def \N {{\Bbb N}}
\define \Proj{{\Bbb P}}
\def \indlim_#1{\underset {\underset #1 \to \rightarrow} \to {\lim}\ }

\def \ratto {\sim\kern-2pt\sim\kern-2pt>}
\def \Irr {\text{Irr}}
\def \orb{{\Cal B}_\R^+}
\def \orbstar{{\Cal B}_\R^*}
\def \orbtwostar{{\Cal B}_\R^{**}}
\def \pvec#1#2{\pmatrix #1\\#2\endpmatrix}

\def \detideal{2.1}
\def \eqbeforezerodiv{2.2}
\def \eqafterzerodiv{2.3}
\def \eqviaVeronese{2.4}
\def \blowupstillsingeq{2.5}

\def \henaffine {3.1}
\def \henhomogeneous {3.2}
\def \veryfirstblowup {3.3}
\def \firstkblowups {3.4}
\def \successiveYs {3.5}
\def \furthersimpleblowups {3.6}
\def \firstsimpleblowup {3.7}

\def \maptohomologyofbu{5.1}

\def \homologyasdirlim{5.2}
\def \exactseqforhomologyeq{5.3}
\def \intersectformseq {5.4}

\def \baddefeltseq {6.1}
\def \gooddefeltseq {6.2}
\def \natprimeend {6.3}
\def \eqafterzerodivagain {6.4}

\def \paramdesceq {7.1}
\def \extendtoblupordeq {7.2}
\def \circactZeq {7.3}
\def \firstwriteofh {7.4}
\def \secondwriteofh{7.5}

\def \defhoninsideeq {7.6}
\def \firstcompcondeq {7.7}
\def \honDone {7.8}
\def \honDtwo {7.9}
\def \honEone {7.10}
\def \honEtwo {7.11}
\def \hextonWone {7.12}
\def \hextonBone {7.13}
\def \hextonBtwo {7.14}
\def \hextonWtwo {7.15}

\def \stereographic{8.1}
\def \stereoinpolar{8.2}
\def \circleaction {8.3}

\def \argumentsconseq{9.1}
\def \argumentsaddeq{9.2}

\def \homologyclassesoffirstblowups{10.1}



\def \firstblowuppic{1}
\def \secondblowuppic{2}
\def \thirdblowuppic{3}
\def \divatinfty{4}
\def \divatinftywithmap{5}

\def \infinitedivisorfig{6}

\def \curveanddeform{7}
\def \finitedivisor{8}
\def \selfintfig{9}
\def \matcompfirstline{10}
\def \matcompsecondline{11}
\def \matcompthirdline{12}
\def \matcompfourthline{13}
\def \graphoftau{14}
\def \openmanfig{15}
\def \graphsetas{16}
\def \picofW{17}
\def \graphofetasinside{18}
\def \threelines{19}
\def \FigVIIetasfordoublepts{20}
\def \picofB{21}
\def \graphaofetasfordblepts{22}
\def \twohopf{23}
\def \rightsecondblowup{24}
\def \secblowup{25}
\def \whichHopfright{26}
\def \firstdblowuppic{27}
\def \toriasrotpics{28}
\def \baseSeifertfig{29}
\def \blowupsforandback{30}
\def \linkingnumber{31}
\def \linkedcirclespic{32}
\def \compHdiv{33}



\define \hyperanddiv{II.1}
\define \divisordef{II.2}
\define \blowupdef{II.3}
\define \blowuppointonsurf{II.4}
\define \therealcase{II.5}
\define \zerodivs{II.6}
\define \universalproperty {II.7}
\define \natblowups {II.8}


\def \pointsofindetlemma{III.1}
\def \computehorribleblowup{III.2}
\def \HenonWellDefined{III.3}
\def \tildeHasratmap{III.4}

\def \wheregraphsmooth {IV.1}
\def \singsingraphs {IV.2}
\def \graphHsmooth {IV.3}
\def \graphasfiberedprod {IV.4}

\def \sequencespace {IV.5}
\def \pointsofinfproduct {IV.6}
\def \badsequencespace {IV.7}
\def \infproductgoodpointssmooth {IV.8}


\def \indlimexample {V.1}
\def \indlimprop {V.2}
\def \indlimextw {V.3}
\def \removepoints {V.4}
\def \removepointsindimthree {V.5}
\def \homologyofblowups {V.6}
\def \blowupsindimoneandthree {V.7}
\def \mapinducedonhomology {V.8}
\def \howindlim {V.9}
\def \DivInfNM {V.10}
\def \homologyoffinblowups {V.11}
\def \computeinclusion {V.12}
\def \embedinproduct {V.13}
\def \exactsequenceforhomology {V.14}
\def \exactsequenceforpair {V.15}
\def \computeintsonvs {V.16}
\def \signthm {V.17}
\def \slightstrengthening {V.18}


\def \axesexampleasrealsubset{VI.1}
\def \primeenddef{VI.2}
\def \polarcoordexample{VI.3}
\def \naturalityofrealblowup{VI.4}
\def \firstdescofrealblowup{VI.5}
\def \blowupandtubnbhd{VI.6}
\def \orientrealblowup{VI.7}
\def \parametrizebyVeronese{VI.8}
\def \polardesc{VI.9}
\def \descoffibersrealblowups{VI.10}
\def \blowupofsmoothcurves{VI.11}
\def \comdescofrealblowups{VI.12}


\def \natylemmaontubneighborhoods{VII.1}
\def \choosetaunaught{VII.2}
\def \divrightasideal{VII.3}
\def \compsimpblowupontorus{VII.4}
\def \blowupofordpoint{VII.5}
\def \controlledblowupordpt{VII.6}
\def \pitildeexistsfordblepoint{VII.7}
\def \compmapontorofdblept{VII.8}
\def \blowupdoublepoint{VII.9}
\def \controlblowupdblept{VII.10}
\def \topofprojlim{VII.11}


\def \hopfisom {VIII.1}
\def \modifiedHopf {VIII.2}

\def \computelinkingnumber{VIII.3}


\def \realextensiontheorem {IX.1}
\def \mainthmrealblowup {IX.2}
\def \computeanglesunderpitilde {IX.3}


\def \solenoidalprop {X.1}
\def \solenoidalmapsconjtotau {X.2}
\def \computepiplusandminus {X.3}
\def \Hatcherthreepointthree {X.4}
\def \seifertfibrationunique {X.5}
\def \surfacesinseifert {X.6}
\def \toroidaldecompofsphere {X.7}
\def \conjugacyofexthens {X.8}
\def \comppisonB {X.9}
\def \comppisonothertori {X.10}
\def \piscoincide {X.11}
\def \toriintersect {X.12}


\def \milnorconj{XI.1}
\def \closureofgraph {XI.2}
\def \structcompactificationforcomp{XI.3}


\def \Douady{Dou} 
\def \Sha {S2} 
\def \Griffiths {Gri} 
\def \Hir{Hir}

\def \BSone{BS1}


\SBIMSMark{1997/11}{September, 1997}{}
\input Veselov.introduction

\input universalblowups

\input Veselov.welldefined.new

\input graphs+infiniteproducts

\input Veselov.compactifications.two

\input Veselov.HenonExample.new

\input Veselov.milnorconj

\newpage
\input bibliography

\end

%% file: imsmark.tex
\def\SBIMSMark#1#2#3{
 \font\SBF=cmss10 at 10 true pt
 \font\SBI=cmssi10 at 10 true pt
 \setbox0=\hbox{\SBF Stony Brook IMS Preprint \##1}
 \setbox2=\hbox to \wd0{\hfil \SBI #2}
 \setbox4=\hbox to \wd0{\hfil \SBI #3}
 \setbox6=\hbox to \wd0{\hss
             \vbox{\hsize=\wd0 \parskip=0pt \baselineskip=10 true pt
                   \copy0 \break%
                   \copy2 \break%
                   \copy4 \break}}
 \dimen0=\ht6   \advance\dimen0 by \vsize \advance\dimen0 by 8 true pt
                \advance\dimen0 by -\pagetotal
 \dimen2=\hsize \advance\dimen2 by .25 true in
%
%
  \openin2=publishd.tex
  \ifeof2\setbox0=\hbox to 0pt{}
  \else 
     \setbox0=\hbox to 3.1 true in{
                \vbox to \ht6{\hsize=3 true in \parskip=0pt  \noindent  
                \input publishd.tex 
                \vfill}}
  \fi
  \closein2
  \ht0=0pt \dp0=0pt
 \ht6=0pt \dp6=0pt
 \setbox8=\vbox to \dimen0{\vfill \hbox to \dimen2{\copy0 \hss \copy6}}
 \ht8=0pt \dp8=0pt \wd8=0pt
 \copy8
 \message{*** Stony Brook IMS Preprint #1, #2 ***}
}

%% file: option_keys
\ifx\optionkeymacros\undefined\else \fi

\catcode`\Œ=\active\defŒ{{\aa}}       
\catcode`\º=\active\defº{\int}        
\catcode`\=\active\def{\c c}        
\catcode`\¶=\active\def¶{\partial}    
\catcode`\Ä=\active\defÄ{\oint}       
\catcode`\Æ=\active\defÆ{\triangle}   
\catcode`\Â=\active\defÂ{\neg}        
\catcode`\µ=\active\defµ{\mu}         
\catcode`\¿=\active\def¿{{\o}}        
\catcode`\¹=\active\def¹{\pi}         
\catcode`\Ï=\active\defÏ{{\oe}}       
\catcode`\§=\active\def§{{\ss}}       
\catcode`\ =\active\def {\dagger}     
\catcode`\Ã=\active\defÃ{\sqrt}       
\catcode`\·=\active\def·{\Sigma}      
\catcode`\Å=\active\defÅ{\approx}     
\catcode`\½=\active\def½{\Omega}      
\catcode`\£=\active\def£{{\it\$}}     
\catcode`\°=\active\def°{\infty}      
\catcode`\¤=\active\def¤{{\S}}        
\catcode`\¦=\active\def¦{{\P}}        
\catcode`\¥=\active\def¥{\bullet}     
\catcode`\»=\active\def»{\leavevmode\raise.585ex\hbox{\b a}}      
\catcode`\¼=\active\def¼{\leavevmode\raise.6ex\hbox{\b o}}        
\catcode`\­=\active\def­{\not=}       
\catcode`\²=\active\def²{\leq}        
\catcode`\³=\active\def³{\geq}        
\catcode`\Ö=\active\defÖ{\div}        
\catcode`\É=\active\defÉ{{\dots}}     
\catcode`\¾=\active\def¾{{\ae}}       
\catcode`\Ç=\active\defÇ{\ll}         
\catcode`\Ò=\active\defÒ{``}          
\catcode`\Á=\active\defÁ{!`}          
\catcode`\¢=\active\def¢{\rlap/c}     
\catcode`\Ô=\active\defÔ{`}           
\catcode`\Õ=\active\defÕ{'}           


\catcode`\=\active\def{{\AA}}       
\catcode`\'=\active\def'{\c C}        
\catcode`\¯=\active\def¯{{\O}}        
\catcode`\¸=\active\def¸{\Pi}         
\catcode`\Î=\active\defÎ{{\OE}}       
\catcode`\®=\active\def®{{\AE}}       
\catcode`\×=\active\def×{\diamond}    
\catcode`\¡=\active\def¡{\accent'27}  
\catcode`\Ó=\active\defÓ{''}          
\catcode`\±=\active\def±{\pm}         
\catcode`\È=\active\defÈ{\gg}         
\catcode`\À=\active\defÀ{?`}          
\catcode`\Ð=\active\defÐ{--}          
\catcode`\Ñ=\active\defÑ{---}         


\catcode`\Š=\active\defŠ{\"a}        
\catcode`\'=\active\def'{\"e}        
\catcode`\•=\active\def•{\"{\i}}     
\catcode`\š=\active\defš{\"o}        
\catcode`\Ÿ=\active\defŸ{\"u}        
\catcode`\Ø=\active\defØ{\"y}        
\catcode`\€=\active\def€{\"A}        
\catcode`\…=\active\def…{\"O}        
\catcode`\†=\active\def†{\"U}        
\catcode`\‡=\active\def‡{\'a}        
\catcode`\Ž=\active\defŽ{\'e}        
\catcode`\'=\active\def'{\'{\i}}     
\catcode`\—=\active\def—{\'o}        
\catcode`\œ=\active\defœ{\'u}        
\catcode`\ƒ=\active\defƒ{\'E}        
\catcode`\ˆ=\active\defˆ{\`a}        
\catcode`\=\active\def{\`e}        
\catcode`\"=\active\def"{\`{\i}}     
\catcode`\˜=\active\def˜{\`o}        
\catcode`\=\active\def{\`u}        
\catcode`\Ë=\active\defË{\`A}        
\catcode`\‹=\active\def‹{\~a}        
\catcode`\–=\active\def–{\~n}        
\catcode`\›=\active\def›{\~o}        
\catcode`\Ì=\active\defÌ{\~A}        
\catcode`\"=\active\def"{\~N}        
\catcode`\Í=\active\defÍ{\~O}        
\catcode`\‰=\active\def‰{\^a}        
\catcode`\=\active\def{\^e}        
\catcode`\"=\active\def"{\^{\i}}     
\catcode`\™=\active\def™{\^o}        
\catcode`\ž=\active\defž{\^u}        

\let\optionkeymacros\null

%% file: Veselov.introduction
\centerline{\bf A Compactification of H\'enon Mappings in $\bold C^2$}
\centerline{\bf as Dynamical Systems}
\bigskip
\centerline{John Hubbard, Peter Papadopol and Vladimir Veselov}
\bigskip
\heading
I. Introduction
\endheading
\bigskip

Let $H:\C^2 \to \C^2$ be a H\'enon mapping:
$$
 H\Bigl(\bvec xy\Bigr) = \bvec {p(x)-ay}x,\ a \ne 0
$$
where $p$ is a polynomial of degree $d \ge 2$, which without loss of generality we may
take to be monic. 

In \cite {HO1}, it was shown that there is a topology on $\C^2 \sqcup S^3$
homeomorphic to a 4-ball such that the H\'enon mapping extends continuously. That paper
used a delicate analysis of some asymptotic expansions, for instance, to understand the
structure of forward images of lines near infinity. The computations were quite difficult,
and it is not clear how to generalize them to other rational maps.

In this paper we will present an alternative approach, involving blow-ups rather than
asymptotics.  We apply it here only to H\'enon mappings and their compositions, but the
method should work quite generally, and help to understand the dynamics of rational maps
$f:\Proj^2
\ratto
\Proj^2$ with points of indeterminacy. The application to compositions of H\'enon maps
proves a result suggested by Milnor, involving embeddings of solenoids in
$S^3$ which are topologically different from those obtained from H\'enon mappings.
In the papers \cite {Ves}, \cite{HHV}, the method is applied to some other
families of rational maps. 

The approach consists of three steps, which we describe below.

{\bf --Resolving points of indeterminacy}  The general theory asserts that  a
``rational map''
$f:\Proj^2 \ratto \Proj^2$ is defined except at finitely many points, and that after a
finite number of blow-ups at these points, the map becomes well-defined [\Sha, IV.3, Thm.
3].  Let us denote by $\tilde X_f$ the space obtained after these blow-ups, and $\tilde
f:\tilde X_f \to \Proj^2$ the lifted morphism. For H\'enon mappings, this is done in
Section III. Section II is an introduction to blow-ups. Every book on algebraic geometry
(\cite{Har} for instance) has a treatment of this subject, but we are aiming the paper at
people who work in Dynamical Systems, and who will often find this literature
inpenetrable; in any case the theory is usually stated in terms of the functor {\bf Proj}
and we find that computations are quite difficult in that language. So we have included
the definitions and basic properties, as well as some fundamental examples.

{\bf --The complex sequence space}  The mapping $\tilde f$ cannot be considered a
dynamical system, since the domain and range are different.  One way to obtain a
dynamical system is to blow up the inverse images of the points you just blew up to
construct $\tilde X_f$, then their inverse images, etc. On the the projective limit, we
finally obtain a dynamical system $f_\infty:X_\infty \to X_\infty$.

The notation to keep track of successive blow-ups grows exponentially and soon becomes
intractable. There is an alternative description of the projective limit in terms of
sequence spaces, which is much easier to describe.  We learned of this construction from
\cite {F1}; something analogous was constructed by Hirzebruch \cite {Hirz} when
resolving the cusps of Hilbert modular surfaces, and also considered by Inoue, Dloussky
and Oeljeklaus \cite{Dl}, \cite{DO}.

The space $X_\infty$ will usually have
a big subset $X^*_\infty$ which is an algebraic variety, but some points which are
extremely singular, in the sense that every neighborhood has infinite dimensional
homology.

We will construct $X_\infty$ for H\'enon mappings in Section IV, and its topology
(homology, etc.) is studied in Section V.  There is a natural way to complete the
homology $H_2(X^*_\infty, \C)$ to a Hilbert space, so that $f_\infty$ induces an
operator, and we hope in a future publication to show that the spectral invariants of this
operator interacts with the dynamics.

{--\bf The real oriented blow-up}   In order to ``resolve'' the terrible singularities of
$X_\infty$, we will use a further {\it real\/} blow-up, in which we consider complex
surfaces as 4-dimensional real algebraic spaces, and divisors as real surfaces in them.
One way of thinking of this blow-up of a surface $X$ along  a divisor $D$  is to take an
open tubular neighborhood $W$ of $D$ in $X$, together with a projection
$\pi:\overline W \to D$. Excise $W$, to form
$X'=X-W$.  If $W$ is chosen properly, $X'$ will be a real 4-dimensional
manifold with boundary $\partial X' =\partial \overline W$.  The interior of $X'$ will be
homeomorphic to $X-D$. Then $X'$ is some sort of blow-up of $X$ along $D$, with
$\pi: \partial X' \to D$ the exceptional divisor. This construction is topologically
correct, but non-canonical; the real oriented blow-up is a way of making it canonical.

We can pass to the projective limit with these real oriented blow-ups, constructing a
space
$\orb(X_\infty)$, which is topologically much simpler and better behaved than the original
compactification, being a 4-dimensional manifold with boundary. The real interest is in
the inner structures of the boundary 3-manifold, where we find solenoids (in this paper
and in
\cite{Ves}), tori with irrational foliations (in \cite{HHV}), etc.  

The definition and first properties of these blow-ups are given in Section VI, and the
methods needed to construct them are given in Section VII. Theorems \blowupofordpoint\ and
\blowupdoublepoint\  are the principal tools for understanding real oriented blow-ups,
and we expect that they will be useful for many examples besides the ones explored in
this paper.  

In Section VIII, we show that the classical Hopf fibration is an example of a
real oriented blow-up, in two different ways. In Section IX, we construct the real
oriented blow-up for H\'enon maps. This turns out to be quite an exciting space, and we
further explore its structure in Section X, using toroidal decompositions.  These results
allow us to prove (Theorem \conjugacyofexthens) that extensions of the H\'enon maps to
their sphere at infinity are not all conjugate, even when they have the same degree; the
conjugacy class of the extension remembers the argument of the Jacobian of the H\'enon
mapping.

In Section XI, we construct the real-oriented blow-up for compositions of H\'enon maps,
obtaining a 3- sphere with two embedded solenoids $\Sigma^+$ and $\Sigma^-$, but such
that the incompressible tori in $S^3-(\Sigma^+\cup \Sigma^-)$ are different from the
incompressible tori we obtain from just one H\'enon mapping.  This difference was first
conjectured by Milnor
\cite{Mi2}.

We wish to thank A. Douady for essential help in the conception of this paper, J. Milnor
for his conjecture, Profs. Soublin (Marseille), A. Sausse
(Sophia-Antipolis), Mark Gross (Cornell) and Mike Stillman (Cornell) for help with the
algebra, Alan Hatcher for help with the topology of Seifert fibrations, and Sa'ar
Hersonsky for suggestions to improve the introduction.

Moreover, P. Papadopol thanks President Bill Williams for his support over many years, and
Cornell University for many years of hospitality.

\subheading{Outline of the paper}

\halign {\quad #&\hfil#:\ &\ #.\hfil\cr  
Section&I& Introduction\cr
Section&II& Blowing up a closed subspace of an analytic space\cr
Section&III& Making H\'enon mappings well-defined\cr
Section&IV& Closures of graphs and sequence spaces\cr
Section&V& The homology of $X_\infty^*$\cr
Section&VI& Real blow-ups of complex divisors in surfaces\cr
Section&VII& The effect of complex blow-ups on real oriented blow-ups\cr
Section&VIII&  Real oriented blow-ups and the Hopf fibration\cr
Section&IX&  Real blow-ups for H\'enon mappings\cr
Section&X& The topology of $\orb(X_\infty,D_\infty)$\cr 
Section&XI& The compactification of compositions of H\'enon mappings\cr}

%% file: universalblowups
\heading
II. Blowing up a closed subspace of an analytic space
\endheading
\bigskip

\subheading{Hypersurfaces and divisors}

Much of this paper will be about making appropriate blow-ups.  In
Section III we will simply be blowing up surfaces at points, but in
Section VI we will need to make much more elaborate blow-ups, of
singular surfaces in 4-dimensional manifolds.  Everything in this
section is presumably standard, but we found it difficult to carry out
our computations in the too elegant language of {\bf Proj}, and we find
the present treatment, in terms of affine equations, better adapted to
our needs. The reader who feels comfortable with the simplest blow-ups,
of surfaces at points, should go directly to the next section, and
refer back when needed (which may well occur in Section VI). It seems
that in the generality in which we will be describing the construction,
the results of this section are due to Hironaka [\Hir].

Let $X$ be an analytic space with structure sheaf $\Cal O_X$,
and $Y\subset X$ a closed subspace. The definition of the blow-up $\tilde
X_Y$ of $X$ along $Y$ requires carefully distinguishing between
hypersurfaces and divisors.

\definition {\hyperanddiv} A hypersurface in $X$ is a subspace
defined by an ideal which locally requires only one generator;
i.e., it is the locus locally defined by a single equation.

\definition {\divisordef} A divisor is a hypersurface locally defined by
a single equation which is not a zero divisor. More precisely, a subspace
$Z \subset X$ is a divisor if every point $x \in X$ has a neighborhood
$U$ such that $Z\cap U$ is defined by the principal ideal generated by
some section $f \in \Cal O_X(U)$ which is not a zero-divisor in the ring
$\Cal O_X(U)$. \medskip

Consider for instance $X = \C \times \{0\} \cup \{0\}\times \C\subset
\C^2$, or alternately, the subspace of $\C^2$ defined by the equation
$xy=0$.  The $x$-axis alone is a hypersurface in $X$, defined by
the equation $y=0$, but it is not a divisor, since $y$ is a zero
divisor in $\Cal O_X(X) = \C[x,y]/(xy)$. This example is typical:
geometrically, saying that a function on a space is a zero divisor is 
saying that the space has several components, and that the function
vanishes identically on some of them. 

\subheading{Definition of the blow-up}

Now we can define $\tilde X_Y$ as the ``universal way'' of
replacing $Y$ by a divisor.

\definition {\blowupdef} The blow-up $\tilde X_Y$ of an analytic space
$X$ along a closed subspace $Y$ is an analytic space $\tilde X_Y$, which
comes with a morphism $\pi:\tilde X_Y \to X$ (called the canonical
projection), such that
$\pi^{-1}(Y)$ is a divisor (called the exceptional divisor) in $\tilde
X_Y$, and if
$g:Z \to X$ is a morphism such that $g^{-1}(Y)$ is a divisor in $Z$,
there exists a unique morphism $\tilde g:Z \to \tilde X_Y$ such that
$\pi \circ \tilde g = g$.

We will prove that the blow-up exists (and is then obviously
unique up to unique isomorphism).
If we can construct the blow-up locally, and prove that it has
the right universal property, then the existence of a global
blow-up will follow, since the local blow-ups will patch
uniquely.

\subheading{Local construction of a blow-up}

So we may assume that $Y \subset X$ is defined by the equation
$Y= f^{-1}(0)$, where
 $f:X \to
\C^{n+1}$ is an analytic mapping (i.e., a section of $\Cal
O_X^{n+1}$).

We define the blow-up $\tilde X_Y$ as follows. Recall that $\Proj^n$
is the space of lines $l \subset \C^{n+1}$ through the origin, and
that it carries the tautological line bundle 
$$
\matrix E&&\hookrightarrow&&\C^{n+1}\times \Proj^n\\
&\SEarrow&&\SWarrow&\\
&&\Proj^n&&\\
\endmatrix
$$
where $E=\set{(x,l) \in \C^{n+1} \times \Proj^n}{x \in l}$.

First consider the locus $X'_Y \subset X \times \Proj^n$ defined as
$$ 
X'_Y = \set {(x,l) \in X \times \Proj^n}{f(x) \in l}.
$$
This sounds set-theoretic rather than ideal theoretic, but we can
easily define $X'_Y$ by equations. Let $[U_0: \dots:
U_n]$ be homogeneous coordinates on $\Proj^n$. Then $X'_Y$ is defined by
the ideal (a mixed ideal, homogeneous with respect to the $U_i$, affine
with respect to the $f_j$), generated by all 
$$  
U_i f_j- U_j f_i \peqno \detideal
$$  
for $i\ne j$.

Let $\pi:X'_Y \to X$ be induced by the projection
onto the first factor. We now will see that locally $\pi^{-1}(Y)$ is
locally defined by a single equation.

Let $E \to \Proj^n$ be the tautological line bundle, and $F \to X'_Y$ be
the pull-back of $E$ by the composition
$$
X'_Y \to X \times \Proj^n \to \Proj^n.
$$
This leads to the following commutative diagram:
$$
\matrix &&X\times \C^{n+1}\times \Proj^n &&\C^{n+1}\times \Proj^n\\
&&\cup&&\cup\\
F&\hookrightarrow& X\times E& \to &E\\
\downarrow&&\downarrow&&\downarrow\\
X'_Y& \hookrightarrow& X\times \Proj^n& \to& \Proj^n
\endmatrix
$$
where both squares are fiber products.

Therefore the line bundle $F$ is a subbundle of the trivial bundle
$X'_Y \times \C^{n+1} \to X'_Y$, and the map $X'_Y \to \C^{n+1}$ induced by
$f$ is a section $f'$ of $F$, so that locally on $X'_Y$ the set $Y'$
defined  by $f'=0$ is a hypersurface.

More precisely, let us denote by $f'':X'_Y \to \C^{n+1}$ the
composition
$$
X'_Y \to X \times\Proj^n \to X \to \C^{n+1}.
$$
In the chart $U_i \ne 0$ on
$\Proj^n$, with affine coordinates $u_j = U_j/U_i,\ j\ne i$, a
non-vanishing section of the tautological bundle is given by   
$$
\sigma_i:\bmatrix u_0\\ \vdots \\u_{i-1}\\u_{i+1} \\ \vdots \\
u_n\endbmatrix   \mapsto \bmatrix u_0\\ \vdots \\u_{i-1}\\1\\u_{i+1} \\
\vdots \\ u_n\endbmatrix
 $$
with the $1$ in the $i$th position.  This section lifts and
restricts to a section $\tau_i$ of $F$ above $X'_Y \cap \{U_i
\ne 0\}$. Evidently, $f'= \tau_i f''_i$ in this chart, where $f_i''$ is
the $i$th coordinate of $f''$, so
$Y'$ is given by the single equation $f''_i=0$ in this chart.

However, in even moderately complicated cases, the space $X'_Y$ has
parasitic components on which $f'$ vanishes identically, so that $f'$ is
a zero divisor, and we must get rid of these components.  

Let $\Cal A_{f'} \subset \Cal O_{X'_Y}$ be the sheaf of ideals generated
by functions which are annihilated by some power of
$f'$, and let $\tilde X_Y\subset X'_Y$ be the subspace defined
by $\Cal A_{f'}$.  This annihilator is the ideal of functions which
vanish identically on the components of $X'_Y$ on which no power
of $f'$ vanishes identically, so $\tilde X_Y$ is obtained from
$X'_Y$ by removing the ``parasitic components'' on which some power
of $f'$ vanishes identically. In the best cases, $\Cal A_{f'}$ will be
the zero sheaf of ideals, so $\tilde X_Y = X'_Y$, but in other cases,
$\tilde X_Y$ will be exactly the union of the components of $X'_Y$ on
which no power of $f'$ vanishes identically.  Thus the restriction
 $\tilde f$ of $f'$ to $\tilde X_Y$ is not a zero divisor,
and the subspace
$\tilde Y \subset \tilde X_Y $ defined by $\tilde f$
is a divisor.

We will prove in \universalproperty\ that the space $\tilde X_Y$ is the
blow-up of $X$ along $Y$, and that $\tilde Y$ is the exceptional
divisor. First we will give some examples.

\subheading {Blowing up a point in a surface}

The easiest example of a blow-up, and also the only one we will use
until Section VI, is the blow-up of a surface at a point.  Because of
the universal property, it is enough to understand the blow-up in one
chart, i.e., to understand the blow-up of $\C^2$ at the origin.

\example \blowuppointonsurf Let $x,y$ be the coordinates of $\C^2$;
the origin is of course defined by the equations
$x=y=0$, so the blow-up is contained in $\C^2\times \Proj^1$.  If we
use
$\bvec {U_0}{U_1}$ as homogeneous coordinates in $\Proj^1$, then the
equation expressing that the point $\bvec xy$ is on the line
corresponding to
$\bvec {U_0}{U_1}$ is
$xU_1=yU_0$.  This equation is not a zero-divisor, so we have computed
the blow-up.

This is covered by two affine coordinate charts:

--One in which
$U_1\ne 0$ and $u_0=U_0/U_1$; in this chart the blow-up is given by the
equation $x=yu_0$; clearly it is a non-singular surface parametrized by
$y$ and $u_0$.  The exceptional divisor is given in this chart by the
single equation $y=0$.

--One in which
$U_0\ne 0$ and $u_1=U_1/U_0$; in this chart the blow-up is given by the
equation $y=xu_1$; clearly it is a non-singular surface parametrized by
$x$ and $u_1$. The exceptional divisor is given in this chart by the
single equation $x=0$.

So in this case the effect of blowing up is to replace the origin by
the projective line $\Proj^1_\C$. More generally, a similar computation
will show that if you blow up a smooth manifold $X$ along a smooth
submanifold $Y \subset X$, then the fibers $\pi^{-1}(x)$ of the
projection
$\pi:\tilde X_Y \to X$ are points if $x \notin Y$, and the projective
space $\Proj(T_xX/T_x Y)$ associated to the ``normal space''  
$T_xX/T_xY$ if $x \in Y$.

\subheading {Two more complicated examples}

\example \therealcase Let $X = \C^2$, and $Y$ be the subset
defined by the equations $\{x^2=0,y=0\}$.  Clearly in this case
$n=1$, and the space $X'_Y$ is the subspace of $X \times
\Proj^1$, given by the equation
$$
x^2U_1=yU_0,
$$
where $[U_0:U_1]$ are homogeneous coordinates on $\Proj^1$.
 
In the chart where $U_0 \ne 0$, if you set $u_1=U_1/U_0$,  the space $
X'_Y$ is defined by the equation $x^2u_1=y$, and
$Y'$ is defined by the single equation $x^2=0$.  In the chart
where $U_1 \ne 0$ and $u_0=U_0/U_1$ is an affine coordinate, $X'_Y$
is the space defined by the equation  $yu_0=x^2$, and $Y'$ is defined by
the single equation $y=0$.  Since neither equation is a zero-divisor, we
see that $\tilde X_Y = X'_Y$, and $\tilde Y = Y'$.  In particular, the
blow-up is a cone, with a singular point at the vertex $x=y=u_0=0$, and
$\tilde Y$ is a ``double line'' on that cone.

Our next example is a special case of the following general fact.  Let
$X$ be an analytic space, $Y\subset X$ be a closed analytic subspace
defined by an ideal $\Cal I_Y$, and denote by $pY$ the subspace defined
by $\Cal I^p$.  Then $\tilde X_{pY}= \tilde X_Y$, and the exceptional
divisor in $\tilde X_{pY}$ is $p\tilde Y$, i.e., if $\tilde Y$ is
locally defined by the single equation $f=0$, then the exceptional
divisor of $\tilde X_{pY}$ is locally defined by $f^p=0$.  This fact is
easy to check from the universal property, but considerably harder to
see from the explicit computation, even in the case when $X=\C^2$ and
$Y$ is the origin.

\example \zerodivs  Let $X = \C^2$, $\Cal J \subset \Cal O_X$ be
the ideal generated by $\{x,y\}$, consider the subset $Y$ defined by the
ideal $\Cal J^2$, or alternately $Y=f^{-1}(0)$, where 
$$
f:\C^2 \to \C^3 \quad \text{is given by}\quad f\left(\bvec xy\right) =
\bmatrix x^2\\xy\\y^2\endbmatrix.
$$
This
time we required three equations, so that $X'_Y \subset X\times \Proj^2$ is
given (in homogeneous coordinates $[U_0:U_1:U_2]$ on $\Proj^2$) by the
equations $x^2U_1=xyU_0, x^2U_2=y^2U_0$ and $xyU_2 = y^2U_1$. 

In the chart $U_0 \ne 0$, using the coordinates $u_1,u_2$ as
above, $X'_Y$ is defined by the equations
$$
x^2u_1=xy\quad,\quad x^2u_2=y^2\quad,\quad xyu_2=y^2u_1,
$$
and it is easy to see that the third equation is a consequence
of the first two.

Further, $Y'$ is defined in this chart by the single equation $x^2=0$.
In this case, $x^2$ is a zero-divisor in $\Cal O_{X'_Y}$, since
$x^2(u_1y-u_2x)=0$. 

Further computations, in which we were helped by J.-P. Soublin (by
hand) and A. Sausse (using REDUCE), shows that $X'_Y$ has two
primary components, one corresponding to the ideal $Q$ generated by the
seven polynomials
$$
x^2U_1-xyU_0\ ,\ x^2U_2-y^2U_0\ ,\ xyU_2 - y^2U_1\ , \ x^3\ ,\ x^2y\ ,\
xy^2\ ,\ y^3\,
$$
with radical $(x,y)$, 
and the other  to the ideal $P$ generated by
$$
U_1^2-U_0U_2\ ,\ U_0y-U_1x\ ,\ U_2x-U_1y,
$$
which is its own radical.  
The annihilator of $f'$ is exactly $P$. So $\tilde X_Y$ is the locus in
$\C^2\times \Proj^2$ defined by the mixed ideal $P$. This is isomorphic
to the blow-up of $\C^2$ at the origin defined by $\Cal J$. Indeed,
consider the mapping $X\times \Proj^1 \to X \times \Proj^2$ given in
homogeneous coordinates by
$$
(x,y,V_0,V_1) \mapsto (x,y, V_0^2, V_0V_1, V_1^2).
$$
The equation $U_0U_2=U_1^2$ defines the image of this embedding, and is
one of the equations of $\tilde X_Y$, so the inverse image of $\tilde
X_Y$ is defined by the mixed ideal generated by
$$
V_0^2y-V_0V_1x\quad,\quad V_1^2x-V_0V_1y,
$$
which is the same as the principal ideal generated by $V_0y-V_1x$.  But
this ideal also defines the blow-up of $\C^2$ at the origin defined by
the ideal $\Cal J$. The proper transform of $Y$, i.e., the locus defined
by
$\tilde f$, is the exceptional divisor $x=y=0$ as a double line.

\subheading{ Proof of the universal property}

We must now prove that our blow-up $\tilde X_Y$ has the right universal property.

\theorem \universalproperty If $g:Z \to X$ is a
mapping such that $g^{-1}(Y)$ is a divisor in $Z$,
then there exists a unique morphism $\tilde g:Z \to
\tilde X_Y$ such that $\pi \circ \tilde g = g$.
\endstat

\proof  By definition, $g^{-1}(Y)$ is the locus
defined by the equation $f\circ g = 0$.  On the other
hand, locally on $Z$, $g^{-1}(Y)$ is defined by a
single equation $h=0$.  More precisely, any point $z_0
\in Z$ has a neighborhood $V \subset Z$ in which $g^{-1}(Y)
\cap V$ is defined by the ideal $h\Cal O_V$, where $h$ is
a function, not a zero-divisor. Since each 
$g\circ f_i$ is in the ideal generated by $h$, we can find
functions
$\phi_i$ such that $ g\circ f_i= \phi_i h$. Moreover,
the $\phi_i$ do not all vanish at any point of $V$
since they generate the unit ideal. So we can find a
morphism $\tilde g_V:V \to X \times \Proj^n$ by the
formula
$$
\tilde g_V(z) = (g(z), [\phi_0(z):\dots:\phi_n(z)]).
$$

First, observe that this morphism does not depend on
the choice of $V$ and $h$.  If $(V',h')$ and
$(V'',h'')$ are two open sets in which $g^{-1}(Y)$
are defined by a single function, then on $V' \cap
V''$, we have that $h'=h'' u$, where $u$ is
invertible (i.e., a unit in $\Cal O(V' \cap
V'')$). With the obvious notation, 
$$
 (\phi'_0,\dots, \phi'_n) = u (\phi''_0,\dots,
\phi''_n),
$$
so the corresponding homogeneous coordinates define
the same mapping into $X \times \Proj^n$.

Thus we have a morphism $\tilde g:Z \to X \times
\Proj^n$; we must show that the image is contained
in $\tilde X_Y$.

First, it is clear that the image lies in $X'_Y$, or
more precisely, that       
$$ 
\tilde g^*(a_i f_j- a_j f_i)= \phi_i(
f_j\circ g)-\phi_j( f_i \circ g) = h((f_i\circ g)(f_j\circ g)-(
f_j\circ g)(f_i\circ g) )=0.
$$

Next, if $\alpha \in \Cal A_{f'}$, we need to know that $\tilde
g^*\alpha = 0$. Since $\alpha (f')^m=0$ for some $m$, we have
$$
(\alpha \circ \tilde g)(f'\circ \tilde g)^m=0,
$$
which certainly implies that $f'\circ \tilde g$ is a
zero-divisor if $\alpha \circ \tilde g \ne 0$. But $f'\circ
\tilde g$ defines $f^{-1}(Y)$, which was assumed to be a
divisor, so defined by a single equation, not a zero-divisor.
\qed

\corollary \natblowups If $f:X_1 \to X_2$ is a morphism of analytic
spaces, $Y_2 \subset X_2$ is an analytic subspace, and
$Y_1=f^{-1}Y_2$, then there exists a unique morphism $\tilde f:\tilde
{(X_2)}_{Y_2} \to \tilde{(X_1)}_{Y_1}$ which makes the diagram
$$
\CD
 \tilde {(X_2)}_{Y_2}  @>{\tilde f}>> \tilde{(X_1)}_{Y_1} \\ 
@VVV   @VVV     \\ 
X_2     @>>>  X_1 
\endCD
$$
commute.
\endstat

\proof Clearly under the composition $\tilde {(X_2)}_{Y_2} \to X_2 \to
X_1$, the inverse image of $Y_1$ is the exceptional divisor, hence we
can apply Theorem \universalproperty.
\qed

We only defined $\tilde X_Y$ when $Y$ is defined by global
equations.  However, the universal mapping properties guarantees
that if we perform blow-ups locally, they will glue together in a
unique fashion, so that this actually constructs a global
blow-up of an arbitrary analytic space along a closed
subspace.

\subheading{An example in $\C^4$}

Let us denote by $x_1,x_2,y_1,y_2$ the coordinates of $\C^4$. In Section
6, we will need to understand the blow-up $\C^4$ along the union 
$$
Y= \C^2 \times \{0\} \cup \{0\} \times \C^2
$$
of the $(x_1,x_2)$ and the $(y_1,y_2)$ coordinate planes.  This also
provides an example where many of the complications of the previous
sections occur and thus illustrates what all these zero-divisors and
annihilators really mean. 

Although $Y$ is of codimension 2, it
requires four equations for its definition:
$$
  x_1y_1=0, x_2y_1=0, x_1y_2=0 \text{  and  } x_2y_2=0,
$$
i.e., $Y=f^{-1}(0)$ where $f:\C^4 \to \C^4$ is given by
$$
f\left(\bmatrix x_1\\x_2\\y_1\\y_2\endbmatrix \right) = \bmatrix
x_1y_1\\x_2y_1\\x_1y_2\\x_2y_2\endbmatrix.
$$
   The fact that it is not
defined by two equations will follow from the fact that the blow-up
cannot be embedded in $X \times \PP^1$, and the fact that four equations
are really required will follow when we see that the blow-up cannot be
embedded in $X \times \PP^2$ either.

Since the ideal defining the union of the coordinate planes is generated
by  $$ 
x_1y_1, x_2y_1, x_1y_2, x_2 y_2,
$$  
the variety $X'_Y\subset \C^4 \times \Proj^3$, is defined by the equation
$f(x) \in \ell$, and if we use homogeneous coordinates $U_1,U_2, U_3,
U_4$ on the second factor, this means that the vectors  
$$  
\bmatrix x_1 y_1\\x_2y_1\\ x_1y_2 \\ x_2y_2\endbmatrix\quad
\text{and}\quad \bmatrix U_1\\ U_2\\ U_3\\U_4  \endbmatrix
$$  
are linearly dependent. Although locally this locus is defined by three
equations,  globally it requires 6 equations: 
$$
\eqalign{ x_1y_1U_4&= x_2y_2U_1 \cr x_2y_1U_4&= x_2y_2U_2 \cr
x_1y_2U_4&= x_2y_2U_3 \cr x_1y_1U_3&= x_1y_2U_1 \cr x_2y_1U_3&=
x_1y_2U_2 \cr x_1y_1U_2&= x_2y_1U_1.} \peqno \eqbeforezerodiv
$$

Let us denote by $\pi:X'_Y \to \C^4$ the projection induced by the
projection $X \times \Proj^3 \to X$ onto the first factor. Locally,
the locus $Y'= \pi^{-1}(Y)$ is defined by a single equation, as it
should.  For instance, in the chart $U_4=1$, the space $X'_Y$ is defined
in $X \times \Proj^3$ by the 3 equations 
$$
\eqalign{ x_1y_1 & = x_2y_2u_1 \cr
 x_2y_1&= x_2y_2u_2 \cr
 x_1y_2&= x_2y_2u_3},
$$ 
where $u_i=U_i/U_4, i=1,2,3$, and $Y'$ is defined by the single equation
$x_2y_2 = 0$.  However, the functions  $$
x_2y_2, x_2y_1, x_1y_1, x_1y_2
$$
 are zero-divisors in $\Cal O(X'_Y)$: for instance,
$$ 
x_2y_2(U_4y_1-U_2y_2)=0
$$ 
in the ring of functions on $X'_Y$.

A careful analysis of equations (\eqbeforezerodiv) shows that $X'_Y$ is
the union of five irreducible components, the four  4-dimensional linear
spaces $Z_1,
\cdots, Z_4$ with equations $$
\eqalign{ x_1=x_2=y_1&=0 \cr x_1=y_1=y_2&=0 \cr x_1=x_2=y_2&=0 \cr x_2=y_1=y_2&=0 }
$$ 
and the 4-dimensional space $Z_5$ with equation
$$
\eqalign{ U_1U_4 &= U_2U_3\cr 
y_1U_4 &= y_2U_2\cr 
x_1U_4 &= x_2U_3\cr 
y_1U_3 &= y_2U_1\cr
x_1U_2 &= x_2U_1}\peqno \eqafterzerodiv 
$$

We can now see that $\tilde X_Y$ is in fact $Z_5$.  If we embed $\PP^1
\times \PP^1$ into $\Proj^3$ by the Veronese mapping
$$
\left(\bvec{X_1}{X_2},\bvec{Y_1}{Y_2}\right)
\mapsto  
\bmatrix X_1Y_1 \\ X_2Y_1 \\ X_1Y_2 \\ X_2Y_2 \endbmatrix,
$$ 
where $\bvec{X_1}{X_2},\bvec{Y_1}{Y_2}$ are homogeneous coordinates in
$\Proj^1
\times \Proj^1$, we observe that the image has equation $U_1U_4 =
U_2U_3$, which is satisfied identically on $\tilde X_Y$. Thus the
blow-up is actually contained in $X \times \left(\PP^1 \right)^2$, and it
is defined by the equations 
$$
\eqalign{ 
y_1X_2Y_2 &=y_2X_2Y_1 \cr
x_1X_2Y_2 &=x_2X_1Y_2 \cr 
y_1X_1Y_2 &=y_2X_1Y_1 \cr 
x_1X_2Y_1 &=x_2X_1Y_1. }\peqno \eqviaVeronese 
$$

From the first equation we can see that either $X_2=0$ or $Y_2y_1=Y_1y_2$.  But if $X_2=0$, then
$X_1 \ne 0$, so we can cancel $X_1$ in the third equation and get the same equation
$Y_2y_1=Y_1y_2$. Therefore the first and third equation are equivalent to the equation
$Y_2y_1=Y_1y_2$.  We can apply the same procedure to the second and fourth equations and see that
they are equivalent to
$X_2x_1=X_1x_2$. Thus equations (\eqviaVeronese) are equivalent to
equations 
$$
\eqalign{ 
y_1Y_2&=y_2Y_1  \cr 
x_1X_2&=x_2X_1}\peqno \blowupstillsingeq
$$ 
The last two equations show that $\tilde X_Y$ is a manifold.  For
instance, in the chart $X_2=Y_2=1$,
$\tilde X_Y$ is given by the equations $x_1 = X_1x_2$ and $y_1=Y_1y_2$,
clearly parametrized by $x_2, y_2, X_1, Y_1$, and $\tilde Y =
\pi^{-1}(Y)$ in this chart is given by the single equation $x_2y_2=0$.
We see in particular that the exceptional divisor $\tilde Y$ is not
smooth, above $(0,0,0,0)$, it has the local structure of $W \times
\C^2$, where $W \subset \C^2$ is the singular curve with equation
$xy=0$, i.e., the union of the axes.

%% file: Veselov.welldefined.new
\heading III. Making H\'enon Mappings well-defined \endheading
 
Consider the H\'enon mapping
$$
H\pmatrix x\\ y\endpmatrix = \pmatrix {p(x)-ay}\\ x\endpmatrix \peqno \henaffine
$$
with $a \ne 0$, which we will consider as a birational mapping $\PP^2 \ratto \PP^2$ given in
homogeneous coordinates as
$$
H\bmatrix x\\y\\z\endbmatrix = \bmatrix \tilde p(x,z)-ayz^{d-1}\\xz^{d-1}\\z^d 
\endbmatrix, \peqno \henhomogeneous
$$
where $\tilde p (x,z) = z^dp(x/z) = x^d + \dots$ is a homogeneous polynomial of
degree $d$ in the two variables $x$ and $z$.

\lemma \pointsofindetlemma a) The mapping $H$ has a unique point of indeterminacy
at $\bold p=\bmatrix 0\\1\\0\endbmatrix$, and collapses the line at infinity $l_\infty$  to the
point
$\bold q= \bmatrix 1\\0\\0\endbmatrix$.

b) The mapping $H^{-1}$ has a unique point of indeterminacy at $\bold q$, and
collapses
$l_\infty$ to the point $\bold p$.
\endstat

\proof A point of indeterminacy of a mapping written in homogeneous coordinates without
common factors is a point where all coordinate functions vanish. In order for this to
happen, we must have $z=0$ of course, and the only remaining term is then $x^d$, so that at
a point of indeterminacy, we also have $x=0$.  Thus $\bold p$ is the unique point of
indeterminacy of $H$.  Clearly any other point of $l_\infty$ is mapped to $\bold q$.  Part (b) is
similar.
\qed

\subheading {Self-intersections} It is often necessary to know the
self-intersection numbers of lines obtained when you make successive blow-ups.
The rules for computing these numbers are simple, in this case where we are
blowing up a surface at a smooth point $\bold x$ [\Sha, IV.3, Thm. 2 and its
Corollaries]:

$\bullet$ The exceptional divisor has self-intersection $-1$;

$\bullet$ The proper transform $C'$ of any smooth curve $C$ passing
through $\bold x$ has its self-intersection decreased by 1.

We will now go through a  sequence of $2d-1$ blow-ups, required to
make the H\'enon mapping (\henhomogeneous) well-defined.  The results are summarized at the end in
Theorem
\HenonWellDefined, which we cannot state without the terminology which we create during the
construction. We will denote $ H_1, \dots, H_{2d-1}$ the extension of the H\'enon mapping to
the successive blow-ups. 

To focus on the point of indeterminacy, we will begin work in the coordinates
$u= x/y, v=1/y$ so that the point of indeterminacy $\bold p$ is the point $u=v=0$.
 The H\'enon mapping, written in the affine  coordinates $u,v$ in the domain and homogeneous
coordinates in the range, is written
$$
H:\pmatrix u\\v\endpmatrix \mapsto \bmatrix \tilde p(u,v)-av^{d-1}\\ uv^{d-1}\\ v^d\endbmatrix.
$$
At a point of indeterminacy all three homogeneous coordinates vanish; we already knew that this
happens only at $\bold p$, but it is clear again from this formula that it happens only at
$u=v=0$.
  
\subheading {The first blow-up}  Blow up $\Proj^2$ at the point $u=v=0$,
look in the chart $v=X_1u$; i.e., use $u$ and $X_1$ as coordinates and discard $v$.

The extension of the H\'enon mapping, written in the affine coordinates $u,X_1$ in the domain and
homogeneous coordinates in the range, is written
$$
H_1: \pmatrix u\\{X_1}\endpmatrix \mapsto \bmatrix \tilde p(u,X_1u)-a(X_1u)^{d-1}\\
u^dX_1^{d-1}\\ u^dX_1^d\endbmatrix=  \bmatrix uq(X_1)-a{X_1}^{d-1}\\ uX_1^{d-1}\\
u{X_1}^d\endbmatrix, \peqno \veryfirstblowup
$$
where we have set $\tilde p(1,X)=q(X)$, so that $q$ is a polynomial of degree $d$ whose constant
term is $1$.

 Again, the only point where all three homogeneous coordinates vanish is where
$u=X_1=0$. We invite the reader to check that the one point of the blow-up not covered by the
chart $u,X_1$ is not a point of indeterminacy. The self-intersection numbers are as indicated in
Figure \firstblowuppic: originally $l_\infty$ has self-intersection number 1, after one blow-up
its proper transform acquires self-intersection number $0$, and the exceptional divisor has
self-intersection number -1. 

\bigskip
\centerline{\BoxedEPSF{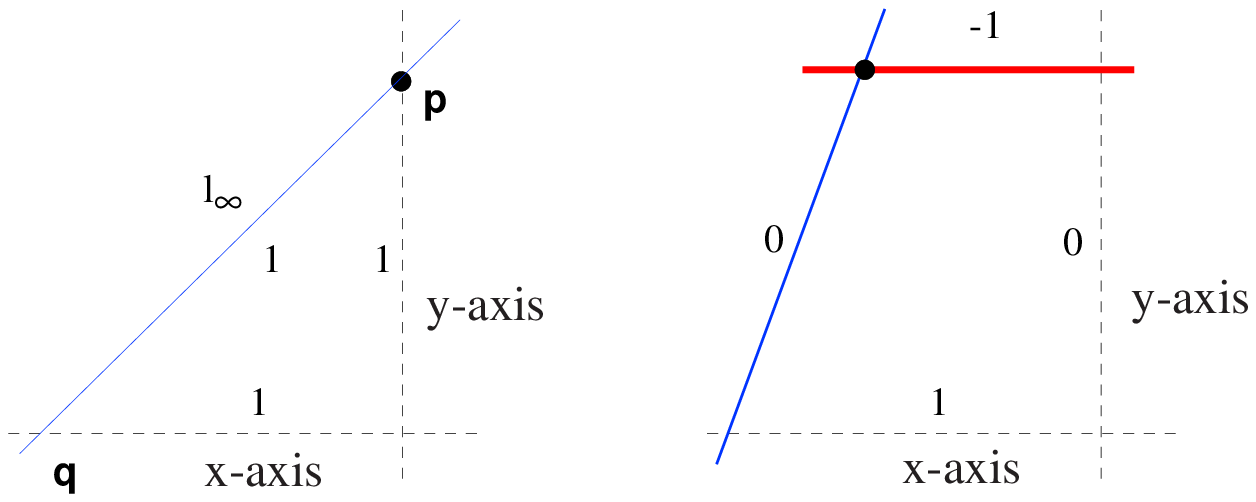 scaled 800}}
\longfigcap \firstblowuppic{The original configuration of the axes (dotted) and the line at
infinity in $\Proj^2$, and the configuration after the first blow-up. The last exceptional
divisor is denoted by a thick line; the heavy dots are the points of indeterminacy of $H$
and $H_1$. The numbers labeling components are self-intersection numbers.}

\subheading {Several more blow-ups} We will now make a sequence of $d-1$ further blow-ups,
setting successively
$$
u=X_1X_2, \ X_2=X_1X_3,\ \dots,\ X_{d-1}=X_1X_d.
$$
The H\'enon mapping in these coordinates is given by the formula
$$
H_k:\pmatrix{X_1}\\{X_k}\endpmatrix \mapsto \bmatrix
X_1X_kq(X_1)-aX_1^{d-k+1}\\X_1^dX_k\\X_1^{d+1}X_k\endbmatrix =
 \bmatrix X_kq(X_1)-aX_1^{d-k}\\X_1^{d-1}X_k\\X_1^dX_k\endbmatrix \peqno \firstkblowups
$$
and we see that when $k<d$, the unique point of indeterminacy is the point $X_1=X_k=0$, but when
$k=d$, the unique point of indeterminacy is the point $X_1=0, X_d=a$.  At each step, there is
one point $X_k=\infty$ which is not in the domain of our chart; we leave it to the reader to
check that this is never a point of indeterminacy.

Again the self-intersection numbers are as indicated in Figure \secondblowuppic. At each step we
are blowing up the intersection of the last exceptional divisor with the proper transform of the
first exceptional divisor. Therefore after this sequence of blow-ups, the last exceptional divisor
(now next-to last) acquires self-intersection -2, and the first exceptional divisor has its
self-intersection number decreased by 1, going from $-1$ to $-d$.

\centerline{\BoxedEPSF{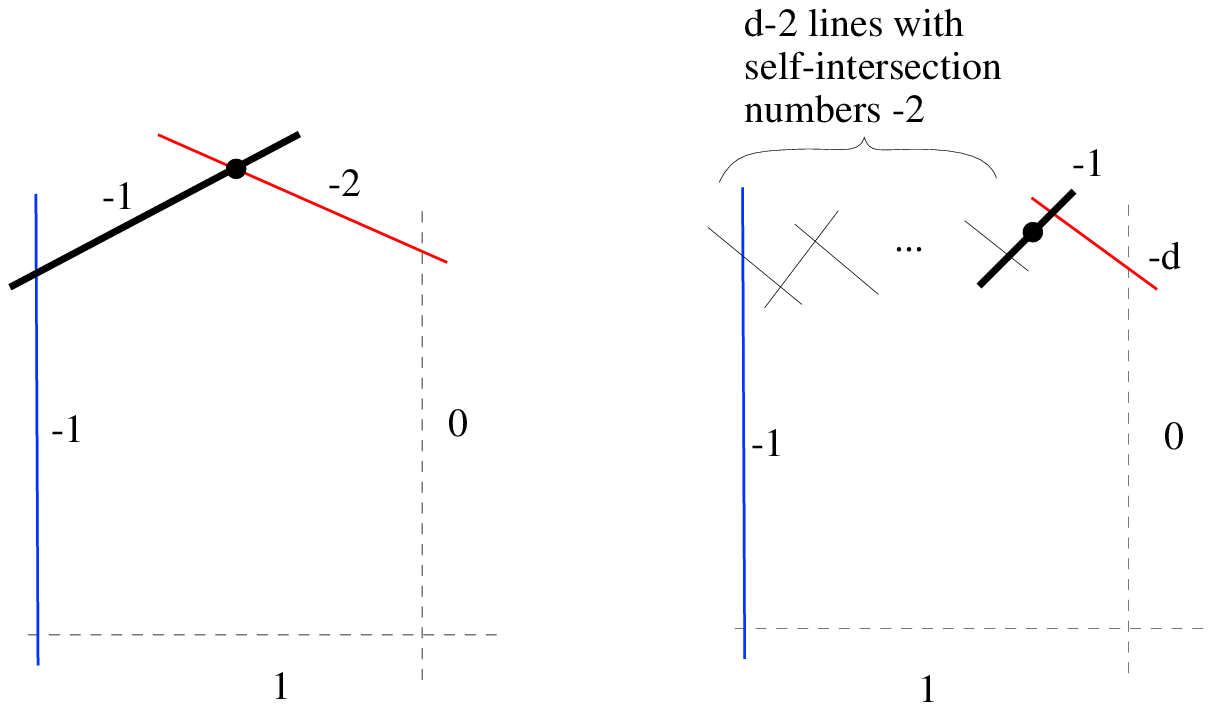 scaled 800}}
\longfigcap \secondblowuppic{The configuration after the second and after the $d$-th of the
blow-up. The heavy dots are the points of indeterminacy of $H_2$ and $H_d$; note that the dot on
the right is an ordinary point; all the earlier ones except the very first were double points.
Again, the numbers labeling components are the self-intersection numbers.}

\subheading {The blow-ups which depend on the coefficients of $p$}

The next $d-1$ blow-ups, although not really more difficult than the earlier ones, have rather
more unpleasant formulas, because each occurs at a smooth point of the last exceptional divisor,
and we need to specify this point. To lighten the notation, we will define by induction the
polynomials
$$
q_0(X) = q(X) \quad \text{and}\quad q_{k+1}(X) = \frac {q_k(X)-q_k(0)}X,\quad k=1,\dots
$$
and the numbers $Q_k$ (they are really coordinates of points on
exceptional divisors)
$$
Q_0=1 \quad \text{and}\quad Q_{k+1} = -\sum_{j=0}^k Q_j q_{k-j+1}(0).
$$

Set $Y_0=X_d$, and make the successive blow-ups 
$$
Y_k-aQ_k = X_1Y_{k+1},\ k=0, \dots, d-2. \peqno \successiveYs
$$

\remark This means that $Y_{k+1}$ is the slope of a line through the point $X_1=0, Y_k=aQ_k$. 
In that sense the number $Q_k$ (or rather $aQ_k$), is really the coordinate on the line
parametrized by $Y_k$ of the next point at which to blow up. 

\lemma \computehorribleblowup a) In the coordinates $X_1, Y_k$, the H\'enon mapping is given by the
formula
$$
H_{d+k}:\pmatrix {X_1}\\{Y_k}\endpmatrix \mapsto \bmatrix X_1Y_kq_1(X_1)+a\sum_{j=0}^{k-1}
Q_jq_{k-j}(X_1)+Y_k\\ X_1^{d-k-1}(X_1^kY_k+a\sum_{j=0}^{k-1}
Q_jX_1^j)\\X_1^{d-k}(X_1^kY_k+a\sum_{j=0}^{k-1} Q_jX_1^j)\endbmatrix, \ k=1, \dots, d-1.\peqno
\furthersimpleblowups
$$

b) The mapping $H_{k+d}$ has the unique point of indeterminacy $X_1=0,\ Y_k=aQ_k$,
for $k=1, \dots, d-2$.

c) The mapping $H_{2d-1}$ has no point of indeterminacy, and maps the
last exceptional divisor to $l_\infty \subset \Proj^2$ by an isomorphism.   
\endstat

\proof This is an easy induction: all the work was in finding the formula. To start the
induction, compute the extension of the H\'enon mapping in the chart $X_d-a=X_1Y_1$ (using
formula (\firstkblowups): 
$$ 
H_{d+1}: \pmatrix {X_1}\\{Y_1}\endpmatrix \mapsto
\bmatrix(a+X_1Y_1)(q(X_1)-1)-X_1Y_1\\X_1^{d-1}(a+X_1Y_1)\\X_1^d(a+X_1Y_1)\endbmatrix= 
\bmatrix X_1Y_1q_1(X_1)+aq_1(X_1)+Y_1\\X_1^{d-2}(X_1Y_1+a)\\X_1^{d-1}(X_1Y_1+a)\endbmatrix,
\peqno \firstsimpleblowup
$$
where we have used $q(X_1)-1=X_1q_1(X_1)$, and factored out $X_1$. Observe that formula
(\firstsimpleblowup) is exactly the case $k=1$ of formula
(\furthersimpleblowups).  

Now suppose Lemma \computehorribleblowup\ is true for $k$, and substitute $Y_k=X_1Y_{k+1}+aQ_k$
from Equation (\furthersimpleblowups).  For the first coordinate of $H_{k+1}$ we find 
$$\eqalign{
X_1(X_1Y_{k+1}&+aQ_k)q_1(X_1)+ a \sum_{j=0}^{k-1} Q_jq_{k-j}(X_1)+X_1Y_{k+1}+aQ_k \cr
&=X_1^2Y_{k+1}q_1(X_1) +aX_1Q_kq_1(X_1)+a \sum_{j=0}^{k-1}
Q_j(q_{k-j}(X_1)-q_{k-j}(0))+X_1Y_{k+1}\cr
&=X_1^2Y_{k+1}q_1(X_1) +aX_1Q_kq_1(X_1)+a X_1 \sum_{j=0}^{k-1}
Q_j q_{k-j+1}(X_1)+X_1Y_{k+1}\cr
&=X_1\left(X_1Y_{k+1}q_1(X_1) +a \sum_{j=0}^{k}
Q_j q_{k-j+1}(X_1)+Y_{k+1}\right).}
$$
The second and third coordinates are similar.  In particular, we see that we can factor out
$X_1$, until $k=d-1$. This proves part (a).

At each step, any points of indeterminacy must be on the last exceptional divisor, of
equation $X_1=0$.  But if we substitute $X_1=0$ in the first coordinate, we find $Y_k=aQ_k$,
so that indeed there is only one point of indeterminacy, proving (b).

The restriction to the last exceptional divisor of the mapping $H_{2d-1}$ is given
by $$
Y_{d-1} \mapsto \bmatrix Y_{d-1}-aQ_{d-1}\\a\\0\endbmatrix
$$
and since $a\ne 0$, the map is well-defined.  Moreover, since $Y_{d-1}$ appears in the first
coordinate with degree 1, this last exceptional divisor maps by an isomorphism to the line at
infinity.
\qed

\bigskip
\centerline{\BoxedEPSF{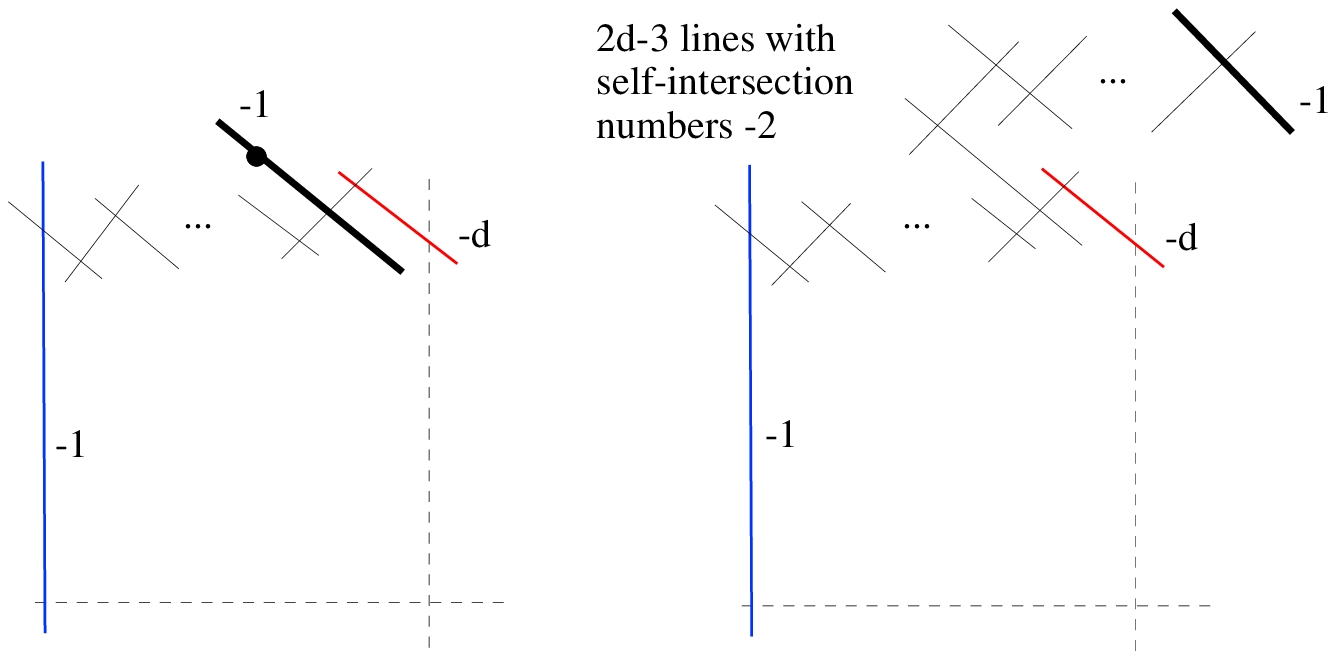 scaled 800}}
\longfigcap \thirdblowuppic {The configuration after the $(d+1)$-st blow-up, and after all
$2d-1$ blow-ups.}

The self-intersection numbers are marked in Figure \thirdblowuppic.  We are now always blowing
up an ordinary point of the last exceptional divisor, so this last exceptional divisor (now
next-to-last) acquires self-intersection -2.  At the end, the very last exceptional divisor
keeps self-intersection -1.  

To summarize, we have proved the following result. Denote by $\tilde X_H$ the space obtained
from $\Proj^2$ be the sequence of $2d-1$ blow-ups described above.

\theorem \HenonWellDefined 
The H\'enon map $H:\C^2 \to \C^2$ extends to a morphism $H_{2d-1}=\tilde H:\tilde X_H \to \Proj^2$,
and maps the divisor at infinity $\tilde {D} =\tilde X_H-\C^2$ to $l_\infty$, mapping all of 
 $\tilde {D}$  to the point $\bold q$ except the last
exceptional divisor which is mapped to $l_\infty$ by an isomorphism.
\endstat

\subheading{Terminology}

To state the next result and later on, we need to give names to the irreducible components of
$\tilde{D}$.  Let us label $A'$ the proper transform of the line at infinity, 
then $B'$ the proper transform of the first exceptional divisor, then in order of creation  $L_1,
L_2, \dots , L_{2d-3}$, and finally $\tilde A$ the components of the divisor $\tilde{D}$. The
line $\tilde A$, i.e. the  last exceptional divisor, will play a special role.

Since $A'$ is the proper transform of $l_\infty$, the projection $A' \to l_\infty$ is an
isomorphism, and we can define
$\bold {p'}, \bold q'$ the points of $A'$ which correspond to $\bold p, \bold q$. Note
that $\{\bold p'\}=L_1 \cap A'$.

The points $\tilde {\bold p}= \tilde H^{-1}(\bold p),\ \tilde {\bold q}= \tilde H^{-1}(\bold q)$
will play a parallel role; note that $\{\tilde {\bold q}\}=L_{2d-3} \cap \tilde A$. This
terminology is illustrated in Figure \divatinfty, which represents $\tilde{D}$, i.e., Figure
\thirdblowuppic (right), redrawn in a more symmetrical way.
\bigskip
\centerline{\BoxedEPSF{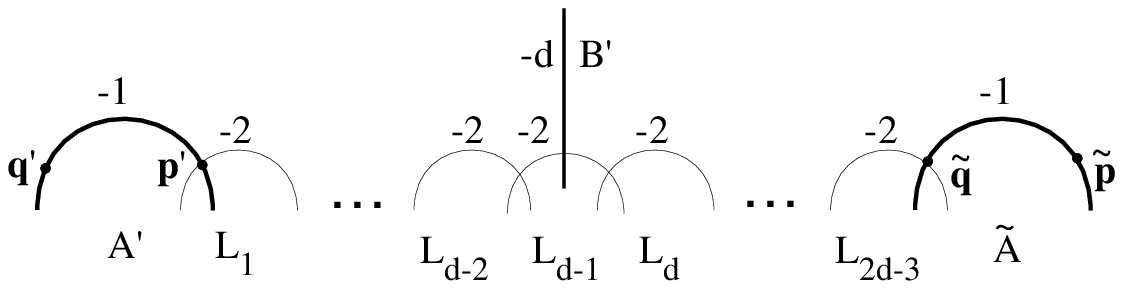 scaled 800}}
\medskip
\longfigcap \divatinfty {The divisor $\tilde{D}$. The numbers labeling the components are
the self-intersection numbers.}

\subheading {The rational mapping $H^\sharp$}

In the next section, we will want to consider $\tilde H$ as a birational map $H^\sharp:
\tilde X_H \ratto \tilde X_H$.

\theorem \tildeHasratmap The rational map $H^\sharp:\tilde X_H \ratto \tilde X_H$ is
defined at all points except $\tilde {\bold p}$.  It collapses $\tilde{D} - \tilde A$ to
$\tilde {\bold q}$, maps $\tilde A-\tilde {\bold p}$ to $A'-\bold p'$ by an isomorphism.

The rational map $(H^\sharp)^{-1}$ is defined at all points except ${\bold q'}$.  It
collapses $\tilde{D}-A'$ to $\tilde {\bold p}$ and maps $A'-{\bold q'}$ to $\tilde
A-\tilde{\bold q}$ by an  isomorphism.
\endstat

\bigskip
\centerline{\BoxedEPSF{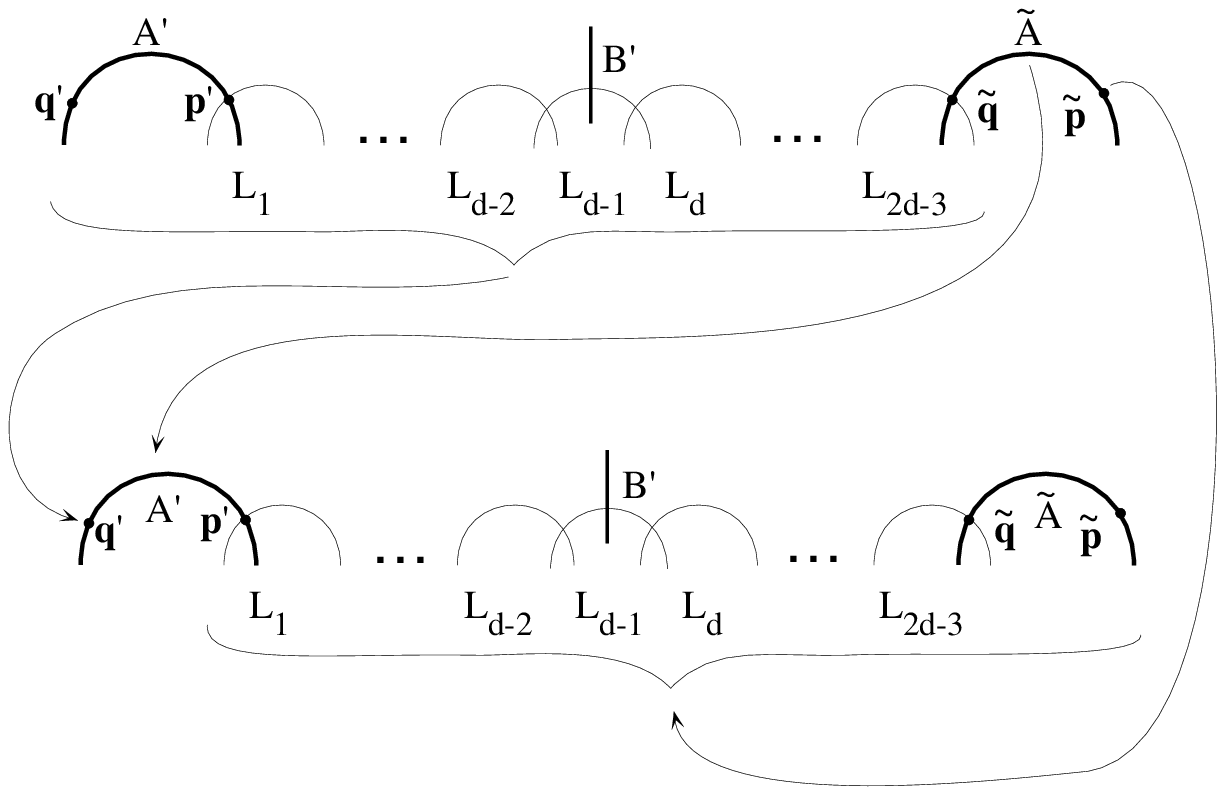 scaled 800}} 
\medskip
\figcap \divatinftywithmap {The ``mapping'' $H^\sharp$ acting on the divisor $\tilde D$.}

\proof This is really a corollary of Theorem \HenonWellDefined. Clearly, $H^\sharp$ is well
defined on $\tilde X_H-\tilde H^{-1}(\bold p)$, i.e., on $\tilde X_H-\{\tilde{\bold p}\}$, and
coincides with $\tilde H$ there.

Further $\tilde H$ is an isomorphism from a neighborhood of $\tilde{\bold p}$ to a neighborhood
of $\bold p$.  So if you perform any sequence of blow-ups at $\bold p$, $\tilde H$ will become
undetermined at $\tilde{\bold p}$.  The statements about the inverse map are similar.
\qed

%% file: graphs+infiniteproducts
\heading IV. Closures of graphs and sequence spaces \endheading

We now have a well-defined map $\tilde H:\tilde X_H \to \Proj^2$, but that does
not solve our problem of compactifying $H:\C^2 \to \C^2$ as a dynamical
system.  We cannot consider $\tilde H$ as a dynamical system, since the domain
and the range are different.  Neither does  $H^\sharp$ solve our problem, since
it still has a point of indeterminacy.  In this section we will show how to
perform infinitely many blow-ups so that in the projective limit we do get a
dynamical system.  We will construct this infinite blow-up as a sequence
space, as this simplifies the presentation and proof (this description was
inspired by Friedland  \cite{Fr1}, who considered the analog in
$\left(\Proj^2\right)^\Z$).  To make this construction, we need to analyze the graph
of $H^\sharp$.

Let $X,Y$ be compact smooth algebraic surfaces, and $f:X \ratto Y$ be a
birational transformation.  Let us suppose that it is undefined at
$\bold p_1,\dots,\bold p_n$, and that $f^{-1}$ is undefined at $\bold
q_1,\dots,\bold q_m$.  Let 
$$
\Gamma_f \subset \left(X-\{\bold p_1,\dots,\bold p_n\}\right) \times Y
$$
be the graph of $f$, and $\overline \Gamma_f \subset X \times Y$ its closure.

\lemma {\wheregraphsmooth} The space $\overline \Gamma_f$ is a smooth
manifold, except perhaps at points $(\x,\y) \in \overline \Gamma_f$ such that 
$$
\x \in \{\bold p_1,\dots,\bold p_n\}\quad \text{and}\quad \y \in \{\bold
q_1,\dots,\bold q_m\}.
$$
\endstat

\proof Clearly $pr_1:\overline \Gamma_f \to X$ is locally an isomorphism near
$(\x,\y)$ if $\x \notin \{\bold p_1,\dots,\bold p_n\}$, and $pr_2:\overline \Gamma_f
\to Y$ is locally an isomorphism near $(\x,\y)$ unless $\y \in \{\bold q_1, \dots,
\bold q_m\}$.
\qed

\example{\singsingraphs} If you have points $(\bold p_i,\bold q_j) \in \overline
\Gamma_f$, they can genuinely be quite singular.  For instance, if 
$X=Y=\Proj^2$ and $f=H$ is a H\'enon mapping, then $f$ (resp. $f^{-1}$) has a
unique point of indeterminacy $\bold p$  (resp. $\bold q$) (see
\pointsofindetlemma). The pair
$(\bold p,\bold q)$ is in $\overline \Gamma_H$, and near $(\bold p,\bold q)$ we can
find equations of $\overline \Gamma_H$ as follows.  

In local coordinates 
$$
u=\frac xy,\  v=\frac 1y\quad\text{near $\bold p$},\quad s=\frac yx,\
t=\frac 1x\quad\text{near $\bold q$}, 
$$
the space $\overline \Gamma_H$ is given by the two equations
$$\eqalign{
v^d&=t\left(\tilde p(u,v)-av^{d-1}\right)\cr
ut&=sv}
$$
which is quite singular at the origin indeed; one way to understand Section III
is as a resolution of this singularity, as Proposition \graphHsmooth\ shows.

Let $H$ be a H\'enon mapping, $\tilde X_H$ be the blow-up on which
$\tilde H:\tilde X_H \to \Proj^2$ is well-defined, and $H^\sharp:\tilde X_H
\ratto \tilde X_H$  be $\tilde H$ viewed as a rational mapping from $\tilde
X_H$ to itself.  

\proposition {\graphHsmooth} The closure $\overline \Gamma_{ H^\sharp} \subset
\tilde X_H \times \tilde X_H$ is a smooth submanifold.
\endstat

\proof The mapping $H^\sharp:\tilde X_H \ratto \tilde X_H$ is birational, and
as we saw in Theorem \tildeHasratmap, it has a unique point of indeterminacy at
$\tilde {\bold p}=\tilde H^{-1}(\bold p)$, and the inverse birational mapping
$(H^\sharp)^{-1}$ also has a unique point of indeterminacy ${\bold q'}$. 
But the point  $(\tilde {\bold p},{\bold q'})$ is {\it not} in $\overline
\Gamma_{H^\sharp}$ , so $\overline \Gamma_{H^\sharp}$ is a smooth (compact)
manifold by Lemma \wheregraphsmooth. \qed

There is another description of $\overline \Gamma_{H^\sharp}$, which we will need
in a moment.

\proposition \graphasfiberedprod The space $\overline \Gamma_{H^\sharp}$, together
with the projections $pr_1$ and $pr_2$ onto the first and second factor
respectively, make the diagram
$$
\CD
 \overline \Gamma_{H^\sharp}  @>{pr_2}>> \tilde X_H  \\ 
 @V{pr_1} VV        @V{\pi}VV   \\ 
  \tilde X_H  @>{\tilde H}>>\Proj^2 
\endCD
$$
a fibered product in the category of analytic spaces \cite{\Douady}.
\endstat

\proof The diagram
$$
\CD
{\tilde X_H} \times \tilde X_H  @>{pr_2}>> \tilde X_H  \\ 
 @V{pr_1} VV        @V{\pi}VV   \\ 
  \tilde X_H  @>{\tilde H}>>\Proj^2 
\endCD
$$
evidently commutes on the graph $\Gamma_H$, and also evidently commutes on a
closed set, hence it commutes on $\overline \Gamma_{H^\sharp}$. 

Since all the
spaces involved are manifolds, it is enough to prove that the diagram is a
fibered product in the category of analytic manifolds, i.e., set-theoretically. 
Since $\pi(\bold y)= \tilde H(\bold x)$ on $\Gamma_H$ this is still true on the
closure $\overline \Gamma_{H^\sharp}$.
\qed

It is time to construct one of our main actors. The space
$X_\infty$, constructed below, is a compact space, which contains $\C^2$ as a
dense open subset, and such that $H:\C^2 \to \C^2$ extends to
$H_\infty:X_\infty \to X_\infty$.

 The locus $D_\infty = X_\infty-\C^2$ is an infinite
divisor at infinity, the geometry of which encodes the behavior of $H$ at
infinity.

\definition \sequencespace Let $X_\infty \subset (\tilde X_H)^\Z$ be the set
of sequences $\underline \x = (\dots, \x_{-1}, \x_0, \x_1, \dots)$ such that
successive pairs belong to $\overline \Gamma_{H^\sharp} \subset \tilde X_H
\times \tilde X_H$ above.  

Let $H_\infty:X_\infty \to X_\infty$ be the shift map 
$$
(H_\infty(\underline {\bold x}))_k = \bold x_{k+1},
$$
where $\underline {\bold x} = (\dots,\x_1, \x_0, \x_1, \dots)$ is a point of
$X_\infty$.

Clearly $X_\infty$ a compact space, since it is a closed subset of a product of
compact sets, and $H_\infty$ is a homeomorphism $X_\infty \to X_\infty$. We will
see below why $H_\infty$ can be understood as an extension of $H$.  

\proposition {\pointsofinfproduct} a) The points of $X_\infty$ are of one of three
types:
\roster
\item Sequences with all entries in $\C^2$;
\item Sequences of the form $(\dots,\tilde {\bold p},\tilde {\bold
p},\bold a,\bold b,\bold q',\bold q',\dots)$ with $\bold a \in \tilde D$, $\bold a
\ne\tilde{\bold p}$; 
\item The two sequences ${\bold p}^\infty = (\dots,\tilde {\bold p},\tilde {\bold p},\dots)$ and 
$\bold q^\infty = (\dots,\bold  q',\bold q',\dots)$. \endroster

b) The sequences of type (1) are dense in $X_\infty$.
\endstat

\proof If a sequence has any entry in $\C^2$, then it is the full orbit of that
point, forwards and backwards.  Otherwise, all entries are in the divisor
$\tilde {D} = \tilde X_H-\C^2$. If these entries are not all $\tilde {\bold p}$, or
all
$\bold q'$, then there is a first entry $\bold a$ which is not $\tilde {\bold p}$;
it must be preceded by all $\tilde {\bold p}$'s.  What follows it is the orbit of
$\bold a$, which is well-defined.  Note that $\bold b=\tilde H(\bold a)$ may be
$\bold q'$ (this will happen unless $\bold a \in \tilde A$), and all the successive
terms must be $\bold q'$. This proves (a).

For part (b), we must show that a sequence $\underline x$ of type 2 or 3 can be
approximated by an orbit, i.e., that for any $\epsilon>0$ and any integer $N$, there
is a point $\bold y \in \C^2$ such that $d(H^n(x) -y_n)<\epsilon$ when $|n|<N$. If
$\underline x$ is of type 2, we may assume that $x(0)\ne \tilde{\bold p}, \bold
q'$.   Then all iterates of $\tilde H$ and of $\tilde H^{-1}$ are defined and
continuous in a neighborhood of $x(0)$, so any point in this neighborhood and close
to $x(0)$ will have a long stretch of forward and backwards orbits close to
$\underline x$; but every neighborhood of $x(0)$ contains points of $\C^2$, which is
dense in $\tilde X$.

Similarly, the orbit of a point with $|x|$ very large and $y=0$ will approximate
$\bold q^\infty$, and a point with $|y|$ large and $x=0$ will approximate $\bold
q^\infty$.
\qed

\example \badsequencespace The fact that $\C^2$ is dense in $X_\infty$ is not
quite so obvious as one might think, and there are examples of birational maps
where it doesn't happen. For instance, consider the mapping 
$$
f:\bvec xy \mapsto \bvec {xy} x,
$$
a priori well defined on $(\C^*)^2$. Denote $\bold p$ and $\bold q$ the points at
infinity on the
$x$-axis and
$y$-axis respectively. Then the pairs $(r,\bold q)$ belong to $\overline
\Gamma_F \subset \Proj^2 \times \Proj^2$ when $r$ is in the line at infinity,
as do the pairs $(\bold q,r)$.  Thus the sequence space contains points like
$$
(\dots, \bold q,r,\bold q,\bold q, r, \bold q, r,  \bold q, \bold q, \dots)
$$
with symbols $\bold q$ and $r$ in any order.  Such sequences cannot be
approximated by orbits in $(\C^*)^2$; they will also form subsets which have
dimension equal to the number of appearances of $r \ne \bold q$, which may be
infinite.  Such sequence spaces are a little scary, as well as
pathological, and irrelevant to the original dynamical system.

In our case, if we had not blown up $\Proj^2$, $\C^2$ would still have been
dense in the sequence space, but it would have had bad singularities.       
\bigskip

\proposition {\infproductgoodpointssmooth} The space $X_\infty^*=X_\infty -
\{\bold p^\infty,\bold q^\infty\}$ is an algebraic manifold.  

More precisely, 
\roster 
\item The projection $\pi_0$ onto the $0$-th coordinate induces an isomorphism of
the space of orbits of the first type to $\C^2$;
\item
If $\underline {\bold x}=(\dots,\tilde{\bold p},\tilde {\bold p},\bold a,\bold
b,\bold q',\bold q',\dots)$ is a point of the second type, and
$\bold a$ appears in the $k$-th position, then the projection $\pi_k$ onto the
$k$-th position induces a homeomorphism of a neighborhood of $\underline {\bold x}$
onto  $\tilde X_H-\{\tilde {\bold p},\bold q'\}$.
\endroster
\endstat

\proof The first part is clear. For the second, if a point
$\underline \y$ satisfies $\y_k \ne\tilde{\bold p},\bold q'$, then the entire
forward and backwards orbit of $\y_k$ is defined:  forwards it will never land on
$\tilde {\bold p}$, and backwards it will never land on $\bold q'$.  

Let us call $\phi_k:\tilde X_H-\{\tilde {\bold p},\bold q')\to X_\infty$ the map
which maps $\x$ to the unique sequence $\underline \x \in X_\infty$ with $\x_k=\x$.
The change of coordinate map $\phi_l^{-1}\circ\phi_k$ is then simply $H^{l-k}$ on
$\C^2$.

This shows that the coordinate changes are algebraic on the intersections of
coordinate neighborhoods, except for one detail. The set $\C^2 \subset X_\infty$
is exactly the intersection of the images of $\phi_k$ and $\phi_l$ when $|l-k| \ge
2$, but when $l=k+1$, the intersection then contains  the sequences
with entries
$(\dots,\tilde{\bold p},\tilde{\bold p},\bold a,\bold b,\bold q',\bold q',\dots)$
with $\bold a \in \tilde A-\{\tilde {\bold p}\}$ and $\bold b \in A'-\{\bold q'\}$,
in this case also the change of coordinates is given by
$H^\sharp$ and is still algebraic.
\qed

We can now see why $H_\infty$ is an extension of $H$.  On the subset isomorphic
to $\C^2$ formed of sequences in $\C^2$, with $\pi_0$ the isomorphism, we have
$$
\pi_0\bigl(H_\infty(\underline \x)\bigr)= H\bigl(\pi_0(\underline \x)\bigr);
$$
i.e., $\pi_0$ conjugates $H_\infty$ to $H$ on that subset.  

\subheading{Notation}

We will systematically identify $\C^2$ with $\phi_0(\C^2) \subset X_\infty$.  With
this identification, $H_\infty$ does extend $H$ continuously, and algebraically in
$X_\infty^*$. Moreover, we will set
$D_\infty = X_\infty-\C^2$; a picture of $D_\infty$ is given in Figure
\infinitedivisorfig.

 The lines denoted by
$A_i,\ i\in \Z$ are formed of those sequences whose $i$-th entry is in $A'$; each
such line connects the sequences whose $i$-th entry is in $L_1$ (denoted by
$L_{i,1}$) with those whose $(i-1)$-st entry is in $L_{2d-3}$.  In particular, the
points $\bold q_0, \bold p_0, \bold q_1, \bold p_1\in X_\infty$ correspond to the
sequences
$$\eqalign{
&\vbox{ \halign{$#$&#&\dots\hfil $#$,&\hfil $#$,&\hfil $#$,&\hfil $#$,&\hfil $#$,\dots&#\cr
\bold q_0= &(& \tilde{\bold p}&\tilde{\bold q}&\bold q'&\bold q'&\bold q'&)\cr
&&-2&-1&\phantom{-}0&\phantom{-}1&\phantom{-}2&\cr}}\qquad
\vbox{ \halign{$#$&#&\dots\hfil $#$,&\hfil $#$,&\hfil $#$,&\hfil $#$,&\hfil $#$,\dots&#\cr
\bold p_0= &(& \tilde{\bold p}&\tilde{\bold p}&\bold p'&\bold q'&\bold q'&)\cr
&&-2&-1&\phantom{-}0&\phantom{-}1&\phantom{-}2&\cr}}
\cr
\noalign{\medskip}
&\vbox{ \halign{$#$&#&\dots\hfil $#$,&\hfil $#$,&\hfil $#$,&\hfil $#$,&\hfil $#$,\dots&#\cr
\bold q_1= &(& \tilde{\bold p}&\tilde{\bold p}&\tilde{\bold q}&\bold q'&\bold q'&)\cr
&&-2&-1&\phantom{-}0&\phantom{-}1&\phantom{-}2&\cr}}
\qquad
\vbox{ \halign{$#$&#&\dots\hfil $#$,&\hfil $#$,&\hfil $#$,&\hfil $#$,&\hfil $#$,\dots&#\cr
\bold p_1= &(& \tilde{\bold p}&\tilde{\bold p}&\tilde{\bold p}&\bold p'&\bold q'&)\cr
&&-2&-1&\phantom{-}0&\phantom{-}1&\phantom{-}2&\cr}}}
$$

\centerline{\BoxedEPSF{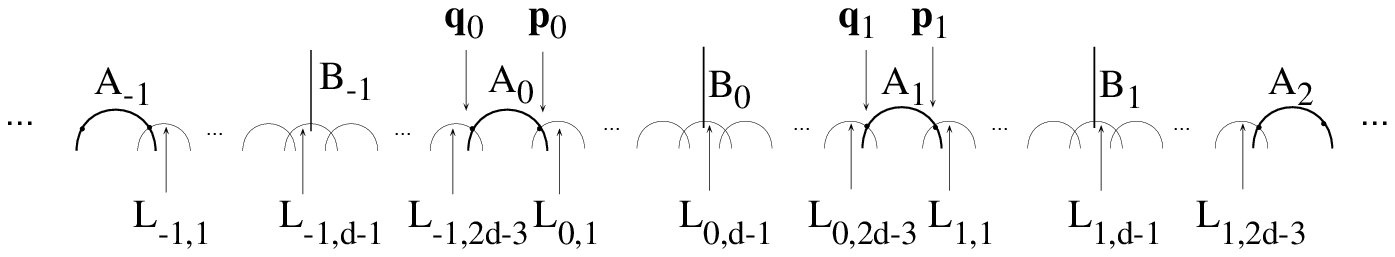 scaled 800}}
\figcap \infinitedivisorfig {The divisor $D_\infty$.}

\heading V. The homology of $X_\infty^*$ \endheading

In this section we will study the homology groups $H_i(X_\infty^*)$, and more
particularly $H_2(X_\infty^*)$ and the quadratic form on it coming from the
intersection product.

Although nasty spaces (solenoids, etc.) are lurking around every corner, here we
will compute only the  homology groups of manifolds, being careful to exclude
the nasty parts.  All homology theories coincide for such spaces, and we may
use singular homology, for instance.  Unless stated otherwise, we use integer
coefficients; at the end we will use complex coefficients.  Using $d$-torsion
coefficients would give quite different results, which can easily be derived
using the universal coefficient theorem.

\subheading{Inductive limits}

It is fairly easy to represent $X_\infty^*$ as an increasing union of subsets whose
homology can be fairly easily computed.  Since inductive limits and homology
commute, it is enough to understand these subsets.  

First some terminology.  If $G$ is an Abelian group, then $G^\N$ is the
product of infinitely many copies of $G$, indexed by $\N$, i.e., the set of all
sequences 
$(g_1,g_2, \dots)$ with $g_i \in G$.  The group $G^{(\N)} \subset G^\N$
is the set of sequences with only finitely many non-zero terms; in the category
of Abelian groups, this is the {\it sum\/} of copies of $G$ indexed by $\N$; it is
also easy to show that it is the inductive limit of the diagram
$$
G \to G^2 \to G^3 \to \dots
$$
where the map 
$$
G^k \to G^{k+1}\quad \text{ is}\quad (g_1, \dots, g_k) \mapsto (g_1, \dots,
g_k,0).
$$

The particular
inductive limit we will encounter is not quite elementary, and we will start
with an example, which has many features in common with our direct limit of
homology groups.

\example {\indlimexample} Consider the inductive system
$$
\Z \overset {f_1} \to {\to} \Z^2 \overset {f_2} \to {\to} \Z^3 \overset {f_3} \to
{\to} \dots,
$$
where $f_n:\Z^n \to \Z^{n+1}$ is defined by
$$
f_n(\bold e_{n,i}) = \cases \bold e_{n+1,i}&\text{if $i<n$}\\
\bold e_{n+1,n}+\bold e_{n+1, n+1}&\text{if $i=n$},
\endcases
$$
using the standard basis $\bold e_{n,1}, \dots \bold e_{n,n}$ of $\Z^n$.

It certainly seems as if the inductive limit of this system should be $\Z^{(\N)}$
of sequences of integers which are eventually $0$.  But this is not true, and
the inductive limit is bigger.

\proposition {\indlimprop} a) If $(v_m,v_{m+1}, \dots)$ represents an element of
$\indlim_n(\Z^n,f_n)$, with $v_m \in \Z^m$, then for any $j$, the coordinate
$(v_m)_j$ is constant as soon as $m>j$.  This defines a map
$$
\indlim_n(\Z^n,f_n) \to \Z^\N
$$
which is easily seen to be injective.

b) The image of $\indlim_n(\Z^n,f_n)$ in $\Z^\N$ consists of the sequences
$(a_j)_{j \in \N}$ which are eventually constant.
\endstat

\proof Any element of the inductive limit has a representative $v_m \in \Z^m$
for some $m$. The $m$-th entry of $v_m$ will be replicated as both the $m$-th and
$(m+1)$-st entry of $v_{m+1}$, and then as the last three entries of $v_{m+2}$,
etc. Clearly the image in $\Z^\N$ will be constant from the $m$-th term on.
\qed

Thus there is an exact sequence
$$
0 \to \Z^{(\N)} \to \indlim_n(\Z^n,f_n) \to \Z \to 0
$$
where the third arrow associates to an eventually constant sequence the value of
that constant.

We see that there is an extra generator to the inductive limit, which one may take
to be the constant sequence of $1$'s in $\Z^\N$.

\example {\indlimextw}  Now let us elaborate our example a little.  Modify $f_n:\Z^n
\to \Z^{n+1}$ so that $f_n(\bold e_{n,n}) = \bold e_{n+1,n} + d \bold
e_{n+1,n+1}$ for some integer $d\ge 1$.

Most of the computation above still holds, except that a sequence
$$
\underline v=(v_1,v_2,\dots) \in \Z^\N
$$
belongs to the inductive limit if and only if it is eventually geometric with
ratio $d$.  We will denote $\Z[1/d]$ the rational numbers with only powers of $d$
in the denominator, i.e., the sub-ring of $\Q$ generated by $\Z$ and $1/d$. If we set
$\underline v^+ = (1,d,d^2,\dots)$, then
$\underline v$ belongs to the inductive limit if and only if there exists $a \in
\Z[1/d]$ such that $\underline v -a \underline v^+$ has only finitely many non-zero
entries.  In other word, there is an exact sequence
$$
0 \to \Z^{(\N)} \to \indlim_N (\Z^n,f_n) \to \Z[1/d] \to 0.
$$

Note that our inductive limit is still a free Abelian group, for there is a
theorem [\Griffiths, Thm. 138] which asserts that a countable subgroup of $\Z^\N$ is
free Abelian. In our case, the elements
$$
(1, d, d^2, d^3, \dots), (0,1,d,d^2,\dots), (0,0,1,d,\dots), \dots
$$
form a basis.  On the other
hand $\Z[1/d]$ is not free (it is divisible by $d$). \bigskip    

\comment
Now to apply this example in our homology computation.  Let $D_\infty =
X_\infty-\C^2$ , and let the subsets 
$D_{(N,M)}\subset D_\infty$ be those
sequences $\underline x$ such that $x_k \ne\tilde{\bold p},\tilde{\bold q}$ for some $k \in
[N,M]$, where $N,M \in \Z \cup\{\pm \infty\}$, and $N \le M$. The subset $\Cal
L_{(N,M)}$ is not closed; let $D_{[N,M]} = \overline{D_{(N,M)}}$. We
invite the reader to verify that
$$
D_{[N,M]} = D_{(N,M)} \cup \{(\dots,\tilde{\bold p}, \underset N-1 \to p',
\tilde q, \dots)\}\cup \{(\dots,\tilde{\bold p}, \underset M+1 \to q',
\tilde q, \dots)\}.
$$
We will set $D_\infty^* = D_\infty-\{p^\infty, q^\infty\}= D_\infty
\cap X_\infty^*
$, observe that with our notation, $D_\infty^* = D_{(-\infty,
\infty)}$. 

Clearly 
$$
X_\infty^* = \bigcup_{N} \left(X_\infty^*-(D_{(-\infty,-N]} \cup
D_{[N,\infty)})\right), $$ so that 
$$
H_2(X_\infty^*) = \indlim_{N} H_2\left(X_\infty^*-(D_{(-\infty,-N]} \cup
D_{[N,\infty)})\right).
$$
\endcomment

\subheading {The homology of blow-ups} 

Before attacking the homology of
$X_\infty^*$, we will remind the reader of some well-known facts about the
homology of algebraic surfaces.

\proposition \removepoints  If $X$ is a smooth algebraic surface (or more generally a
four-dimensional topological manifold), and $Z \subset X$ is a finite subset, then
the inclusion
$X-Z \hookrightarrow X$ induces an isomorphism on 1- and 2-dimensional homology.
\endstat

\proof Consider the long exact sequence of the pair $(X,X-Z)$, which gives
in part $$
\dots \to H_3(X,X-Z) \to H_2(X-Z) \to H_2(X) \to H_2(X,X-Z) \to \dots;
$$
it is enough to show that the first and last term vanish.  Let $(U_{\bold
z})_{{\bold z} \in Z}$ be a set of neigh\-borhoods of the points of $Z$
homeomorphic to 4-balls; by excision, we have 
$$
H_k\bigl(X,X-\{{\bold z}\}\bigr) = 
\bigoplus_{{\bold z}\in P} H_k\bigl(U_{\bold z},U_{\bold z}-\{{\bold z}\}\bigr).
$$
The long exact sequence of such a pair $(U_{\bold z},U_{\bold z}-\{{\bold z}\})$
gives in part 
$$\eqalign{
\dots \to &H_3(U_{\bold z}) \to H_3\bigl(U_{\bold z},U_{\bold z}-\{{\bold
z}\}\bigr)
\to  H_2(U_{\bold z}-\{{\bold z}\}) \to \cr H_2(U_{\bold z})  \to & 
H_2\bigl(U_{\bold z},U_{\bold z}-\{{\bold z}\}\bigr) \to H_1\bigl(U_{\bold
z}-\{{\bold z}\}\bigr)
\to
\dots} 
$$
The first, third, fourth and sixth terms vanish, since $U_{\bold z}$ is contractible and
$U_{\bold z}-\{{\bold z}\}$ has the homotopy type of a 3-sphere. The result for
2-dimensional homology follows; the proof for one-dimensional homology is similar.
\qed

\proposition \removepointsindimthree  If $X$ is a smooth algebraic surface (or more
generally an orientable 4-dimensional topological manifold), and $Z \subset X$ is a
finite subset, then the inclusion
$X-Z \hookrightarrow X$ induces an exact sequence
$$
0 \to H_4(X) \to \Z^Z \to  H_3(X-Z) \to H_3(X) \to  0.
$$
In particular, if $X$ is compact and $Z$ is a single point, then the
inclusion induces an isomorphism $H_3(X-Z) \to H_3(X)$. 
\endstat

The proof comes from considering the same exact sequences as above; we omit it.

We will now see that if you blow up a point of a surface, you increase the
2-dimensional homology by the class of the exceptional divisor. 

Let $X$ be a surface, and ${\bold z}$ a smooth point.  Let $\pi:\tilde X_{\bold z}
\to X$ be the canonical projection, and $E = \pi^{-1}(\bold z)$ be the
exceptional divisor.  

Consider the homomorphism
$$
i:H_2(X) \to H_2(\tilde X_{\bold z}) \peqno \maptohomologyofbu
$$
given by the composition 
$$
H_2(X) \to H_2\bigl(X-\{{\bold z}\}\bigr) \to H_2(\tilde X_{\bold z});
$$
first the inverse of the isomorphism $H_2\bigl(X-\{{\bold z}\}\bigr) \to H_2(X)$ in
Proposition \removepoints, followed by the map induced by inclusion.
  
\proposition \homologyofblowups The map 
$$
H_2(X) \oplus \Z \to H_2(\tilde X_{\bold z})
$$
given by $(\alpha,m) \mapsto i(\alpha)+m[E]$ is an isomorphism. 
\endstat

\proof This is an application of the Mayer-Vietoris exact sequence. Let $U$ be an
open neighborhood of ${\bold z}$ in
$X$ homeomorphic to a 4-ball. Clearly, $\pi$ is a homeomorphism  $\tilde X_{\bold
z}-E \to X-\{{\bold z}\}$. So applying Mayer-Vietoris to the open cover $\tilde
X_{\bold z} -E$ and
$\tilde U_{\bold z}$ of $\tilde X_{\bold z}$ gives in part
$$
\dots \to H_2\bigl(U-\{{\bold z}\}\bigr) \to H_2\bigl(X-\{{\bold z}\}\bigr) 
\oplus H_2(\tilde U_{\bold z}) \to H_2(\tilde X_{\bold z}) \to
H_1\bigl(U-\{{\bold z}\}\bigr) \to \dots.
$$
The first and last terms are zero, because $U-\{{\bold z}\}$ has the topology of a
3-sphere, so the middle mapping is an isomorphism
$$
H_2\bigl(X-\{{\bold z}\}\bigr) \oplus H_2(\tilde U_{\bold z}) \to H_2(\tilde X_{\bold
z}).
$$      
The result follows by applying Proposition \removepoints, since  $\tilde U_{\bold z}$
deforms onto the exceptional divisor, which is a projective line. 
\qed

The same Mayer-Vietoris exact sequence, together with Proposition 
\removepointsindimthree, will also show the following result.

\proposition \blowupsindimoneandthree  If $X$ is compact, the canonical projection
induces isomorphisms\break $H_1(\tilde X_{\bold z}) \to H_1(X)$ and $H_3(\tilde
X_{\bold z}) \to H_3(X)$.
\endstat

The next proposition will be the key to most of our computations.

\proposition \mapinducedonhomology Consider the composition
$$
H_2(X) \to H_2\bigl(X-\{{\bold z}\}\bigr) \to H_2(\tilde X_{\bold z}).
$$
Let $C$ be a curve in $X$ which has $m$ smooth branches through ${\bold z}$.  Then
the image of $[C]$ in $H_2(\tilde X_{\bold z})$ is $[C'] + m [E]$, where
$C'$ is the proper transform of $C$ in $\tilde X_{\bold z}$. 
 \endstat

\proof Let $\hat C$ be the normalization of $C$, which is in particular a
smooth 2-dimensional differentiable manifold, and $f:\hat C \to C$ the
normalizing map. The mapping $f$ can be deformed (differentiably, but perhaps not
analytically) to a map $f':\hat C  \to X$  which avoids ${\bold z}$; so $f'$ lifts to
a map $\tilde f':\hat C \to \tilde X_{\bold z}$.

\medskip
\centerline{\BoxedEPSF{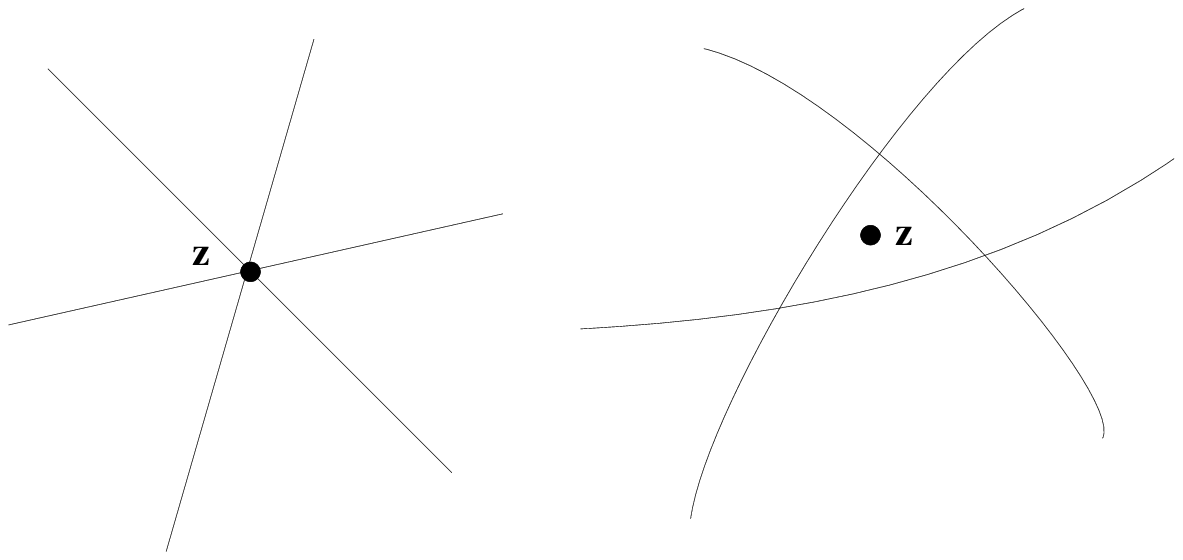 scaled 500}}
\medskip 
\longfigcap \curveanddeform {A curve with 3 smooth branches through ${\bold z}$, 
and a deformation which avoids ${\bold z}$.} 
\bigskip
Then $[\tilde f'(\hat C)]$ is the image of $[C]$ in $H_2(\tilde
X_{\bold z})$.  The homology class $[\tilde f'(\hat C)]$ is of the form $[\tilde
C']+n[E]$ for some
$n$: indeed, $f'$ can be chosen so that  $[\tilde f'(\hat C)]$ is
contained in a small neighborhood of $C' \cup E$, which will retract onto
$C' \cup E$, and hence whose 2-dimensional homology is generated by
$[C']$ and $[E]$. We discover what $n$ is by observing that the
intersection number $[\tilde f'(\hat C)]\cdot[E]$ vanishes, since the
corresponding cycles are disjoint. So
$$
0=[\tilde f'(\hat C)]\cdot[E] = \bigl([C'] + n[E]\bigr)\cdot [E] = m-n,
$$
since each branch of $C$ through $\bold z$ contributes 1 to $[C']\cdot[E]$.
\qed

\subheading {The finite approximations to $X_\infty$}

Now let us consider the set 
$$
X_{[N,M]} \subset \prod_{i=N}^M \tilde X_H, \quad N \le M,
$$
of finite sequences $(\x_N,\x_{N+1}, \dots, \x_M)$ with pairs of successive
points in $\overline \Gamma_{H^\sharp}$, and $D_{[N,M]}\subset X_{[N,M]}$ the
subset with all coordinates in $\tilde{D}$. 

The set $D_{[N,M]}$ contains the point ${\bold
p}^{[N,M]}$ all of whose coordinates are $\tilde {\bold p}$, and the point $\bold
q^{[N,M]}$ all of whose coordinates are $\bold q'$. Let us set
$$
X_{[N,M]}^*= X_{[N,M]}-\left\{{\bold p}^{[N,M]},{\bold q}^{[N,M]}\right\} \quad
\text{and}
\quad D_{[N,M]}^* = X_{[N,M]}^*\cap D_{[N,M]}.
$$

\proposition {\howindlim} If $-\infty \le N'\le N\le M \le M'\le \infty$, the natural
projection 
$$
X_{[N',M']} \to X_{[N,M]}
$$
has an inverse 
$$
X_{[N,M]}^* \to X_{[N',M']}
$$
defined on $X_{[N,M]}^*$.
\endstat

\proof For any point of $X_{[N,M]}^*$, the $N$-th coordinate is not $\bold q'$, so
has a well-defined backwards orbit, and the $M$-th coordinate is not $\tilde{\bold
p}$, so it has a well-defined forwards orbit.  These orbits define an inclusion of
$X_{[N,M]}^*$ into $X_\infty$.  \qed 

The point of this proposition is that we can compute the homology of $X_{[N,M]}^*$. 
If $V$ is an algebraic variety, let $\Irr(V)$ denote the set of irreducible
components of $V$. 

\proposition \DivInfNM a) The space $X_{[N,M]}$ is a smooth algebraic
surface, and $D_{[N,M]}$ is a divisor in $X_{[N,M]}$.

b) The divisor $D_{[N,M]}$ consists of $M+1-N$ ordered blocks, each
consisting of $2d$ projective lines, with the last line of one block coinciding
with the first of the next, as in Figure \finitedivisor.
\endstat

\bigskip
\centerline{\BoxedEPSF{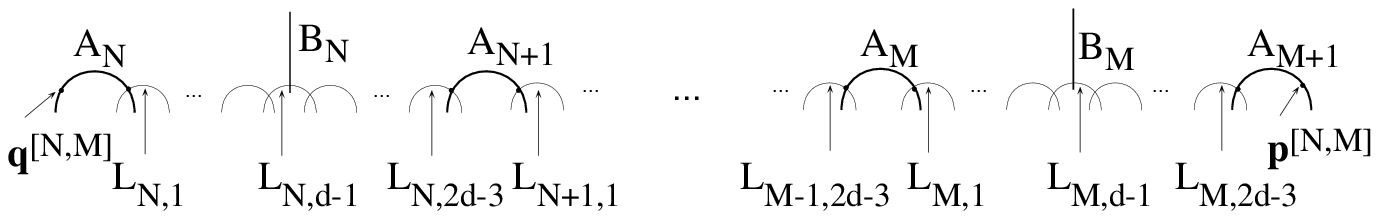 scaled 800}}
\figcap \finitedivisor {The divisor $D_{[N,M]}$.}

\proof Part (a) is more or less obvious, except perhaps for the points 
$\bold p^{[N,M]}, \bold q^{[N,M]}$. The projection onto the $M$-th coordinate gives
an isomorphism of a neighborhood of $\bold p^{[N,M]}$ onto a neighborhood of
$\tilde{\bold p}$, and the projection onto the $N$-th coordinate works for $
\bold q^{[N,M]}$.

A point of $D_{[N,M]}$ will be a sequence of points at infinity in $\tilde
X_H$.  Such a sequence will consist of either

$\bullet$ all $\tilde {\bold p}$ or all $\bold q'$, or 

$\bullet$ a certain number of $\tilde {\bold p}$'s (perhaps none), then a first element
different from $\tilde {\bold p}$, then something (perhaps $\bold q'$), then all
$\bold q'$s.

Let us denote by $D_k$ the $k$-th block
$$
D_k=\set{\underline \x \in D_{[N,M]}}{\x_k \in (\tilde D-\tilde
A)\cup\{\tilde{\bold q}\}\ \text{and}\ \x_{k-1} \ne \bold q'},
$$
for $N\le k \le M$ (if $k=N$, then condition $x_{k-1}\not=\bold q'$ is void).
This set is parametrized by $\x_k \x_k \in (\tilde D-\tilde
A)\cup\{\tilde{\bold q}\}$, and every point of $D_{[N,M]}$ belongs to precisely
one $D_k$, except as follows. 

$\bullet$ The points whose $M$-th coordinate belong to $\tilde A-\{\tilde{\bold
q}\}$; these form a projective line denoted $A_{M+1}$.

$\bullet$ The points $\bold q_k,\ k=N+1, \dots, M$ whose $k$-th coordinate is $\bold
q'$ and whose $k-1$-st coordinate is $\tilde{\bold q}$. The point $\bold q_k$ is
simultaneously the left-most point of $D_k$ and the right-most point of $D_{k-1}$;

$\bullet$ The point $\bold q_{M+1}=A_{M+1} \cap D_M$.
\qed

All lines have the same
self-intersection numbers as the corresponding lines have in $\tilde D$, except
for the connecting lines, i.e., the lines $A_k,\ k=N+1, \dots,M$ where $\x_k \in A'$
and $x_{k-1} \ne \bold q'$, which have self-intersection $-3$, as indicated in the
Figure \selfintfig. This is proved as part of the proof of \homologyoffinblowups.

\bigskip
\centerline{\BoxedEPSF{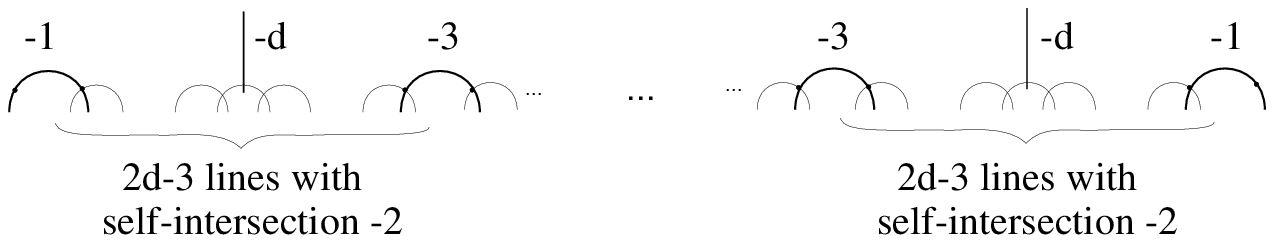 scaled 800}}
\figcap \selfintfig {The self-intersections of the components of $D_{[N,M]}$.}

\proposition \homologyoffinblowups a) The map which associates to each
irreducible component of $D_{[N,M]}$ the 2-dimensional homology class
which it carries induces an isomorphism
$$
\Z^{\Irr(D_{[N,M]})} \to H_2\bigl(X_{[N,M]}\bigr),
$$
when $-\infty<N\le M<\infty$.

b) The inclusion $X_{[N,M]}^* \to X_{[N,M]}$ induces an isomorphism on 2-dimensional
homology. \endstat

Before we prove this proposition, note that it represents the second homology group
of
$X_\infty^*$ as 
$$
\indlim_N \Z^{\Irr(D_{[-N,N]})}, \peqno \homologyasdirlim
$$
since an increasing union of open sets is a inductive limit in the category of
topological spaces, and homology commutes with inductive limits (\cite{Spa}, Chap.
IV, 1.7). This is very similar to example \indlimextw\ above, and we will need to
look carefully at the inclusions.

\proof (a) With a slightly different definition of $X_{[N,M]}$, this follows from
Proposition \homologyofblowups.

We need to know that $X_{[N,M]}$ is obtained from the projective plane
$\Proj^2$ by a sequence of blow-ups, corresponding naturally to the irreducible
components of  $D_{[N,M]}$.  First notice that we may assume that $N=0$:
clearly shifting the indices gives an isomorphism $X_{[N,M]} \to
X_{[0,M-N]}$.  

Next observe that $X_{[0,0]}=\tilde X_H$, which as we saw is obtained from
$\Proj^2$ by a sequence of blow-ups, each of which creates one component of
$\tilde D = D_{[0,0]}$ other than $A'=A_0$.  The component $A'$ is
the proper transform of $l_\infty$ which was there to begin with
and generated the homology $H_2(\Proj^2)$. So the theorem is true when $M=0$.

If $M=1$, notice that $X_{[0,1]}= \Gamma_{H^\sharp}$, so the diagram
$$
\CD
 X_{[0,1]}  @>{pr_2}>> X_{[0,0]}  \\ 
 @V{pr_1} VV        @V{\pi}VV   \\ 
  X_{[0,0]}  @>{\tilde H}>>\Proj^2 
\endCD
$$
is a fibered product by Proposition \graphasfiberedprod.  But the bottom mapping
$\tilde H$ is an isomorphism on a neighborhood of $\tilde {\bold p}$, mapping $\tilde
{\bold p}$ to $\bold p=[0:1:0]$.  Thus the inverse image by $pr_1$ of this
neighborhood maps under $pr_1$ to its image just as the inverse image of the
neighborhood of $\bold p$ maps under $\pi$.

This same argument shows that the component $A_1$ of $D_{[0,1]}$ has
self-intersection 3.  Indeed, the line $A_1\subset D_{[0,0]}$ has
self-intersection $-1$, and the first two blow-ups required to build $X_{[0,1]}$
are blow-ups of points of $A_1$.  

Apply the same argument, using the diagram
$$
\CD
 X_{[0,2]}  @>>> X_{[1,2]}  \\ 
 @VVV        @VVV   \\ 
  X_{[0,1]}  @>>>  \tilde X_H
\endCD
$$
to show that $X_{[0,2]}$ is constructed from $X_{[0,1]}$ by a sequence of
blow-ups, etc.

As above, we see that $A_2$, which has self-intersection $-1$ in $X_{[0,1]}$, has
twice a point blown up, and has self-intersection $-3$ in $X_{[0,2]}$; by induction
$A_k$ will have self-intersection $-1$ in $X_{[0,k+1]}$ and self-intersection $-3$
in $X_{[0,k+2]}$.

(b) This follows immediately from Propositions \removepoints\ and \DivInfNM.
\qed

Next, we need to compute the homomorphism $H_2\left(X_{[-N,N]}\right) \to
H_2\left(X_{[-(N+1),N+1]}\right)$ induced by the composition of the
isomorphism  
$$
H_2\left(X_{[-N,N]}\right) \to
H_2\left(X_{[-N,N]}^*\right) 
$$ 
and the mapping  
$$
H_2\left(X_{[-N,N]}^*\right) \to
H_2\left(X_{[-(N+1),(N+1)]}\right)
$$
induced by the inclusion.  

\proposition{\computeinclusion} The homomorphism
$$
i_N:H_2\left(X_{[N,M]}\right) \to H_2\left(X_{[(N-1),M+1]}\right)
$$
described above is given by the following formula: 
$$\eqalign{
i_N[C] &= [C] \quad\text{if $C \ne A_{M+1}, A_N$} \cr
i_N[A_{M+1}]&= [A_{M+1}]+[B_{M+1}] + 2 [L_{2d-3,M+1}]+3 [L_{2d-4,M+1}]+\dots+\cr
&\hskip 3cm d \Bigl([L_{d-1,M+1}]+ [L_{d-2,M+1}]+ \dots +
[L_{1,M+1}]+[A_{M+2}]\Bigr),\cr
i_N[A_N]&= [A_N]+[B_{N-1}]+ 2[L_{1,N-1}] + 3[L_{2,N-1}]\dots +\cr 
&\hskip 3cm d\Bigl([L_{d-1,N-1}] + [L_{d,N-1}] + \dots
\dots + [L_{2d-3,N-1}]+[A_{N-1}]\Bigr).}
$$
\endstat

\proof  This is a straightforward verification, using Proposition
\mapinducedonhomology. The following sequence of figures should explain exactly the
sequence of blow-ups.
\bigskip
\centerline{\BoxedEPSF{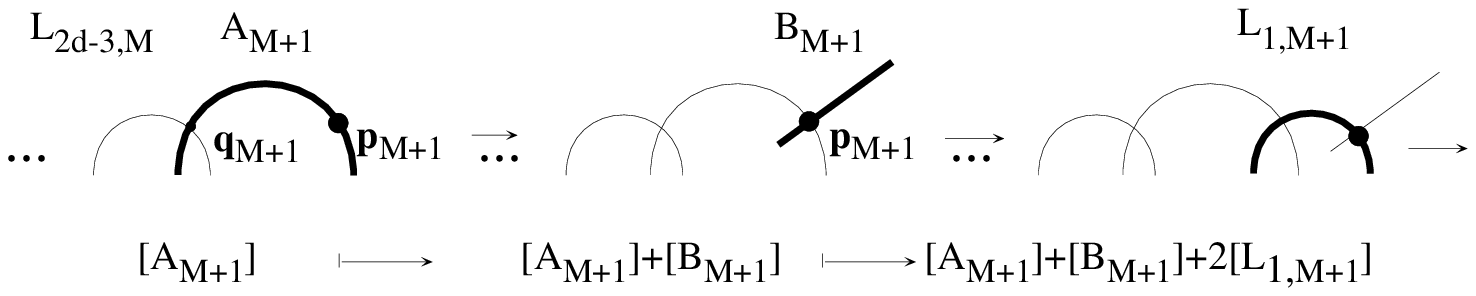 scaled 800}}
\longfigcap \matcompfirstline {The first two blow-ups performed on $A_{M+1}
\subset X_{[N,M]}$. Note that each of $A_{M+1}$ and $B_{M+1}$ contribute 1 to the
coefficient of 
the exceptional divisor $L_{1,M+1}$.} 

\centerline{\BoxedEPSF{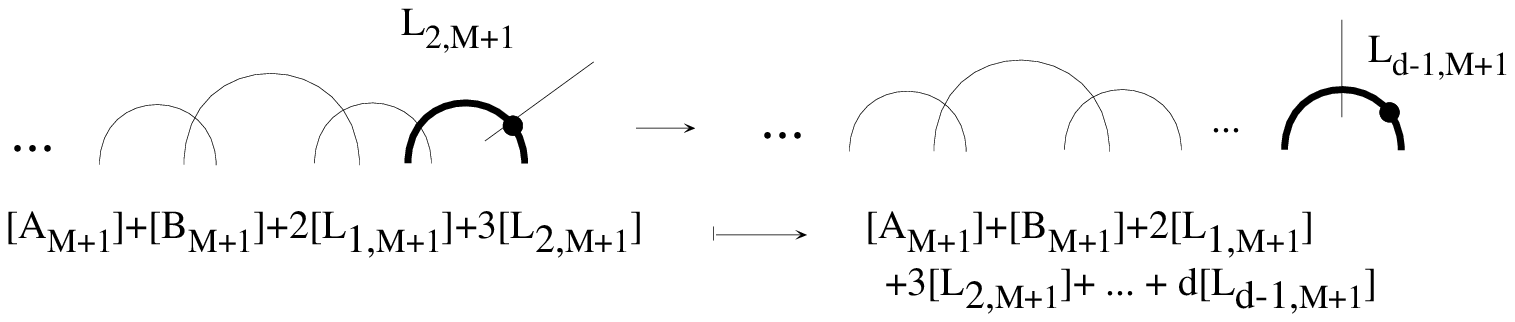 scaled 800}}
\longfigcap \matcompsecondline {The next blow-up and the configuration after $d$
blow-ups. For the figure on the right,  $B_{M+1}$ contributes 1 and $2L_{1,M+1}$
contributes 2 to the coefficient of the exceptional divisor $L_{2,M+1}$.} 

\centerline{\BoxedEPSF{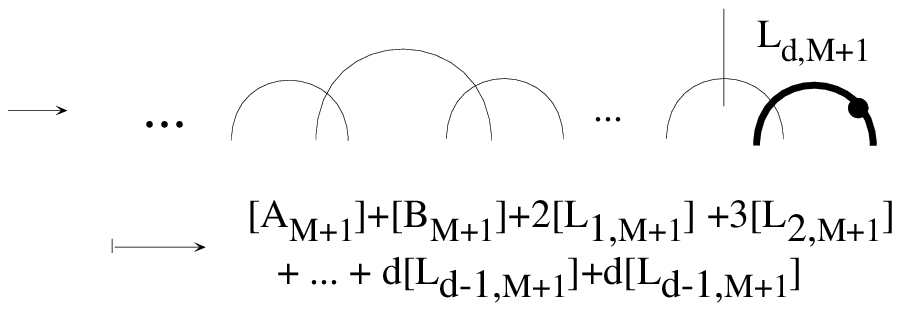 scaled 800}}
\longfigcap \matcompthirdline {The configuration after $d+1$
blow-ups. $dL_{d-1,M+1}$ contributes $d$ to the coefficient of the exceptional
divisor
$L_{d,M+1}$, and it is the only contribution since this time we are blowing up an
ordinary point.}

\centerline{\BoxedEPSF{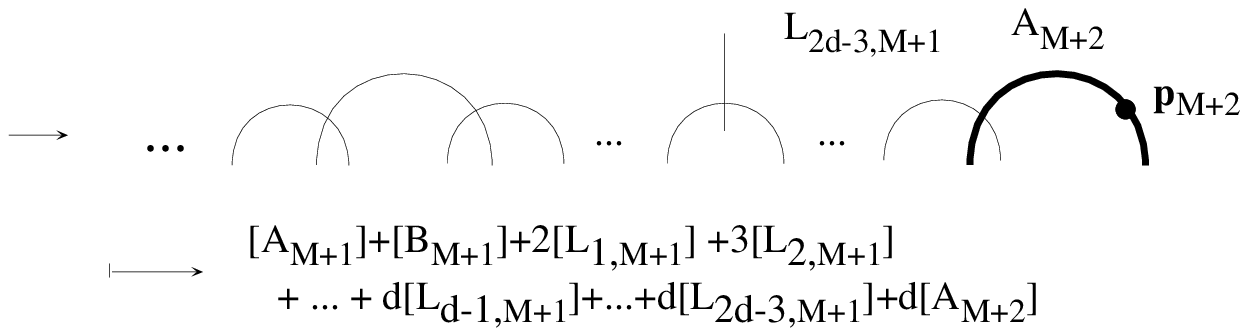 scaled 800}}
\longfigcap \matcompfourthline {The configuration after all the blow-ups required to
pass from $X_{[N,M]}$ to $X_{[N,M+1]}$ have been made. We have blown up ordinary
points on lines with weight $d$, so the new exceptional divisor always has weight
$d$.} \qed

\theorem \embedinproduct a) The inductive limit
$$
\indlim_NH_2(X_{[-N,N]})=H_2(X_\infty^*)
$$
embeds naturally in $\Z^{\Irr(D^*_\infty)}$.

b) If $\underline v \in \Z^{\Irr(D^*_\infty)}$ is an element of $H_2(X_\infty^*)$, then
the limits
$$
\nu^+(\underline v) = lim_{n \to \infty} \frac {\underline a(A_n)}{d^n}\quad
\text{and}\quad \nu^-(\underline v) =lim_{n \to \infty} \frac {\underline
a(A_{-n})}{d^n}
$$
both exist, since the sequences are eventually constant.

c) The sequence
$$
0 \to \Z^{(\Irr(D_\infty^*))} \to H_2(X_\infty^*) \overset{(\nu^+, \nu^-)} \to {\to}
Z[1/d]\oplus Z[1/d] \to 0  \peqno \exactseqforhomologyeq
$$
is exact.
\endstat

\proof a) Any element $\underline v$ of the inductive limit is the image of some 
$\underline v_N\in H_2(X_{[-N,N]})= \Z^{\Irr(D_{[-N,N]})}$ for all sufficiently
large 
$N$, and the coefficient $\underline v_N(L)$ of any irreducible component $L\in
\Irr(D_{[-N,N]})$ is then the same as the coefficient $\underline v_{N'}(L)$ for
all $N'\ge N$ by Proposition \computeinclusion.  This proves (a).

The element $\underline v$
of the limit is entirely determined by the corresponding element of
$v_N \in H_2(X_{[-N,N]})$.  In particular,
$v_N$ assigns some integer weights $\alpha$ to $\left[A_{-N}\right]$ and $\beta$ to
$\left[A_{N+1}\right]$. Then, again by Proposition \computeinclusion, we see that
$$\eqalign{
&\underline v(A_{-N-1}) = d\alpha,\ \underline v(A_{-N-2})= d^2\alpha, \dots\cr
&\underline v(A_{N+2}) = d\beta,\ \underline v(A_{N+3})= d^2\beta, \dots.}
$$
In particular, the sequences defining $\nu^-$ and $\nu^+$ are constant after $N$,
so the limits exist. This proves (b).

For any element $\underline v \in \Z^{(\Irr(D_\infty^*))}$, there exists $N$ such
that
$\underline v$ has coefficient
$0$ for all irreducible components $L \in \Irr(D_\infty^*$ which do not belong to
$D_{[-N,N]}$. Then $\underline v$ is in the image of $H_2(X_{[-(N+1), N+1]}$, and
we see that $\Z^{(\Irr(D_\infty^*))}$ is included in $H_2(X^*_\infty)$. Clearly it
is the kernel of the mapping $(\nu^-, \nu^+)$, which is surjective.
\qed

The exact sequence (\exactseqforhomologyeq) is naturally split: 
call $v^-,\ v^+  \in H_2\left(X_\infty^*\right)$ the images of $[A']$ and
$[\tilde A]$ in $H_2(X_{[0,0]})= H_2(\tilde X_H)$ under the inclusions
$$
 H_2(X_{[0,0]}) \to H_2(X_{[-1,1]}) \to H_2(X_{[-2,2]}) \to \dots.
$$
Then $\nu^+(v^+)= \nu^-(v^-)=1$, and another way of stating part (c) is that 
for any element $\underline v$ of the inductive limit
$\indlim_NH_2(X_{[-N,N]})$ there exists a unique pair \
$(a,b)=(\nu^-(\underline v), \nu^+(\underline v)) \in \Z[1/d]\oplus \Z[1/d]$ such
that 
$ v - a v^- - b  v^+$ belongs to $\Z^{(\Irr(D_\infty))}$.

\theorem \exactsequenceforhomology  The H\'enon mapping $H_\infty:X_\infty^* \to
X_\infty^*$ induces a commutative diagram 
$$
\CD
0     @>>>  \Z^{(\Irr(D_\infty^*))}  @>>> H_2(X_\infty^*) 
@>(\nu^-,\nu^+)>> \Z[1/d]\oplus \Z[1/d]  @>>>  0  \\ 
@VVV  @V{\alpha} VV        @VH_2(H_\infty)VV  @VV{\beta}V @VVV    \\ 
0     @>>>  \Z^{(\Irr(D_\infty^*))}  @>>> H_2(X_\infty^*) 
@>>> \Z[1/d]\oplus \Z[1/d]  @>>>  0, 
\endCD
$$
where $\alpha$ is the shift 
$$
\alpha([A_k])= [A_{k-1}]\quad,\quad  \alpha([B_k])= [B_{k-1}]\quad,
\quad  \alpha([L_{i,k}])= [L_{i,k-1}]
$$
and $\beta$ is the mapping $\beta(a,b)=(a/d, bd)$.
\endstat

\proof  The action of $H$ on the homology is induced by
shifting (to the left) by one block in $\Z^{\Irr(D_\infty^*)}$.
Clearly, this induces the same shift on $\Z^{(\Irr(D_\infty^*))}$,
and the statement about $\alpha$ is true. The see that $\beta$ is
correct, consider a homology class $ x \in \Z^{\Irr(D_\infty^*)}$
in the image of $H_2\left(X_{[-N,N]}\right)$. It will satisfy  $$ x_{A_{N+1}}=b,\
x_{A_{N+2}}=db,\ x_{A_{N+3}}=d^2b, \dots, $$
for some $b\in \Z$, and $\nu^+( x) = b/d^N$. The sequence
$(H_2(H_\infty))( x)$ is the same sequence shifted, so that 
$$
\bigl(H_2(H_\infty)\bigr)( x)_{A_{N+1}}=db,\ \bigl(H_2(H_\infty)\bigr)(
x)_{A_{N+2}}=d^2b, \ \bigl(H_2(H_\infty)\bigr)( x)_{A_{N+3}}=d^3b, \dots,
$$
and 
$$
\nu^-\bigl(H_2(H_\infty)\bigr)( x)= \frac {db}{d^N} = d \nu^-( x).
 $$
The computation for $\nu^+$ is similar.
 \qed

One way of understanding the exact sequence (\exactseqforhomologyeq) is as part of
the homology exact sequence of the pair $D_\infty^* \subset X_\infty^*$.

\proposition \exactsequenceforpair a) There exists a unique isomorphism
$$
\Z[1/d]\oplus \Z[1/d] \to H_2(X_\infty^*,D_\infty^*)
$$
which makes the diagram
$$
\CD
H_2(X_\infty^*) @>>> \Z[1/d]\oplus \Z[1/d]\\
@V{Id}VV  @VVV\\
H_2(X_\infty^*) @>>> H_2(X_\infty^*,D_\infty^*)
\endCD
$$
commute.

b) Both $H_3(X_\infty^*)$ and $H_3(X_\infty^*, D_\infty^*)$ are isomorphic
to $\Z$, and the canonical map
$$
H_3(X_\infty^*) \to H_3(X_\infty^*, D_\infty^*)
$$
is an isomorphism.

c) Both $H_1(X_\infty^*)$ and $H_1(X_\infty^*, D_\infty^*)$ are zero.
\endstat

\remark We will see in Section IX that the homology group
$H_2(X_\infty^*,D_\infty^*)$ can also be understood as $H_1(S^3-(\Sigma^+ \cup
\Sigma^-))$, where $\Sigma^+$ and $\Sigma^-$ are solenoids embedded in a
3-sphere obtained by an appropriate real oriented blow-up.  A classical result of
algebraic topology asserts that for the standard $d$-adic solenoid $\Sigma_d$
embedded in the 3-sphere in the standard way, $H_1(S^3-\Sigma_d)= \Z[\frac 1d]$.  
This explains why these bizarre groups are appearing in this complex-analytic
setting, by making precise sense of the sentence ``at infinity,
$D^*$ has two solenoids''.

\proof  The exact
sequence of the pair $D_\infty^* \subset X_\infty^*$ reads in part
$$ 
H_2(D_\infty^*) \to H_2(X_\infty^*) \to H_2(X_\infty^*,\Cal
D_\infty^*)\to H_1(D_\infty^*),
$$
Clearly the first term is precisely $\Z^{(\Irr(D_\infty^*))}$, and the
last terms vanishes, since $D_\infty^*$ is a union of 2-spheres identified
at points, with the quotient topology from the disjoint union. This proves
(a). 

Another part of the long exact sequence reads
$$
H_3(D_\infty^*) \to H_3(X_\infty^*) \to H_3(X_\infty^*, D_\infty^*)
\to H_2(D_\infty^*) \to H_2(X_\infty^*)
$$
and the fact that the last map is injective and that $H_3(\Cal
D_\infty^*)=0$ says that the canonical map $$ H_3(X_\infty^*) \to
H_3(X_\infty^*, D_\infty^*) $$
is an isomorphism.

To see what it is an isomorphism between, notice first that
$H_3(\Proj^2)=0$, and it then follows from \blowupsindimoneandthree\ that
$H_3(X_{[N,M]})=0$ for all
$N,M$. Next, the inclusion
$$
X_{[N,M]}^* \in X_{[N,M]}
$$
induces (still by \blowupsindimoneandthree\ for the final $0$) an exact sequence
$$
0 \to \underset {\overset \parallel \to {\Z}} \to{ H_4(X_{[N,M]})} 
{\to}
\underset{\overset \parallel \to{\Z^2}}
\to {H_4(X_{[N,M]},X_{[N,M]}^* )}
\to  H_3(X_{[N,M]}^*) \to
\underset{\overset \parallel \to {0}}
\to {H_3(X_{[N,M]})}.
$$
Thus $H_3(X_{[N,M]}^*)$ is canonically the quotient of $\Z^2$ by the image of $\Z$
under the diagonal map, i.e., it is isomorphic to $\Z$.

The argument for (c) is similar but easier.
\qed 

\subheading{The intersection form on the homology}

The homology space $H_2\left(X_\infty^*\right)$ carries a quadratic form
coming from intersection.  We can make it explicit as follows.

\proposition \computeintsonvs On $\Z^{\Irr(D_\infty^*)}$, the quadratic form is
determined by the self-intersections and mutual intersections of the irreducible
components of $D_\infty^*$.

The classes $v^+$ and $v^-$ satisfy the following rules:
$$\eqalign{
v^+ \cdot v^+ = v^- \cdot v^- &= -1\quad,\ v^+ \cdot v^- = 0\cr
v^+ \cdot [L_{2d-3,0}] =v^-\cdot [L_{1,0}] &=1,\ v^+\cdot [A_1]=v^-\cdot
[A_{0}]=-1}  \peqno \intersectformseq
$$
with all other intersections $0$.
\endstat

\proof The statement about $\Z^{\Irr(D_\infty^*)}$ should be clear. 

For the other classes, one way to do it is to construct a differentiable surface
$C^+\subset \tilde X_H$ (not an algebraic curve), which represents $\tilde A$, and
which avoids $\tilde {\bold q}$ and $\bold q'$. Note that $C^+$ cannot be algebraic
(or analytic): the self-intersection of $\tilde A$ is $-1$ so it is {\it rigid\/} as
an algebraic curve. The curve $C^+$ is then contained in $X_\infty^*$ and represents
$v^+$. But a neighborhood of the curve $C^+$ is also contained in $X_\infty^*$, so
$v^+$ only intersects curves of $D_\infty$ which $C^+$ intersects. Thus $v^+
\cdot v^+ = C^+ \cdot [A_0]=C^+\cdot C^+= -1$, and $v^+ \cdot [L_{2d-3,0}]=C^+ \cdot
[L_{2d-3,0}]=1$.

Similarly, construct a differentiable surface $C^-\subset \tilde X_H$ which is a
deformation of $A'$; clearly we can take $C^+ \cap C^-=\emptyset$.
  \qed

 Of course the quadratic form is invariant under the action of $H_\infty$, since
this is a homeomorphism of $X_\infty^*$.  It certainly isn't obvious from the
formulas. Let us check one case, just for consistency's sake.  Take $d=2$. 
Since $H_\infty$ induces the shift,  we see that
$$
(H_\infty)_*(v^+) = 2 v^+ +[A_{0}]+[B_0]+2[L_{1,0}];
$$

The intersection product gives
$$\eqalign{
\bigl((H_\infty)_*(v^+)\bigr)^2 &= 4(v^+)^2+(A_{0})^2+(B_0)^2+4(L_{1,0})^2+\cr
&\qquad 4B_0 \cdot
L_{1,0}+4A_{0}\cdot L_{1,0}+8v^+\cdot L_{1,0}\cr
&\quad  = -4-3-2-8+4+4+8=-1}
$$
as it should.

This quadratic form on $H_2(X_\infty^*)$ is of course neither positive nor negative
definite.  For instance $\Delta$, the closure of the diagonal of $\C^2$ in
$X^*_\infty$, has self-intersection $+1$, whereas all the irreducible components of
$\Cal D_\infty^*$ have negative self-intersection. The following proposition says
that the form is mainly negative.

\theorem \signthm The intersection form is negative definite on
$\Z^{(\Irr(D_\infty^*))}$.
\endstat

\proof  We will give two proofs, one conceptual and one computational.  Each
proves a stronger (but not the same stronger) result.

\noindent{\bf First proof.} An element $ v \in \Z^{(\Irr(\Cal
D_\infty^*))}$ comes from the homology of some $X_{[N,M]}$ which assigns
coefficient $0$ to the first and the last exceptional divisor. Its self-intersection
in $X_{[N,M]}$ and in $X_\infty^*$ coincide. We will in fact prove that if $ v$ has
coefficient $0$ with respect to one of these, then
$( v\cdot  v) <0$ unless $ v = 0$.

Indeed, the complement of the last exceptional divisor in $D_{[N,M]}$,
can be blown down to a point, so by a theorem of Grauert [G], the intersection
matrix of this divisor is negative definite, hence $ v
\cdot  v$ is negative if $v \ne 0$.

\noindent{\bf Second proof.}
  Let us call
$a_n,b_n, x_{i,j}$ the coefficients of $A_n,B_n, L_{i,j}$ respectively.  Thus we are
considering the quadratic form
$$\eqalign{
\dots&-3a_{n}^2+2a_nx_{n,1}-2x_{n,1}^2+2x_{n,1}x_{n,2} + \dots
-2x_{n,d-2}^2+2x_{n,d-2}x_{n,d-1}\cr
&-d b_n^2 +2b_nx_{n,d-1}
-2x_{n,d-1}^2+2x_{n,d-1}x_{n,d}+\dots +2x_{n,2d-3}a_{n+1}-3a_{n+1}^2 + \dots.}
$$
It is clearly enough to show that the quadratic form obtained by allocating half the
coefficient $a_n$ to the next term and half to the previous term is negative
definite, i.e., that the quadratic term in $2d$ variables
$$\eqalign{
&-\frac 32 a_0^2+2a_0x_{1}-2x_{1}^2+2x_{1}x_{2} + \dots
-2x_{d-2}^2+2x_{d-2}x_{d-1}\cr
&-d b^2 +2bx_{d-1}
-2x_{d-1}^2+2x_{d-1}x_{d}+\dots +2x_{2d-3}a_1-\frac 32 a_1^2}
$$
is negative definite. This is something like working in one block at a time.

If we isolate the terms containing $b$, and complete squares, we find that this
quadratic form can be written
$$\eqalign{
&-\Bigl(d b^2 -2bx_{d-1}+\frac 1d x_{d-1}^2\Bigr)\cr
&-\frac 32 a_0^2+2a_0x_{1}-2x_{1}^2+2x_{1}x_{2} + \dots
-2x_{d-2}^2+2x_{d-2}x_{d-1}\cr
&-\Bigl(2-\frac 1d \Bigr)x_{d-1}^2+2x_{d-1}x_{d}-2 x_d^2 \dots
+2x_{2d-3}a_1-\frac 32 a_1^2.}
$$

If we complete squares from both ends, we can write this as
$$\eqalign{
&-\Bigl(d b^2 -2bx_{d-1}+\frac 1d x_{d-1}^2\Bigr)\cr
&-\Bigl(\frac 32 a_0^2-2a_0x_{1}+\frac 23 x_1^2\Bigr)-\Bigl(\frac 32
a_1^2-2a_1x_{2d-3}+\frac 23 x_{2d-3}^2\Bigr)\cr 
&-\Bigl(\frac 43 x_{1}^2-2x_{1}x_{2} + \frac 34 x_2^2\Bigr)  -\Bigl(\frac 43
x_{2d-3}^2-2x_{2d-3}x_{2d-2} +
\frac 34 x_{2d-2}^2\Bigr) \cr
&-\dots - \dots\cr
&-\Bigl(\frac {d+1}{d}x_{d-2}^2-2 x_{d-2}x_{d-1} +\frac {d}{d+1}x_{d-1}^2\Bigr)
-\Bigl(\frac {d+1}{d}x_{d}^2-2 x_{d}x_{d-1} +\frac {d}{d+1}x_{d-1}^2\Bigr)
\cr
&-\frac {d-1}{d(d+1)} x_{d-1}^2.}
$$
\qed
\goodbreak
It works, with a tiny bit to spare, so we actually get a slightly stronger result:

\proposition \slightstrengthening There exists $K$ depending only on $d$
such that for any $ v \in \Z^{(\Irr(D_\infty^*))}$, we have
$$
\frac 1K \sum v_i^2 \le - v \cdot  v \le K \sum v_i^2.
$$
\endstat

Thus we can complete $\C^{(\Irr(D_\infty^*))}$ with respect to the
intersection inner product, to get a Hilbert space, which we denote  $\hat
H_2^-(X_\infty-\{\bold p^\infty,\bold q^\infty\};\C)$. By Proposition
\slightstrengthening, this intersection product norm is equivalent to the
$l_2$ norm on the space of sequences.

The exact sequence
$$
0 \to \Z^{(\Irr (D_\infty^*))} \to H_2(X_\infty^*) \to \Z[1/d]\oplus \Z[1/d]
\to 0
$$
gives, tensoring with $\C$,
$$
0 \to \C^{(\Irr (D_\infty^*))} \to H_2(X_\infty^*;\C) \to \C\oplus
\C \to 0
$$
so it is natural to complete the entire homology, i.e., to set
$$
\hat H_2(X_\infty^*;\C) = \hat H_2^-(X_\infty^*;\C)\oplus \C v^+ \oplus \C v^-.
$$
On this completed homology space (unlike homology with infinite chains, dual of
cohomology with compact supports), the inner product is still defined (for
instance by the formulas (\intersectformseq)).

Clearly the subspace $\C^{(\Irr(D_\infty^*))}\subset H_2(X^*_\infty;\C)$ is
invariant under the H\'enon mapping $H$, which is simply a shift in $D_\infty$,
so it induces a unitary operator on the Hilbert space $\hat H:
\hat H_2^-(X_\infty-\{p^\infty,q^\infty\};\C)$. This unitary operator has
only continuous spectrum, on the unit circle, and with spectral density
$2d-1$.  There are in addition two eigenvectors of 
$$
(H_\infty)_*:\hat H_2(X_\infty^*;\C) \to \hat H_2(X_\infty^*;\C),
$$
one with eigenvalue $d$ and one with eigenvalue $1/d$.  One way of defining
them is as
$$
w^+ = \lim_{n \to \infty}\frac 1{d^n} (H_\infty)_*^n(v^+) \ \text{and}\   
w^- = \lim_{n \to \infty}\frac 1{d^n} (H_\infty)_*^{-n}(v^-).
$$
These do belong to the completed homology (but not to the homology), since
$w^+$ is $v^+$ on the positive part of $\Cal \L_\infty^*$, and decreases like
a geometric series on the negative part.

These homology classes are already well known in the theory:  they are the
homology classes of the currents $\mu^-$ and $\mu^+$, as defined by \cite{\BSone}.

%% file: Veselov.compactifications.two
\heading 
VI. Real oriented blow-ups of complex divisors in surfaces
\endheading

In this section we will describe a technique which will allow us, among other things, to
understand the structure of $X_\infty$ at the bad points $\bold p^\infty, \bold q^\infty$. Given
a divisor $D \subset X$ in a surface, the idea is to take an open tubular neighborhood $W$
of $D$ in $X$, together with a projection $\pi:\bar W \to D$. Excise $W$, to form
$X'=X-W$.  If $W$ was chosen properly, $X'$ will be a real 4-dimensional
manifold with boundary $\partial X' =\partial \bar W$.  The interior of $X'$ will be
homeomorphic to $X-D$. We will think of $X'$ as some sort of blow-up of $X$ along $D$, with
$\pi: \partial X' \to D$ the exceptional divisor.

This picture is the right one to keep in mind, but it is non-canonical. Further,
the spaces $\partial X', X'$ have no natural algebraic structure (or even analytic). In
this section we will formalize the construction above when $W$ is the ``first-order
infinitesimal tubular neighborhood'' of $D$. This will involve considering $X$ and $D$ as
{\it real algebraic objects}, a 4-dimensional smooth variety $X_\R$ and a singular surface
$D_\R \subset X_\R$.

In our definition \blowupdef\ of blow-ups, and in the associated construction, we
never used that we were working over an algebraically closed field; section II
applies just as well if the underlying field is $\R$ (or any field). Let $Y
\subset X$ be a complex algebraic manifold and a subspace. We will call $X_\R$ the manifold
$X$ seen as a real algebraic manifold, $Y_\R \subset X_\R$ the subspace defined by the ideal
$I_\R(Y)$ of real-algebraic functions on $X_\R$ which vanish on the real-algebraic subset $Y$.

\remark This construction of $Y_\R$ does not coincide with the more obvious idea of taking a
chart $U \subset  X$ with $U \subset \C^n$, and complex equations $f_j:U \to \C$ such that $Y \cap
U$ corresponds to the ideal generated by $f_1, \dots, f_k$.  Then consider the subset of $U
\subset \R^{2n}$ with coordinates $x_1+iy_1=z_1, \dots, x_n+iy_n=z_n$ and the locus defined by the
$2k$ equations $\re f_1=0, \im(f_1)=0, \dots, \re(f_k)=0, \im(f_k)=0$, or equivalently by 
the ideal generated by the $\{\re(f_j),\ \im(f_j)\}$.  Even in simple cases, these two ideals are
different.

\example \axesexampleasrealsubset Consider $X=\C^2$, with $Y$ given by the equation $xy=0$.  Then
the construction above yields the locus in $\R^4$ defined by the equations
$$
x_1x_2-y_1y_2 = 0\quad \text{and} \quad x_1y_2+x_2y_1=0. \peqno \baddefeltseq
$$

But the algebraic set $Y \subset \R^4$ is the union of the $x_1,y_1$-plane and the $x_2,y_2$-plane;
the ideal $I_\R(Y)$ is generated by the four elements
$$
x_1x_2,\ x_1y_2,\ y_1x_2, \ y_1x_2. \peqno \gooddefeltseq
$$
In Section II we saw that this ideal cannot be generated by fewer elements, since the blow-up of
$\C^4$ along the subspace it defines is not a subset of $\Proj^1$ or $\Proj^2$.  Thus the ideals
are different, and in fact the subset defined by the elements (\baddefeltseq) consists in $\C^4$ of
four planes, two of which are our real planes, and the others intersect the real locus only at the
origin. Our space $Y_\R$ is the one given by the ideal generated by the elements (\gooddefeltseq).

The ideal $Y_\R$ appears to be quite difficult to compute in general.  For instance, if $X=\C^2$ and
$Y$ is the curve of equation $x^2+y^3=0$, we do not know how to write generators for the
corresponding ideal $I_\R(Y)$. \endremark

We will study the blow-up of
$X_\R$ along
$Y_\R$. The objects we obtain this way are very different from the blow-up of $X$ along $Y$.  In
particular, when $Y$ is a divisor in $X$, the complex blow-up doesn't change it at all, whereas as a
real variety it is of codimension 2, and the real blow-up is quite different from the
original space, as we will see. We will denote this blow-up by
$\Cal B_\R(X,Y)$, rather than
$(\widetilde{X_\R})_{Y_\R}$, both to lighten notation, and to emphasize the
different, more topological way we think of real blow-ups.

\remark In section II, we were working with analytic spaces, with structure sheaves which
might contain nilpotents, etc.  This is if anything even more important here: the origin in
$\R^2$ defined by the two equations $x=0,\ y=0$ and the origin defined by the single equation
$x^2+y^2=0$ are completely different objects, and their blow-ups have no obvious relation. \endremark

The space $\Cal B_\R(X,Y)$ is a real-algebraic space, but it is not the boundary of a
tubular neighborhood we are after; we will need to modify it a bit.  To avoid interrupting
the formal definition, we isolate a purely topological construction we will need: {\it
prime ends}.

\subheading{Prime Ends}

\definition \primeenddef Let $Z \subset U$ be a topological space and a closed subset.
Define the set of {\it prime ends} $P_z(U,Z)$ of $U-Z$ at $z \in Z$  to be the projective
limit $$
P_z(U,Z)=\underset \leftarrow \to{\lim}\  \pi_0(V-Z)
$$
where the projective limit is over all neighborhoods $V$ of $z$ in $U$ with respect to
inclusions, and
$\pi_0$ is the functor which associates to a topological space its set of connected
components.  

Now consider the {\it prime end modification} of $U$ along $Z$, which is the disjoint
union
$$
\Cal P(U,Z) = (U-Z) \sqcup \coprod_{z\in Z} P_z(U,Z),
$$
where the topology on $U-Z$ in
unchanged, and a basis of neighborhoods of a point $\tilde z \in P_z(U,Z)$ is obtained
by choosing a basis of neighborhoods of $z$ in $U$, and taking the closure in $U$
of the component of $U-Z$ corresponding to $\tilde z$. This space comes with an
obvious projection
$\Cal P(U,Z) \to U$, which is continuous. With slight abuse of notation, we will denote 
$\Cal P(Z) \subset \Cal P(U,Z)$ the inverse image of $Z$ in $\Cal P(U,Z)$.

The following properties of the prime end modification are left to the reader:

$\bullet$ The prime end modification of a Hausdorff space is Hausdorff.

$\bullet$  When $U$ is a locally finite simplicial complex, and $Z \subset
U$ is a simplicial subset, $P_z(U,Z)$  is finite, and the projective limit is achieved by
taking $V$ the star of $z$, i.e., the union of the open simplices with $z$ in their closure. 

$\bullet$ The prime end modification of a simplicial complex along a subcomplex is a
simplicial complex.

$\bullet$ If $F:U_1 \to U_2$ is continuous, $Z_2$ is closed in $U_2$ and
$Z_1=F^{-1}(Z_2)$, then there exists a unique mapping 
$$
\Cal P(F):\Cal P(U_1,Z_1) \to \Cal P(U_2,Z_2)  \peqno \natprimeend
$$
making the diagram
$$
\CD
\Cal P(U_1,Z_1)  @>{\Cal P(F)}>> \Cal P(U_2,Z_2) \\ 
@VVV   @VVV     \\ 
U_1    @>>>  U_2 
\endCD
$$ 
commute.

\subheading{Definition of the real oriented blow-up}

For real algebraic varieties, it makes sense to speak of the oriented blow-up. This takes
two steps:

\noindent{\bf Step 1.}  Let $X$ be an analytic space, and $Y$ a closed subspace. If $Y =
f^{-1}(0)$, where
$f:X \to \C^{n+1}$ is analytic, recall from Section II that the real blow-up is a  subspace $\Cal
B_\R(X,Y)
\subset X_\R\times
\Proj^{2n+1}_\R$, locally defined by a section $\tilde f$ of the tautological line bundle
over  $X_\R\times \Proj^{2n+1}_{\R}$ which is not a zero-divisor. We will be mainly
concerned with the case
$n=0$, so
$Y$ is a divisor in $X$, but $Y_\R$ is of codimension $2$ in $X_\R$.

 The sphere $S^{2n+1}$ is naturally a double cover of
$\Proj_\R^{2n+1}$: to every $s\in S^{2n+1}$ we can consider the line 
$l_s=s\R \in \Proj^{2n+1}$, and
the double cover is simply $s \mapsto l_s$. Consider the inverse image $\orbstar (X,Y)$
of
$\Cal B_\R(X,Y)$ in $X_\R\times S^{2n+1}$, now defined by a section $f^*$ of the
tautological bundle on $X\times S^{2n+1}$ (note that the tautological line bundle over
$S^{2n+1}$ is simply the normal bundle, so that $f^*$ is a normal vector field on $S^{2n+1}$).
Since this bundle is oriented (in fact trivial), it makes sense to consider the subset 
$$
\orbtwostar (X,Y) \subset\orbstar (X,Y)
$$
defined by $f^* \ge 0$, with $\orbtwostar(Y) \subset \orbtwostar(X,Y)$ defined by
$f^*=0$ as the exceptional divisor.

\remark The object of passing to the double cover is to make the tautological bundle oriented,
so as to be able to speak of $f^* \ge 0$.
\endremark 

\noindent{\bf Step 2.} The pair 
$$
\orbtwostar(Y) \subset \orbtwostar(X,Y)
$$
is almost the oriented real blow-up, but sometimes it has singularities which we
want to eliminate, using prime ends. Define the {\it oriented real blow-up} $\orb(X,Y)$
and the exceptional divisor $\orb(Y) \subset \orb(X,Y)$ to be the prime end modification
$$
\orb (Y)=\Cal P(\orbtwostar(Y)) \subset \Cal P \bigl(\orbtwostar(X,Y),
\orbtwostar(Y)\bigr)=\orb(X,Y).
$$
of $\orbtwostar(X,Y)$ along $\orbtwostar(Y)$.

We will denote by $p_{(X,Y)}:\orb(X,Y) \to X$ the canonical projection, and we will omit the
subscript when there is no ambiguity. 

\example \polarcoordexample Let us construct the real oriented blow-up of the origin $Y$ in
$\C$. Set $z=x+iy$, so that as a real algebraic subspace the origin $Y_\R$ is defined by
the equations $x=y=0$, and $\Cal B_\R(\C,Y) \subset \C
\times
\Proj_\R^1$ is the space
$$
\set {\pmatrix x\\y\endpmatrix, \bvec uv \in \R^2 \times \Proj_\R^1}{\ xv=uy}.
$$
The set $\orbstar(\C,Y)$ comes from replacing  $\Proj_\R^1$ by $S^1$, i.e., it is the set
$$
\orbstar(\C,Y)=\set {\pmatrix x\\y\endpmatrix, \pmatrix u\\v\endpmatrix \in \R^2 \times
\R^2}{\ u^2+v^2=1 \ \text{and}\ xv=uy}. 
$$
The set $\orbtwostar(\C,Y)\subset \orbstar (\C,Y)$ is the subset where 
$\pmatrix x\\ y\endpmatrix$ is a non-negative multiple of $\pmatrix
u\\v\endpmatrix$, in other words, it is the set
$$
\orbtwostar(\C,Y)=\set {\pmatrix x\\y\endpmatrix, \pmatrix u\\v\endpmatrix \in \R^2 \times
\R^2}{\ u^2+v^2=1,\ xv=uy \ \text{and}\ ux\ge 0}. 
$$
The mapping $(\R/2\pi \Z) \times [0,\infty)\to \orbtwostar(\C,Y)$ given by
$$
\left(\theta, r\right) \mapsto \left(r\pmatrix \cos \theta\\ \sin \theta\endpmatrix,
\pmatrix \cos \theta\\ \sin \theta\endpmatrix\right) 
$$
is a homeomorphism; so every point of the exceptional divisor $(\R/2\pi \Z) \times \{0\}$
corresponds to a unique prime end, and $\orb(\C,Y) =\orbtwostar(\C,Y)$. You might
prefer to think of the real oriented blow-up as passing to polar coordinates, without
identifying the different angles when the radius is 0. 

We will refer to points of the exceptional divisor $\R/2\pi Z$ as points ``infinitesimally
near the origin'', whose modulus is zero, but which still have a polar angle. 

In this case, the exceptional divisor is the circle $\R/2\pi \Z$, but the identification
with the circle depends on the chart. In general, if you take the real oriented blow-up of
a Riemann surface at a point, the circle $\R/2\pi \Z$  will act freely and transitively
on the exceptional divisor, but there will be no distinguished point on this divisor.  More
generally yet, the real oriented blow-up of a surface along a smooth curve will be a
principle circle-bundle, but not naturally trivial (and usually topologically non-trivial:
see Section VIII).
 
\subheading{Naturality of the real oriented blow-up}

This definition seems (and perhaps it is) a little ad-hoc. Still: it has some nice
properties.  For instance, the real oriented blow-up is always homeomorphic to a
simplicial complex.  Indeed, the space $\orbtwostar(X,Y)$ is obviously subalgebraic,
hence triangulable by a theorem of Hironaka \cite{Hi2}.  Since the prime-end modification
of a simplicial complex along a subcomplex is still a simplicial complex, the result
follows.

Next, notice that it does have a certain weak naturality.

\proposition \naturalityofrealblowup Let $F:X_1 \to X_2$ be an algebraic isomorphism of complex
algebraic spaces, $Y_2 \subset X_2$ an algebraic subspace, and
$Y_1=F^{-1}(Y_2)$.  Then there exists a unique homeomorphism
 $\orb(F):\orb(X_1,Y_1) \to \orb (X_2,Y_2)$ such that the diagram
$$
\CD
\orb(X_1,Y_1)  @>{\orb(F)}>> \orb (X_2,Y_2) \\ 
@V{p_1}VV   @V{p_2}VV     \\ 
X_1    @>F>>  X_2 
\endCD
$$ 
commutes, where the $p_i$ are the canonical projections.
\endstat

\remark The mapping $\orb(F)$ is actually an algebraic isomorphism, after an appropriate
definition has been given.

\proof It is enough to prove the statement locally on $X_1$, so we may suppose that $(Y_2)_\R$ is
given by equations.  But clearly if $p \in I_\R(Y_2)$, then $p\circ F \in I_\R(Y_1)$, and since
$F$ is invertible, this means that $F^*(I_\R(Y_2)) = I_\R(Y_1)$. Then  
Corollary \natblowups\ gives a mapping 
$$
\Cal B_\R(F):\Cal B_\R(X_1,Y_1) \to \Cal B_\R(X_2,Y_2)
$$
which lifts uniquely to the double cover 
$$
\orbstar(F):\orbstar(X_1,Y_1) \to \orbstar(X_2,Y_2)
$$
so that oriented lines are mapped to oriented lines. Thus the region where
$f^*\ge 0$ is preserved and the map lifts to 
$$
\orbtwostar(F):\orbtwostar(X_1,Y_1) \to \orbtwostar(X_2,Y_2).
$$

Finally, we find our lift $\orb (F)$ from the naturality of the prime end modification (see
\natprimeend).
\qed

The authors do not know in what generality it is true that the real oriented blow-up of a
complex manifold along a closed subspace yields a PL manifold with boundary: is
this perhaps always true? In any case, we will require the result only when $X$ is a
complex surface and $Y=D = \cup D_i$ is a divisor consisting of a locally finite union of curves
intersecting transversally, i.e., a divisor with normal crossings. In that case it is
true: the real blow-ups admit the following detailed description.

\theorem {\firstdescofrealblowup} a) The real blow-up $\Cal B_\R(X,D)$ is a 4-dimensional
manifold, and $\orbtwostar(D)$ is a singular 3-dimensional subvariety.

b) The real oriented blow-up $\orb(X,D)$ is a 4-dimensional manifold with boundary, and the
boundary is a manifold with corners.
\endstat

Before we give the proof, note the following corollary.

\corollary \blowupandtubnbhd There is a tubular neighborhood $W$ of $D$ in $X$, together
with a projection $\bar p:\overline W \to D$ extending the identity on $D$, and a
homeomorphism $h:\orb(D) \to \partial W$ making the diagram
$$
\CD
\orb(D)  @>h>> \partial W \\ 
@V p VV   @V {\bar p} VV     \\ 
D    @>>>  D 
\endCD
$$ 
commute.
\endstat

\proofof \blowupandtubnbhd Since $\orb(X,D)$ is a manifold with boundary, there exists
a continuous mapping $\phi:\orb(D) \times [0,1] \to \orb(X,D)$ which is a homeomorphism
onto its image, which forms a tubular neighborhood of the boundary.  Let $\overline W$ be the image
of the projection of
$\phi(\orb(D)
\times [0,1])$ into $X$; it is a neighborhood of $X$ which retracts onto $D$, i.e., a tubular
neighborhood of $D$, and clearly its boundary is $\orb(D) \times \{1\}$, hence
naturally homeomorphic to $\orb(D)$.
\qed

\remarkno \orientrealblowup It follows from \firstdescofrealblowup\ that $\orb(D)$ is an
orientable manifold, being the boundary of an oriented 4-manifold.  However, because we
tend to think of
$\orb(D)$ rather as the boundary of a tubular neighborhood of $D$ than as the boundary of
its complement, we will give $\orb(D)$ the opposite of the boundary orientation of
$\partial
\orb(X,D)$. 

\proofof \firstdescofrealblowup The theorem only needs to be proved locally, so
it is enough to prove it when $X= \C^2$ and $D$ is the union of the coordinate axes. 
However, remember that we are considering these as real algebraic varieties, so $X= \R^4$,
and $D = \R^2 \times \{0\} \cup \{0\} \times \R^2$.  We have already computed equations for
this blow-up, in equations \eqafterzerodiv:
$$
\eqalign{ U_1U_4 &= U_2U_3\cr 
U_4y_1 &= U_2y_2\cr 
U_4x_1 &= U_3x_2\cr 
U_3y_1 &= U_1y_2\cr
U_2x_1 &= U_1x_2.} \peqno \eqafterzerodivagain
$$ 
When we computed these equations, we were working in the complex, but they are just as true
in the real.

In particular, since $\widetilde {\C^4}_{\C^2\times 0\cup 0 \times \C^2}$ is smooth, so is
$\Cal B_\R(X,D)$, and $\Cal B_\R(D)$ is singular, having locally the structure of two
3-dimensional submanifolds intersecting transversally along a (real) surface in $\R^4$ (see
\blowupstillsingeq).

Now, to understand $\orbstar(X,D)$, we simply have to consider the locus with the same
equation, but think of $(U_1,U_2, U_3, U_4)$ as the coordinates of a point in $S^3$
rather than homogeneous coordinates of a point in $\Proj^3_\R$.  Note that
$\pmatrix U_1\\U_2\\U_3\\U_4\endpmatrix $ is a positive section of the tautological (i.e.,
normal) bundle to $S^3$, and $\orbstar(X,D)$ is the locus where
$$ 
\pmatrix x_1y_1\\x_2y_1\\x_1y_2\\x_2y_2 \endpmatrix \quad \text{is a multiple of}\quad
\pmatrix U_1\\U_2\\U_3\\U_4\endpmatrix
$$
whereas $\orbtwostar(X,D)$ is the locus where this multiple is $\ge 0$.

\lemma \parametrizebyVeronese The mapping $\Phi:[0,\infty)^2 \times (\R/2\pi\Z)^2 \to \R^4
\times S^3$ which sends
$$
\bigl(r_1,r_2 \in [0,\infty), \theta_1, \theta_2 \in \R/2\pi\Z\bigr)
\quad \text{to}\quad \left\{\eqalign{
x_1= r_1 \cos \theta_1 \quad&\quad x_2 = r_1 \sin \theta_1 \cr
y_1= r_2 \cos \theta_2 \quad&\quad y_2 = r_2 \sin \theta_2 \cr
U_1= \cos\theta_1\cos\theta_2 \quad&\quad U_2= \sin\theta_1\cos\theta_2 \cr
U_3= \cos\theta_1\sin\theta_2 \quad&\quad U_4= \sin\theta_1\sin\theta_2 }\right\}
$$
maps surjectively to $\orbtwostar(X,D)$, and only the points
$$
(0,0,\theta_1,\theta_2)\quad\text{and}\quad (0,0,\theta_1+\pi,\theta_2+\pi)
$$
are identified.
\endstat

\proofof \parametrizebyVeronese Clearly the equations \eqafterzerodivagain\ are satisfied
on the image of $\Phi$. Moreover, on the image of $\Phi$, we have 
$$ 
\pmatrix x_1y_1\\x_2y_1\\x_1y_2\\x_2y_2 \endpmatrix = r_1r_2\pmatrix
U_1\\U_2\\U_3\\U_4\endpmatrix
$$
and since $r_1r_2\ge 0$, the image of $\Phi$ is contained in $\orbtwostar(X,D)$.

We will try to define an inverse $\Phi^{-1}:\orbtwostar(X,D) \to [0,\infty)^2 \times
(\R/2\pi\Z)^2$; our attempt will show the remainder of the lemma.

Clearly if $(x_1,x_2) \ne (0,0)$ and $(y_1,y_2) \ne (0,0)$, it is easy to reconstruct
$r_1,r_2,\theta_1,\theta_2$ by representing these points in polar coordinates.

You can still compute the polar coordinates if one of $(x_1,x_2)$ or $(y_1,y_2)$ is the origin
of $\R^2$.  Suppose for instance it is the first. The difficulty is then to compute
$\theta_1$, but this can be found from the second lot of equations, for instance $\cos
\theta_1 = U_1/\cos \theta_2$. But when both points are the origin, then the second set
of equations only determine $(\theta_1,\theta_2)$ up to the given ambiguity, since the real
``Veronese'' mapping
$$
V:\R/2\pi\Z\times \R/2\pi\Z \to S^3 \quad,\quad(\theta_1,\theta_2) \mapsto \pmatrix
\cos\theta_1\cos\theta_2\\ \sin\theta_1\cos\theta_2\\
\cos\theta_1\sin\theta_2\\
\sin\theta_1\sin\theta_2 \endpmatrix
$$
is a double cover of its image.
\qedof \parametrizebyVeronese

This shows that the local structure of $\orbtwostar(X,D)$ is the union of two opposite
quadrants, multiplied by a smooth surface, since the domain of $\Phi$ is the first
quadrant times a torus, and  (corner $\times$ torus) is a double cover of its image.  The
open quadrants correspond to the locus
$r_1r_2>0$, and the prime end closure consists exactly of doubling the singular locus,
to make the oriented blow-up a manifold with boundary, the boundary of which is a manifold
with corners.  In particular, the map $\Phi$ lifts to a homeomorphism of the oriented
blow-up.
\qed

We saw that we could understand the real oriented blow-up of $\C$ at $0$ in terms of polar
coordinates.  We will now show that this is still true for the real oriented blow-up of $\C^2$
along $D= \C \times \{0\}\cup \{0\} \times \C$. We will denote $\C^*= \C-\{0\}$. 

\corollary \polardesc The real oriented blow-up 
$\orb(\C^2,D)$ is the closure of the set 
$$
\set {(w_1,w_2,\theta_1, \theta_2) \in (\C^*)^2 \times (\R/2\pi \Z)^2}{\frac
{w_1}{|w_1|}=e^{i\theta_1},\ \frac {w_2}{|w_2|}=e^{i\theta_2}}
$$
in $\C^2 \times (\R/2\pi \Z)^2$.
\endstat

\proof We saw in Lemma \parametrizebyVeronese\ that the real oriented blow-up is parametrized
by $[0,\infty)^2 \times (\R/2\pi \Z)^2$.
\qed 

From this description, we see that in the real oriented blow-up, when one of the variables is
$0$, its polar angle remains.  More precisely, let $p:\orb(D) \to D$ be the projection of the
exceptional divisor.

\corollary {\descoffibersrealblowups} Let $D\subset$ of $X$ be a divisor with normal crossings in a
surface, and
$p:\orb (X, D) \to X$ be the canonical projection. The inverse images of points by $p$ admit the following
canonical description:

a) If $z \in D$ is an ordinary point, $p^{-1}(z)$ is canonically homeomorphic to 
$\Bbb S(T_zX/T_zD)$, where $\Bbb S(E)$, the sphere of $E$, is the space of oriented lines
in
$E$ (i.e., the double cover of $\Proj_\R(E)$).  In particular, it is homeomorphic to a circle, and is
principal under the circle $\R/2\pi\Z$ acting by rotations.  

b) Let $z \in D$ be a double point, say locally $z=D_1 \cap D_2$. Then $p^{-1}(z)$ is
canonically homeomorphic to
$$
\Bbb S\bigl(T_zX/T_zD_1\bigr)\times\Bbb S\bigl(T_zX/T_zD_2\bigr).
$$
In particular, it is homeomorphic to a torus, and principal under
$(\R/2\pi\Z)\times(\R/2\pi\Z)$ acting coordinatewise on the two factors.
\endstat

\corollary \blowupofsmoothcurves If $X$ is a complex surface, and $D \subset X$ is a
smooth curve, then the canonical mapping $p:\orb (D) \to D$ is a principal circle-bundle.
\endstat

\remarksno \comdescofrealblowups 1) We are using the
fact that $T_zX/T_zD$ is a complex vector space of dimension 1 and not just a real vector
space of dimension 2 to make the rotations above well-defined. 

2) In part (b), we could also say that $p^{-1}(z)$ is canonically homeomorphic to 
$\Bbb S(T_zD_2)\times\Bbb S(T_zD_1)$, since there is a canonical homeomorphism
$T_z D_2 \to T_zX/T_zD_1$. This seems simpler, but it doesn't help when we want to think of
a double point as a limit of ordinary points.

3) Let $D'$ denote $D$ with the double points removed. It follows from this theorem
that $p^{-1}(D')$ is a principal circle-bundle over $D'$.

\heading 
VII. The effect of complex blow-ups on real oriented blow-ups
\endheading

Let $z \in D \subset X$ be a point of a divisor in a surface.
We now want to show exactly how a complex blow-up of $X$ at $z$ affects the
real oriented blow-up of the divisor.  This will relate $\orb (X, D)$ and $\orb (\tilde X_z,
\tilde D)$, where $\tilde D=\pi^{-1}(D)$ and $\pi:\tilde X_z \to X$ is the canonical projection. 
Clearly  $\orb (X, D)$ and $\orb (\tilde X_z,
\tilde D)$ have the same overall topology: the boundary of a tubular
neighborhood of the divisor
$D$ is still the boundary of a tubular neighborhood of $\tilde D$.  What is changed is
the way the circle acts.

The description will be given in Theorems
\blowupofordpoint\ and \blowupdoublepoint, which deal with blowing up a smooth point
and a double point of the divisor respectively. Before we state these results, we will
isolate a result in differential topology, which is somehow almost obvious, but which we
found surprisingly hard to prove.  It will be essential to the proof of Theorem
\topofprojlim.

\subheading{Opening Manifolds} In this section, we will show that you can ``open a
manifold along a submanifold'', squeezing the outside arbitrarily little. The difficulty
of the results comes from the fact that the opening is constrained to respect certain
coordinates, whereas the metric with respect to which we are measuring the squeezing is
unrelated to those coordinates.

Let $M$ be a compact smooth $k$-dimensional manifold, $E$ a vector space with a
Euclidean norm, and let $\mu$ be a continuous Riemannian metric on $M \times E$; we will
denote $|\xi|_\mu$ the length of tangent vectors with respect to this metric,
and
$d_\mu(u,v)$ the distance between $u \in M \times E$ and $v \in M \times E$. Note
that we are not assuming any particular relation between the norm $|\ |$ on $E$
and the Riemannian metric $\mu$. 

Let us denote by
$N_r$ the closed subset
$$
N_r= \set{(m,x) \in M \times E}{\ |x|\le r}
$$
of $M \times E$, so that $N_0=M \times \{0\}$, and $N_r$ is a closed tubular
neighborhood of $N_0$ when $r>0$. 

\lemma {\natylemmaontubneighborhoods} For any $\alpha<1$ and any $R>0$, there
exist
$r>0$ and a
$C^\infty$ diffeomorphism
$$
F:(M\times E)-N_0 \to (M\times E)-N_r
$$
such that
\roster
\item
$F$ maps the fiber $\{m\} \times E$ to itself for every $m \in M$;
\vskip 2 pt
\item
$F$ is the identity outside $N_R$;
\vskip 2 pt
\item
$d_\mu(F(u),F(v)) \ge \alpha d_\mu(u,v)$ for all pairs of points $u, v$ in
$M\times E-N_0$.
\endroster
\endstat

We will refer to such a mapping $F$ as ``opening'' $M\times E$ along $M \times \{0\}$.

\proof  Fix $\alpha$, $0<r_0<R$.  Choose a family of increasing
diffeomorphisms
$\tau_\rho:(0,\infty)\to (\rho,\infty)$ for $0< \rho$ sufficiently small, such
that 

(i) All $\tau_\rho(t) = t$ when $t>R$;

(ii) $\tau_\rho$ converges to the identity in the $C^1$ topology on compact subsets of
$(0,\infty)$ as $\rho \to 0$.

(iii) $\tau'_\rho(t)>1$ for $t \le r_0$.

\centerline{\BoxedEPSF{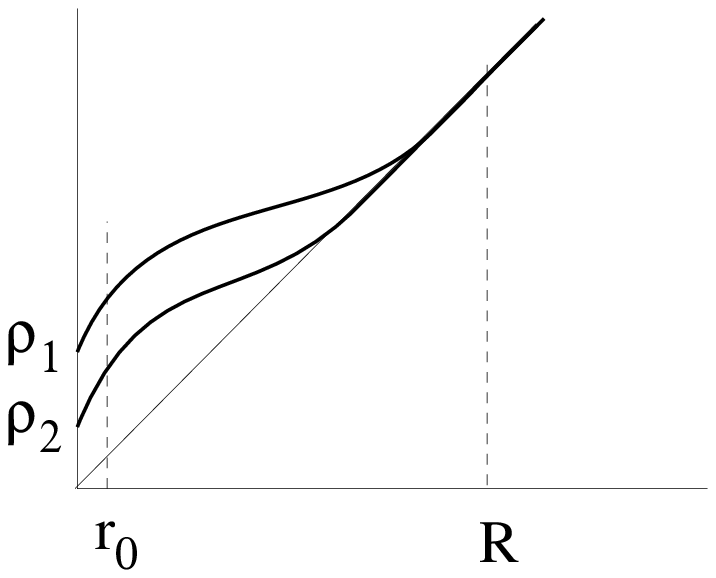 scaled 500}}
\nobreak
\longfigcap \graphoftau {The graphs of $\tau_{\rho_1}$ and $\tau_{\rho_2}$.}

Consider
$$
F_\rho:M\times E -N_0\to M\times E-N_\rho,\quad(m,x) \mapsto 
\left(m,  \tau_\rho(|x|)\frac x{|x|}\right).
$$

We claim that $F_\rho$ will do the trick for $\rho$ sufficiently small.  Only
the third condition poses any problem.  

\sublemma \choosetaunaught  Recall that $\alpha < 1$ and $R>0$ are fixed. There exists $\rho_0>0$
such that
$$d_\mu\bigl(F_\rho(u),F_\rho(v)\bigr) \ge \sqrt\alpha d_\mu(u,v)
$$
when $u,v \in N_{\rho_0}-N_0$ and for all $\rho<\rho_0$. 
\endstat

\proofof \choosetaunaught Choose a number $\beta<1$
such that
$\beta^2>\sqrt \alpha$. For each point  $(m,0)\in M \times \{0\}$, find
$U\subset \R^k, V \subset E$ and a local coordinate
$\phi:U \times V \to M \times E$ at $(m,0)$ such that the diagram
$$
\CD
U\times V  @>>> M\times E \\ 
@V{pr_1}VV   @VV{pr_1}V     \\ 
U   @>>>  M
\endCD
$$
commutes.

By making our coordinate neighborhood small enough, we can assume that the
ratio of the constant (Euclidian) Riemannian metric $\mu_0$ on
$U\times V$, given by the Riemannian metric on $T_{(m,0)}M \times E$, and the
given Riemannian metric $\mu$ (viewed as a metric on $U\times V$ via
$\phi$) satisfy
$$
\beta< \frac {|\xi|_{\mu_0}}{|\xi|_{\mu}} < \frac 1\beta,
$$ 
for all   vectors  $\xi$ tangent to $M \times E$ at a point in $\phi(U \times V)$.

By taking $\rho$ sufficiently small, we can assume that  there exists $V'
\subset V$, such that the mapping
$$
\phi^{-1}\circ F_\rho\circ \phi:U\times V'\to U \times V
$$ is defined and
non-decreasing on tangent vectors, for the metric $\mu_0$, so it will decrease
the length of tangent vectors at most by a factor of $\beta^2$ for the metric
$\mu$. Since $M$ is compact, such coordinate neighborhoods will cover $N_{\rho_0}$
for
$\rho_0>0$ sufficiently small.
\qedof {\choosetaunaught}

Now take $u,v \in M \times E$, and let $\delta$ be a length-minimizing geodesic joining
$F(u)$ to $F(v)$ (which we will temporarily assume exists). Let $\gamma=F^{-1}(\delta)$; of
course, there is no reason to think that $\gamma$ is length-minimizing.

Consider $\gamma_1$ the
part of $\gamma$ in $N_{\rho_0}$ and $\gamma_2$ the remainder. Then the $\mu$-length of
$F_\rho(\gamma_1)$ is contracted at most by $\sqrt\alpha$ when $\rho<\rho_0$,
whereas the
$\mu$-length of $F_\tau(\gamma_2)$ is contracted arbitrarily little when $\rho$
is sufficiently small, since on $M\times E-M_{\rho_0}$ the mappings $F_\rho$
converge to the identity in the $C^1$ topology. Thus for any $\epsilon$, we may assume that
$\gamma$ is contracted by $\sqrt \alpha-\epsilon$, which we may take to be $>\alpha$ by
taking $\epsilon$ sufficiently small.  Then we have
$$
d_\mu \bigl(F(u),F(v)\bigr) = l_\mu(\delta) \ge \alpha l_\mu(\gamma) \ge \alpha d_\mu(u,v).
$$

\centerline{\BoxedEPSF{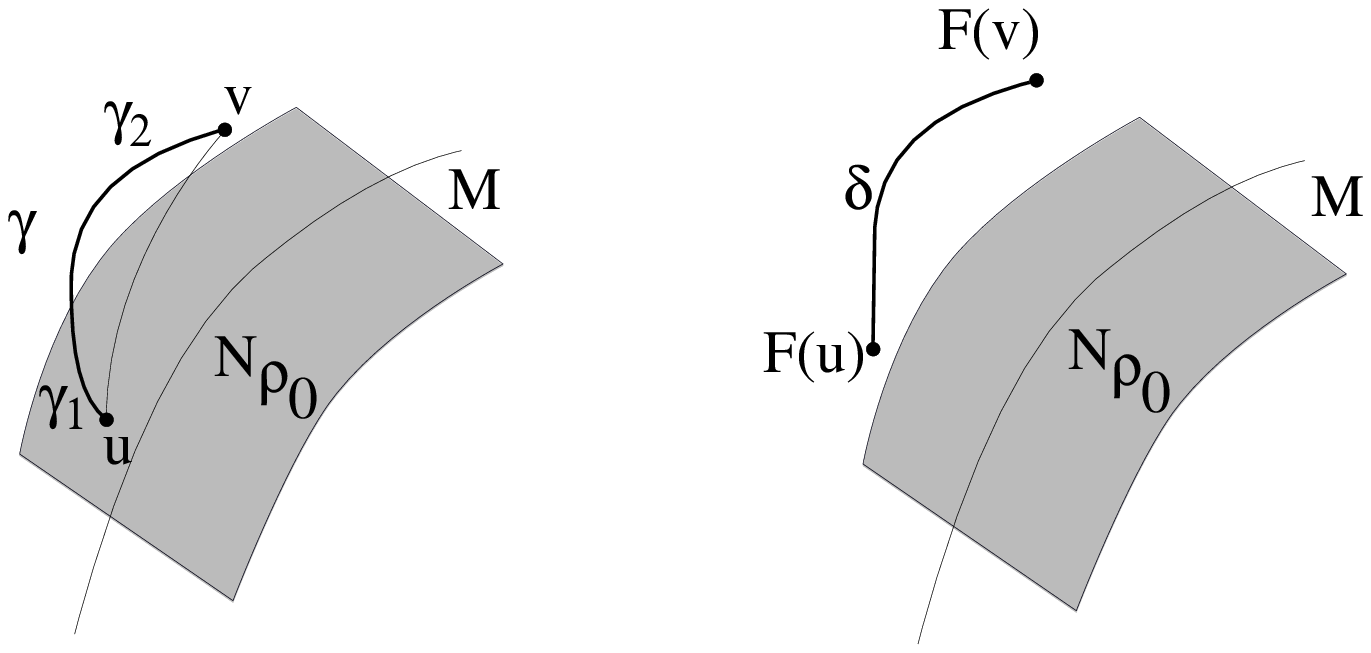 scaled 600}}

\longfigcap \openmanfig {A picture of $\gamma = \gamma_1 \cup \gamma_2$, together with the
true shorter geodesic joining $u$ to $v$, and the image $\delta = F(\gamma)$.}

If geodesics do not exist, it is easy to modify the argument above, using curves which
almost minimize distances.
\qed

\subheading{Blowing up smooth points of divisors}  Let $D \subset X$ be a divisor with
normal crossings in a surface, and $\orb (D) \subset \orb(X,D)$ be the
real oriented blow-up, so that the diagram
$$
\CD
\orb (D)  @>>> \orb(X,D) \\ 
@VpVV   @VpVV     \\ 
D   @>>>  X 
\endCD
$$
commutes.  Let $z\in D$ be a point, $\tilde X_z$ be the complex blow-up of $X$ at $z$,  $\pi:\tilde
X_z\to X$ the canonical projection, and set $\tilde D = \pi^{-1}(D)$.
The divisor
$\tilde D$ is the union of the exceptional divisor $E$ and the proper transform $D'$ of $D$; we will
call their point of intersection $\tilde z = E \cap D'$.

 Let $w_1, \ w_2$ be local coordinates on $X$ at $z$, mapping a neighborhood $W$ of $z$
isomorphically onto the region $V_R \times V_R$ of $\C^2$, where $V_R$ is the open disc of radius
$R$ in $\C$. To lighten notation, we set $V_R^*=V_R-\{0\}$.
Within $W$, we will assume that $D$ is given by the equation $w_2=0$, so that $D$ has no double
points in $W$. Note that $\tilde W_z=\pi^{-1}(W)$ can be understood either as the inverse image of
$W$ in $\tilde X_z$, or as $W$ blown-up at $z$.

Recall (see Example \blowuppointonsurf) that  $\tilde W_z \subset
W\times\Proj^1_\C$ is defined by the equation
$w_1U_2=w_2U_1$, where $(w_1,w_2)\in W$ and $\bvec{U_1}{U_2}$ are homogeneous coordinates in
$\Proj^1_{\C}$.  We will also use the affine coordinates $u_1=U_1/U_2$ and $u_2=U_2/U_1$ on
appropriate parts of $\tilde W_z$; note that $w_1$ and $u_2$ are local coordinates near $\tilde z$.

Denote by $\tilde p:\orb(\tilde X_z, \tilde D) \to \tilde X_z$ the canonical
projection from the real oriented blow-up, and use the same name for the restriction to $\orb(\tilde
D)$ to the boundary.

We will need to parametrize $\orb(W,D)$ and $\orb(\tilde W_z, \tilde D)$.  By Corollary \polardesc\
the choice of
$w_1, \ w_2$ induces a homeomorphism of
$p^{-1}(W)$ with the closure of the set 
$$
\set{(w_1,w_2,\theta_2) \in V_R \times V_R^* \times \R/2\pi \Z\ }{\frac{w_2}{|w_2|}=e^{i\theta_2}}
\peqno \paramdesceq 
$$
in $V_R\times V_R \times \R/2\pi \Z$. In this closure, when $w_2=0$, its polar angle remains, so
$(w_1, |w_2|, \theta_2)$ parametrize $\orb(W,D)$.

By \polardesc\ again, a neighborhood of $\tilde p^{-1}(\tilde z)$ is homeomorphic to
the closure of the set of 
$$
\set {(w_1,u_2, \theta_1, \phi_2) \in V_R^* \times \C^*\times (R/2\pi \Z)^2}{\theta_1 = \frac
{w_1}{|w_1|},\ \phi_2 = \frac{u_2}{|u_2|}}
$$
in $V_R \times \C \times (R/2\pi \Z)^2$. So $(|w_1|, |u_2|, \theta_1, \phi_2)$ parametrize a
neighborhood of $\tilde z$ in $\orb(\tilde W_z, \tilde D)$; in particular, by Corollary
\descoffibersrealblowups, the torus
$\tilde p^{-1}(\tilde z)$ is parametrized by $\theta_1$ and $\phi_2$.

\proposition \divrightasideal The mapping $\pi:\tilde X_z \to X$ lifts to a unique mapping $\tilde
\pi:\orb(\tilde X_z, \tilde D) \to \orb(X,D)$ such that the diagram
$$
\CD
\orb (\tilde X_z,\tilde D)  @>{\tilde \pi}>> \orb(X,D) \\ 
@V{\tilde p}VV   @V{p}VV     \\ 
\tilde X_z  @>{\pi}>>  X 
\endCD
$$
commutes. 
\endstat

\proof It is enough to extend $\pi$ above $W$. Let us compute $\tilde \pi$ near $\tilde z$. The space
$\tilde W_z$ is parametrized by
$w_1$ and $u_2$, and $\pi(w_1,u_2)=(w_1, w_1u_2)$. 

Thus if we set $w_1=r_1e^{i\theta_1},\ u_2=\rho_2e^{i\phi_2}$, then the map $\pi$ is 
$$
\pi:\pvec {r_1e^{i\theta_1}}{\rho_2e^{i\phi_2}} \mapsto 
\pvec {r_1e^{i\theta_1}}{r_1\rho_2e^{i\phi_2+\theta_1}}.
$$
As a map of real oriented blow-ups, the domain is parametrized by $(r_1, \theta_1, \rho_2, \phi_2)$
and the range by $(w_1, |w_2|\in [0,\infty), \arg w_2\in \R/2\pi \Z)$, and the map is written
$$
(r_1, \theta_1, \rho_2, \phi_2) \mapsto (r_1e^{i\theta_1}, r_1\rho_2, \theta_1+\phi_2) \peqno
\extendtoblupordeq
$$
and is obviously continuous. The proof near the one omitted point $u_2=\infty$ is easier.
\qed

\proposition \compsimpblowupontorus The mapping $\tilde \pi$ maps the torus $\tilde p^{-1}(\tilde
z)$( parametrized by $(\theta_1, \phi_2)$)  to the circle $p^{-1}(z)$ (parametrized by $\theta_2$),
by the mapping
$$ (\theta_1, \phi_2) \mapsto \theta_1+ \phi_2.
$$
\endstat

\proof Since the projection $\pi$ is simply $(w_1, u_2) \mapsto (w_1, w_1 u_2)$, and our variable
$\theta_2$ corresponds to the argument of the second entry, the result follows.
\qed

We will now show how the real oriented blow-up $\orb(\tilde X_z, \tilde D)$ can be constructed from
$\orb (X,D)$. Recall that $\orb(X,D)$ is a 4-dimensional manifold with boundary, and that the
boundary is a 3-manifold containing tori, corresponding to the double points of $D$, and with a
circle action on the complement of these tori.  We need to understand the corresponding structure on
$\orb(\tilde X_z,\tilde D)$.
                                                                                                            
Roughly speaking, it can be understood as follows. If $z\in D$ is a smooth point of $D$,
then the fiber $p^{-1}(z)$ is a circle. Thicken this circle to make a solid torus,
invariant under the existing circle action. Keep the old circle action on the outside of the solid
torus, and modify it inside, so that the oriented circle orbits on the boundary of the solid
torus are the ``sums'' of the old ones and of boundaries of discs $\Delta$ in the solid torus,
which are oriented so that
$p:\Delta\to D$ is orientation preserving.

This description is too imprecise to be proved in this form: Theorem
\blowupofordpoint\ makes it precise.

Call $Z$ the space $\orb(X,D)$, but with a modified circle action on part of the boundary $\orb(D)$.
Choose $r<R$, set $\Delta_r \subset D$ the region $|w_1| < r$, and in $p^{-1}(\Delta_r)$ use
the circle action 
$$
\Theta*(w_1,0,\theta_2) = (w_1 e^{i\Theta},0,\theta_2+\Theta), \peqno \circactZeq
$$
where we are using the coordinates (\paramdesceq) in $\orb(X,D)$.

\theorem \blowupofordpoint There exists a homeomorphism $h:\orb (\tilde X_z, \tilde D) \to Z$,
which coincides with $\tilde \pi$ outside of $\tilde p^{-1}(\tilde W_z)$, and which respects the
circle actions. \endstat

\proof  We only need to construct $h$ in $\tilde p^{-1}(\tilde W_z) \subset \orb(\tilde X_z,
\tilde D)$, so long as $h$ coincides with $\tilde \pi$ outside a compact subset of $\tilde
p^{-1}(\tilde W_z)$. We will begin by defining our homeomorphism on $\partial \orb(\tilde X_z, \tilde
D) = \orb (\tilde D)$, and then extend it to the interior.

  A neighborhood of $\tilde
p^{-1}(E-\{\tilde z\})$ can be described as the closure of the set of
$$
\set {(u_1,w_2,  \theta_2) \in \C \times V_R^* \times \R/2\pi \Z}
{\ e^{i\theta_2}=\frac{w_2}{|w_2|}}
$$
in $\C \times V_R\times \R/2\pi \Z$, and $\tilde p^{-1}(E)$ is the subset of
equation $w_2=0$.  Again the polar angle of $w_2$ remains when $w_2=0$;
we will write
$u_1=|u_1|e^{i\phi_1}$, but when $u_1=0$ it does not have an argument. 

Choose a homeomorphism $\eta_1:[0,\infty) \to [r,\infty)$, such that $\eta_1(t)=t$
for $t>R$; it will need to satisfy further properties in order for Theorem 
\controlledblowupordpt\ to be true, which we will spell out when we come to it.
Let us set (see Figure \graphsetas)
$$
\eta_2(t)=\frac {tr}{t+r}.
$$

\centerline{\BoxedEPSF{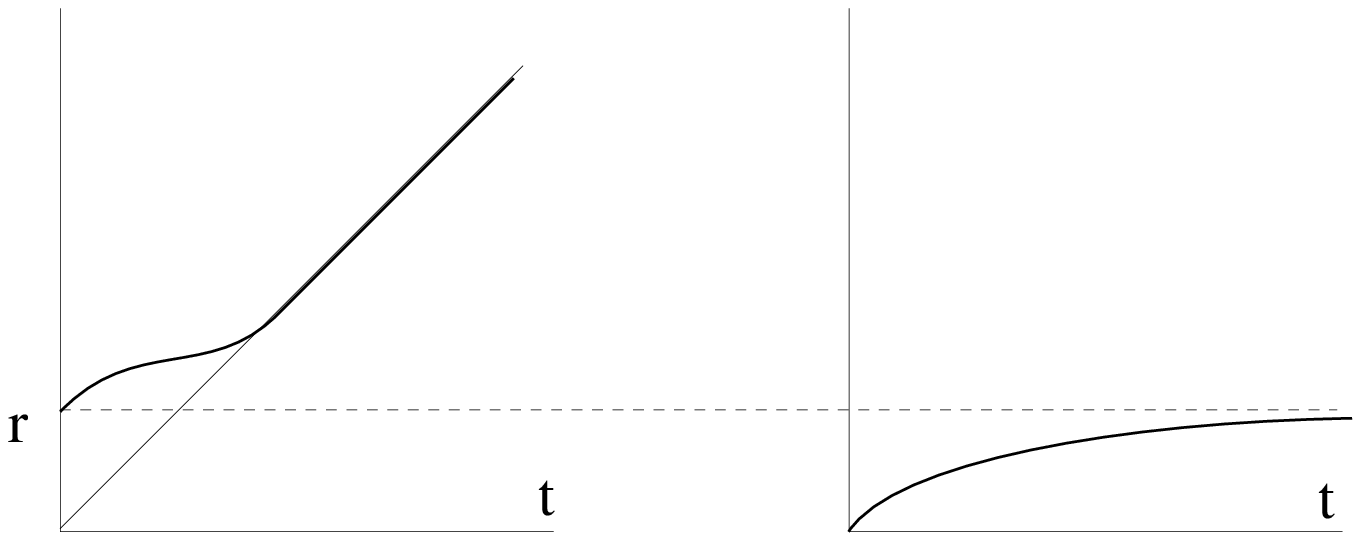 scaled 600}}

\longfigcap \graphsetas {The graph of $\eta_1$ (left) and the graph of $\eta_2$ (right).}

Define $h:\tilde p^{-1}(D'-\{\tilde z\}) \to \orb(D)$ by the formula
$$
h:(w_1, 0, \theta_1, \theta_2) \mapsto \bigl(\eta_1(|w_1|)e^{i\theta_1}, 0, \theta_2\bigr).
\peqno \firstwriteofh
$$

As a set, $\partial
Z=\orb(D)$.  Using parametrization (\paramdesceq) in the range, we define $h:\tilde
p^{-1}(E-\{\tilde z\})
\to \partial Z$ by the formula 
$$
h:(u_1, 0, \phi_1, \theta_2) \mapsto \bigl(\eta_2(|u_1|)e^{i(\phi_1+\theta_2)}, 0,
\theta_2\bigr). \peqno \secondwriteofh
$$

To check the continuity of these two definitions on $\tilde p^{-1}(\tilde z)$, we choose local
coordinates $(u,v)$ on $\orb(\tilde X_z, \tilde D)$ near $\tilde z$, related to the previous
coordinates by
$$\eqalign{
\left\{\eqalign{u&=w_1\cr v&=\frac {w_2}{w_1}\ \ }\right.\ &\qquad\text{on a neighborhood of
$\tilde p^{-1}(D'-\{\tilde z\})$}\cr 
\left\{\eqalign{u&=w_2u_1\cr v&=\frac {1}{u_1}}\right.\ &\qquad\text{on a neighborhood of
$\tilde p^{-1}(E-\{\tilde z\})$}} $$

By Corollary \polardesc, in the coordinates $u,v$, the space $\orb(\tilde X_z, \tilde
D)$ is canonically homeomorphic to the closure of the set 
$$
\set {(u,v, \beta,\gamma) \in V_R^*\times \C^* \times  (\R/2\pi \Z)^2}{e^{i\beta}=
\frac u{|u|},\  e^{i\gamma}= \frac v{|v|}}
$$
in $V_R\times \C \times  (\R/2\pi \Z)^2$; when the variables $u$ or $v$ or both vanish, their
polar angles both remain. The locus $\tilde p^{-1}(\tilde D) \cap \tilde W_z$ is given by the
equation
$uv=0$.

The expression (\firstwriteofh) of $h$ on $p^{-1}(D'-\{\tilde z\})$ becomes, in these
coordinates,
$$
h: (u,0,\beta, \gamma) \mapsto (\eta_1(|u|) e^{i\beta}, 0, \beta+\gamma).
$$
Similarly, the expression (\secondwriteofh) of $h$ on $\tilde p^{-1}(E-\{\tilde z\})$ becomes
$$
h: (0,v,\beta, \gamma) \mapsto \left(\eta_2\Bigl(\frac 1{|v|}\Bigr)
e^{-i\gamma}e^{i(\beta+\gamma)}, 0, \beta+\gamma\right)= \left(\eta_2\Bigl(\frac
1{|v|}\Bigr) e^{i\beta}, 0, \beta+\gamma\right). 
$$
We have used the computation
$$
\theta_2= \arg w_2 = \arg uv = \beta+\gamma
$$
of $\theta_2$ in the coordinates $u,v$.

Now to check the continuity:
$$\eqalign{
\lim_{\rho \to 0} h(\rho e^{i\beta}, 0, \beta, \gamma) = (re^{i\beta}, 0, \beta+\gamma)\cr  
\lim_{\rho \to 0} h(0,\rho e^{i\gamma}, \beta, \gamma) = (re^{i\beta},0,\beta+\gamma).}
$$
It works; our homeomorphism is well-defined on $\orb (\tilde D)=\partial \orb(\tilde X_z,
\tilde D)$.

We need to check that $h$ transforms the circle action of $\orb(\tilde D)$ into the circle action of
$Z$.  The only place that needs checking is in $\tilde p^{-1}(E)$. Using the notation of
Equation (\secondwriteofh), we see that in the domain the circle action is
$$
\Theta *(u_1, 0, \phi_1, \theta_2)= (u_1, 0, \phi_1, \theta_2+\Theta);
$$
by (\circactZeq), the circle action in the range is 
$$
\Theta*\Bigl(\eta_2(|u_1|)e^{i(\phi_1+\theta_2)}, 0, \theta_2\Bigr) =
\Big(\eta_2(|u_1|)e^{i(\phi_1+\theta_2+\Theta)}, 0,
\theta_2+\Theta\Bigr).
$$
Clearly $h$ transforms one into the other.

%

Now we extend $h$ to the interior.  Let us consider the subset $ W' \subset \tilde W_z$
given by the inequality
$$
 W'= \set {(w_1,w_2)}{|w_1|^2+(|w_2|-r)^2<r^2}, \quad 2r < R,
$$
as shown in Figure \picofW.

\centerline{\BoxedEPSF {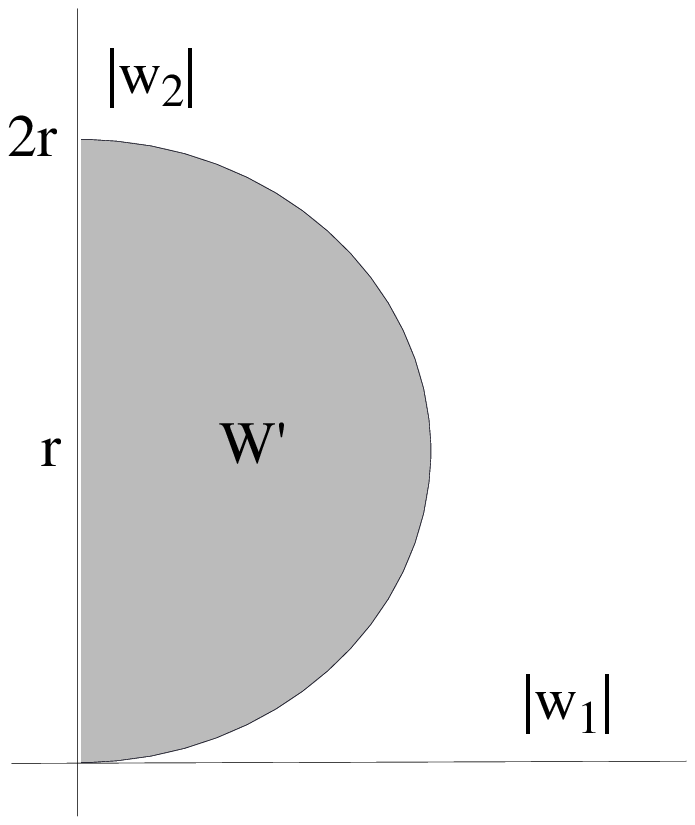 scaled 600}}
\figcap \picofW {The region in the $|w_1|,\ |w_2|$-plane defining $W'$.} 

The set $W'$ can be considered as a subset of $W$, or of $\tilde W_z$, or of $\orb(W,D)$, or of
$\orb(\tilde W_z, \tilde D)$, since the canonical projections are all isomorphisms on the inverse
image of $W'$. As a subset of
$\tilde W_z$, it is a neighborhood of
$E-\{\tilde z\}$, and we can use
$w_1/w_2, w_2$ as parameters in $W'$, which extend the parameters $u_1,w_2$ above. As a subset
of $\orb(\tilde W_z, \tilde D)$, it is a neighborhood of $\tilde p^{-1}(E-\{\tilde z\})$.

We can now extend $h$ by the formulas
$$
h:\cases
(w_1,w_2) \mapsto (\eta_{1, |w_2|}(|w_1|) \frac {w_1}{|w_1|}, w_2) &\quad\text{if $(w_1,w_2)
\notin W'$}\\
(u_1, w_2) \mapsto (\eta_2(|u_1|) \frac {u_1}{|u_1|}\frac {w_2}{|w_2|}, w_2) &\quad\text{if
$(u_1,w_2) \in W'$},
\endcases \peqno \defhoninsideeq
$$
where we are using the coordinates $(w_1,w_2)$ in the range in both cases.  
The second line above could also be written
$$
(w_1, w_2) \mapsto (\eta_2\left(\frac{|w_1|}{|w_2|}\right) \frac {w_1}{|w_1|}, w_2).
$$
We opted to write it in terms of $u_1$ instead because that variable extends to the boundary. Note
that since
$W'$ does not intersect the boundary, we do not need the polar angles as variables in either domain
or range; of course, they will return when we check continuity.

The functions $\eta_{1,\rho}$ and $\eta_2$ are illustrated in Figure \graphofetasinside.

We need to check the continuity of the mapping, both on the  part of $\partial W'$
where $w_2 \ne 0$, and on $\partial Z$.

\centerline{\BoxedEPSF{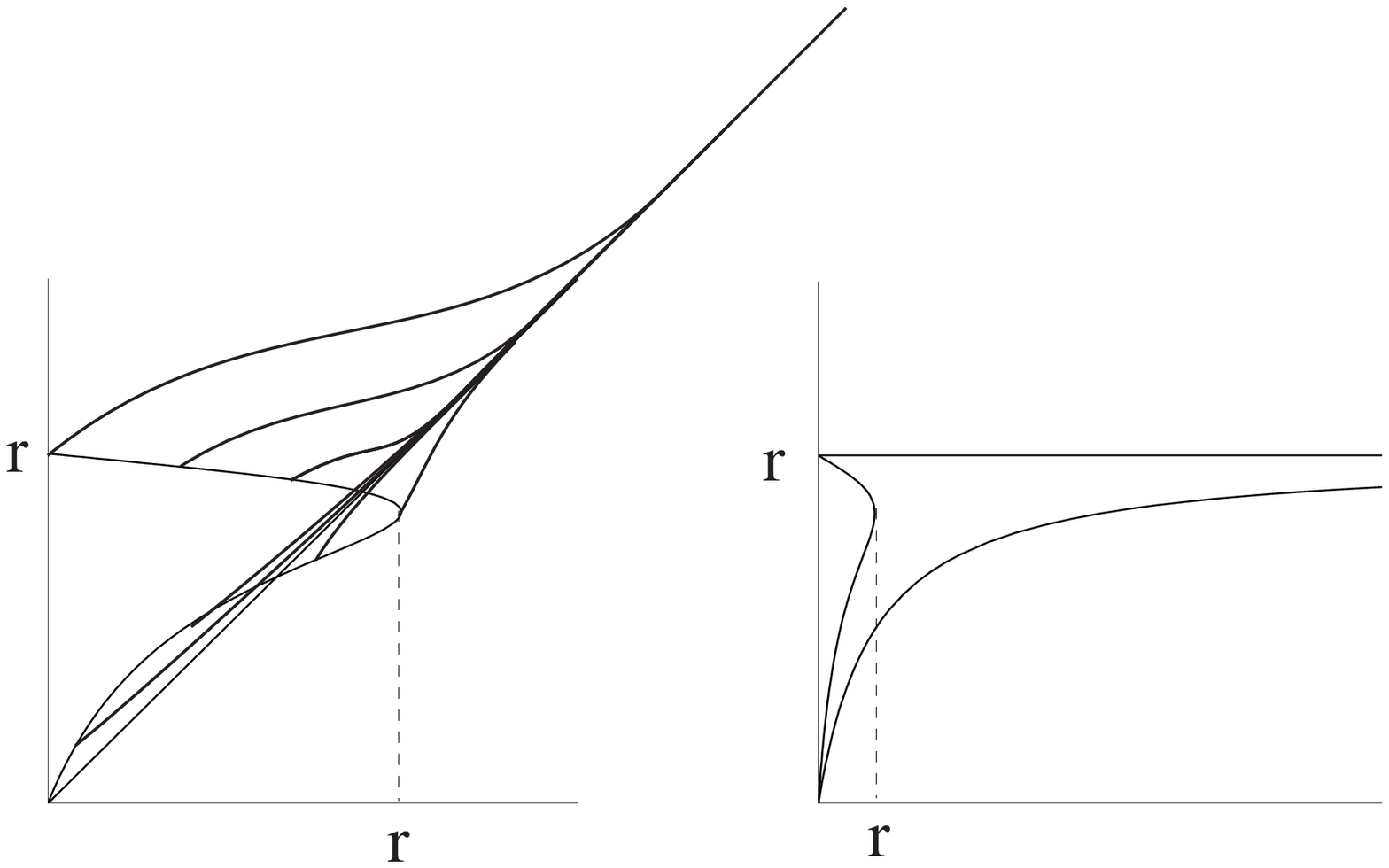 scaled 400}}
\nobreak
\longfigcap \graphofetasinside {On the left, the graphs of functions $\eta_{1,\rho}$ for
various values of $\rho \in [0, 2r]$.  They are defined for $t \ge
\sqrt{r^2-(\rho-r)^2}$, and the domain is mapped homeomorphically to
$$
\bigl[\eta_2(\sqrt{r^2-(\rho-r)^2}/\rho),\infty\bigr);
$$
thus all the graphs start on the curve given
parametrically by  
$$\rho \mapsto
\left(\sqrt{r^2-(\rho-r)^2},\,\eta_2(\sqrt{r^2-(\rho-r)^2}/\rho)\right)
$$
which is also represented on the right, where the horizontal scale and the vertical scale are
different.}

The second will be true if $\eta_{1,0}(t)=\eta_1(t) $, as is clear from Equations
(\defhoninsideeq) and (\secondwriteofh), since the argument of $w_2$ remains as a variable when
$w_2=0$. For the first, we have
$|w_1|^2 = r^2-(|w_2|-r)^2
$, and we see that we must have $$
\left(\eta_{1,|w_2|}(|w_1|) \frac{w_1}{|w_1|},w_2\right)=
\left(\eta_2\left(\left|\frac{w_1}{w_2}\right|\right)\frac {w_1}{|w_1|},w_2\right) 
$$
when $|w_1|^2=r^2-(|w_2|-r)^2$.  Setting $\rho=|w_2|$ to clarify the notation, we need
$$
\eta_{1,\rho}(\sqrt{r^2-(\rho-r)^2}) = \eta_2\left(\frac{\sqrt{r^2-(\rho-r)^2}}{\rho}\right).
\peqno \firstcompcondeq $$

It is now easy to choose $\eta_{1,\rho}(t)$, defined for $t\ge \sqrt{r^2-(\rho-r)^2}$, so
that 
\roster
\item Equation (\firstcompcondeq) is satisfied;
\item for each $\rho$, the function $\eta_{1,\rho}(t)$ is monotone increasing;
\item $\eta_{1,0}(t)=\eta_1$;
\item $\eta_{1,2r}(t) = t$.
\endroster
These conditions will guarantee that the mapping $h$ defined by (\defhoninsideeq) is a
homeomorphism as required.
\qed

To prove Theorem \topofprojlim, we will need some control on how the homeomorphism
$h:\orb(\tilde X_z, \tilde D) \to Z$, constructed above, distorts distances. Note that $h$ is not
canonical: it depends on the local coordinates $w_1,w_2$, and on the choice of $r$, $\eta_1$
and $\eta_2$.

\theorem \controlledblowupordpt If $\mu$ is any
Riemannian metric on $\orb(D)$ which is continuous in the $(w_1,\theta_2)$ coordinates, then
for any $\epsilon>0$, $\alpha<1$ and any neighborhood $W$ of $z$ in $X$,  the homeomorphism
$h$ can be chosen to coincide with $\tilde \pi$ except on $\orb\left(\tilde W_z, 
\tilde D \cap \tilde W_z \right)$, and such that for any $x,y \in \orb(\tilde D)$, 
 $$ 
d_\mu(h(x),\tilde \pi(x)) <\epsilon 
$$
and
$$
d_\mu\bigl(h(x), h(y)\bigr)\ge \alpha d_\mu\bigl(\tilde \pi(x),\tilde \pi(y)\bigr) .
$$
\endstat

Note that all the distances are being measured in $\orb(D)=\partial Z$ with respect to the
Riemannian metric $\mu$.

\proof  This is a matter of carefully choosing the number $r$ and the function $\eta_1$.
Clearly $d_\mu(h(x),\tilde \pi(x))$ will be arbitrarily small if $r$ is chosen sufficiently
small.  The second inequality follows from Lemma \natylemmaontubneighborhoods, on the outside
of $\tilde \pi^{-1}(E)$, we have precisely a map of the form specified there, and we have
proved  the needed result. 

With our choice of coordinates, the region $p^{-1}(\Delta_R)$ is canonically $V_R \times S^1$,
parametrized by $w_1$ and $\theta_2$, with $|w_1| \le R$. We can modify the Riemannian metric
$\mu$ on this region by setting $\mu_1(w_1,\theta_2)=\mu(0,\theta_2)$.  With respect to $\mu_1$,
the projection
$\bar \pi:(w_1,\theta_2) \mapsto (0,\theta_2)$ is obviously length-decreasing.  But since $\mu$
and
$\mu_1$ coincide on $p^{-1}(z)$,  we have
$$
\alpha\le \left|\frac{\mu}{\mu_1}\right| \le \frac 1\alpha
$$
on $p^{-1}(\Delta_r)$ for $r$ sufficiently small (we only need the inequality on the right). So
the projection
$\bar \pi$ satisfies $\alpha d_\mu(\bar \pi(u),\bar \pi(v))\le d_\mu(u,v)$ for $u,v \in
p^{-1}(\Delta_r)$.

Now take the value of $r$ in the definition of $h$ to be the $r$ above. Since $h$
conjugates $\tilde \pi$ with $\bar \pi$, we see that the second requirement is satisfied
when $x,y \in \tilde p^{-1}(E)$, so that $h(x), h(y) \in p^{-1}(\Delta_r)$.   

If $x \in p^{-1}(E)$ and $y \notin p^{-1}(E)$, it suffices to put the two arguments above
together.  Join $h(x)$ to $h(y)$ by a length-minimizing geodesic $\gamma$, and consider the
parts $\gamma_1 = \gamma \cap p^{-1}(\Delta_r)$ and $\gamma_2=\gamma - \gamma_1$. Then we have
the string of inequalities
$$
\eqalign{d_\mu\bigl(h(x),h(y)\bigr) &= l_\mu(\gamma_1)+ l_\mu(\gamma_1) \ge
\alpha\Bigl(l_\mu\bigl(\tilde \pi \circ
h^{-1}(\gamma_1)\bigr)\Bigr)+\alpha\Bigl(l_\mu\bigl(\tilde \pi \circ
h^{-1}(\gamma_2)\bigr)\Bigr)\cr
&\ge
\alpha_\mu\bigl(\tilde
\pi(x), \tilde \pi(y)\bigr).}
$$
\qed 

\subheading{Blowing up double points of divisors}  We now will give the corresponding
description for double points of divisors. Let 
 $$
D \subset X,\ \orb(D) \subset \orb(X,D), \ p:\orb(X,D) \to X
$$
be as above.  Let $z \in D$ be a double point, and denote by
$$
\pi:\tilde X_z\to X,\  \tilde p:\orb(\tilde X_z, \tilde D) \to \tilde X_z
$$
the canonical projections (complex and real respectively).

Let $w_1$ and $w_2$ be local coordinates centered at $z$, giving an isomorphism of a
neighborhood $W$ of $z$ with the bidisc $V_R \times V_R$ so that
$D\cap W$ is given by the equation $w_1w_2=0$; in particular, in $W$ there is no other
double point of $D$. Let $D_1 = W \cap \{w_2=0\},\ D_2=W\cap \{w_1=0\}$ be the two irreducible
components of $D$ in $W$. 

The space  $\tilde W_z$ is the locus in $W \times \Proj^1$ given by the
equation $w_1U_2=w_2U_1$, where $\bvec {U_1}{U_2}$ are homogeneous coordinates; again we will use
the affine coordinates $u_1=U_1/U_2$ and $u_2=U_2/U_1$ where they are defined. The divisor $\tilde
D= \pi^{-1}(D)$ is the union of the proper transforms $D_1'$ of $D_1$ and $D_2'$ of $D_2$, and the
exceptional divisor $E$, as shown in Figure \threelines.
\medskip
\centerline{\BoxedEPSF{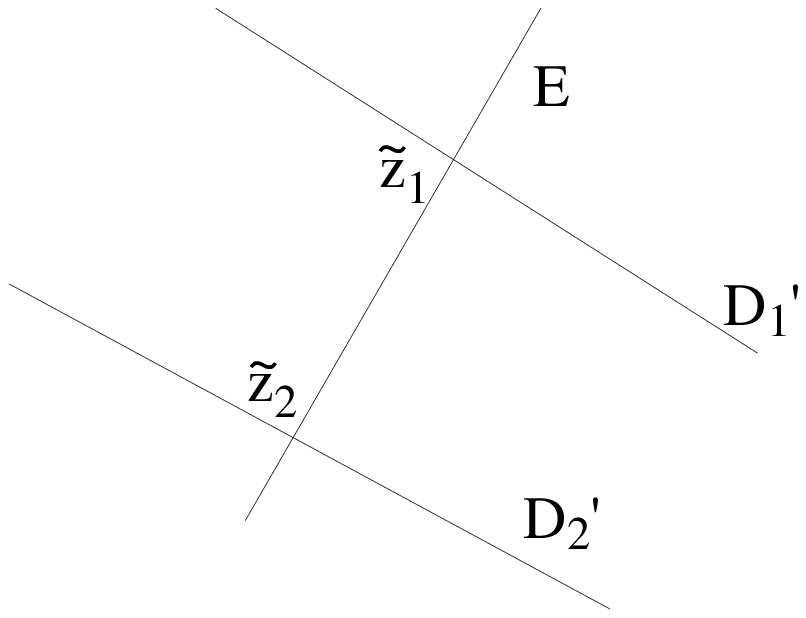 scaled 600}}
\figcap \threelines {The configuration in $\tilde W_z$.}

By Corollary \polardesc, we can identify
$\orb(W,D \cap W )=p^{-1}(W)$ with the closure of the set 
$$
\set {(w_1,w_2,\theta_1,\theta_2)\in (V_R^*)^2 \times (\R/2\pi \Z)^2}{e^{i\theta_1}=\frac
{w_1}{|w_1|},\ e^{i\theta_2}=\frac {w_2}{|w_2|}
}$$
in  
$$
V_R^2 \times (\R/2\pi \Z)^2
$$
so that when each coordinate function vanishes, its polar angle remains. This is a manifold
with boundary, where the boundary is defined by $w_1w_2=0$.

The fiber  $p^{-1}(z)$ is a torus, which is canonically principal under the group $\R/2\pi\Z \times
\R/2\pi \Z$, and which is parametrized by $\theta_1, \theta_2$.

Similarly, a neighborhood of $\tilde z_1\in \orb(\tilde W_z, \tilde D \cap \tilde W_z)$ is
parametrized by 
$$
|u_2| \in [0,\infty),\ \phi_2= \arg u_2, |w_1|\  \text{and}\ \theta_1=\arg w_1,
$$
with the torus $\tilde p^{-1}(\tilde z_1)$ parametrized by $\phi_2$ and $\theta_1$. The torus
$\tilde p^{-1}(\tilde z_2)$ can be parametrized in a similar way.

Proposition \pitildeexistsfordblepoint\ is very similar to Proposition \divrightasideal.

\proposition \pitildeexistsfordblepoint There exists a unique lift $\tilde \pi:\orb(\tilde
X_z,\tilde D) \to \orb (X,D)$ such that the diagram
$$
\CD
\orb (\tilde X_z,\tilde D)  @>{\tilde \pi}>> \orb(X,D) \\ 
@V{\tilde p}VV   @V{p}VV     \\ 
\tilde X_z  @>{\pi}>>  X 
\endCD
$$
commutes. 
\endstat

\proof Begin by working near $\tilde z_1$. The computation starts as in Proposition \divrightasideal,
but in the range of equation  (\extendtoblupordeq), we must use $(r_1, \theta_1, r_2,
\theta_2)$  as coordinates, since  the arguments of $w_1$ and $w_2$ both remain when $r_1$ or $r_2$
tend to $0$.  So the  mapping $\tilde \pi$ in these coordinates becomes (see Equation
(\extendtoblupordeq))
$$
(r_1, \theta_1, \rho_2, \phi_2) \mapsto (r_1, \theta_1, r_1\rho_2, \theta_1+\phi_2);
$$
again this is obviously continuous.  The computation near $\tilde z_2$ is similar.
\qed

\proposition \compmapontorofdblept The mapping $\tilde p$, mapping the torus $\tilde
p^{-1}(z_1)$, parametrized by $\theta_1, \phi_2$, to the torus $p^{-1}(z)$, parametrized by
$\theta_1$ and
$\theta_2$, is given by the formula
$$
\bvec {\theta_1}{\phi_2} \mapsto \bvec {\theta_1}{\theta_1+\phi_2}.
$$
\endstat

\proof Near $\tilde z_1$, we can use the local coordinates $w_1$ and $u_2$; again the
mapping $\pi$ is given by the formula
$$
(w_1, u_2) \mapsto (w_1, u_2 w_1).
$$
Thus the argument of the first coordinate is $\theta_1$, and the argument of the second is
$\theta_1 + \phi_2$.
\qed

As in the case of simple points, we will need a much more precise description of the real oriented
blow-up
$\orb(\tilde X_z, \tilde D)$: it can be constructed from
$\orb(X,D)$ as follows. Thicken the torus $p^{-1}(z)$, so that the thickened torus
is invariant under the circle actions. Keep the old
circle actions on the outside of the thickened torus, and modify it inside, so that the
oriented circle orbits on the boundary torus are the ``sums'' of an orbit of $\R/2\pi\Z\times
\{0\}$ and an orbit of $\{0\}\times \R/2\pi \Z$.

Again, this description is too imprecise to be proved; Theorem \blowupdoublepoint\
should make everything precise.

Let $Z$ be $\orb(X,D)$, with a new circle action the boundary $\partial Z = \orb (D)$ defined as
follows. Let $P$ be the region 
$$
P= \{|w_1|\le r, w_2=0\} \cup \{w_1=0, |w_2|\le r\}
$$
and on this region, define the circle action
$$\eqalign{
\Theta *(w_1,0,\theta_1,\theta_2) = (w_1,0,\theta_1+\Theta,\theta_2+\Theta)\cr
\Theta *(0,w_2,\theta_1,\theta_2) = (0,w_2,\theta_1+\Theta,\theta_2+\Theta).}
$$

Keep the previous circle action outside $P$; this gives as it should two circle actions on the
two tori $|w_1|=r$ and $|w_2|=r$.

\theorem \blowupdoublepoint  There exists a homeomorphism of $h: \orb (\tilde X_z, \tilde D)\to
Z$ which respects the circle actions on the boundaries.
\endstat

\proof The proof is quite similar to that of \blowupofordpoint, and we will be a bit sketchier. As
before, we will begin by defining the map $h$ on the boundary, and then extend it to the interior. 
The divisor
$\tilde D\cap \tilde W_z$ consists of three components: the proper transforms $D_1', D_2'$ of $D_1$
and
$D_2$, and the exceptional divisor $E$. Correspondingly, the boundary of $\orb(\tilde W_z,
\tilde D \cap \tilde W_z)$ consists of
$\tilde p^{-1}(D_1'),\ \tilde p^{-1}(D_2')$ and $\tilde p^{-1}(E)$. Let
$\tilde z_1=E\cap D'_1$ and $\tilde z_2=E\cap D'_2$, as shown in Figure \threelines.

Define $h$ on $\tilde p^{-1}(D_1') \cup \tilde p^{-1}(D_2')$ by the formulas
$$
\halign{\hskip 2cm \hfil$#$& = $\Bigl($\hfil$#$\hfil, &\hfil$#$\hfil,
&\hfil$#$\hfil, &\hfil$#$\hfil$\Bigr)$&\quad #&\quad$#$\hfil&\hskip 1
cm\hfill(#)\cr  h(w_1,0,\theta_1, \theta_2) & \eta_1(|w_1|) \frac {w_1}{|w_1|}&
0 &
\theta_1& \theta_2& in&  \tilde p^{-1}(D_1')&\honDone\cr  
h(0, w_2,\theta_1, \theta_2) & 0& \eta_1(|w_2|) \frac {w_2}{|w_2|} &
\theta_1& \theta_2& in&  \tilde p^{-1}(D_2'),&\honDtwo\cr}
$$
where $\eta_1$ is the same function which we used in the proof of Theorem \blowupofordpoint\
(see Figure \FigVIIetasfordoublepts).  This leaves $\tilde p^{-1}(E)$. Recall that 
$\tilde W_z$ is the locus in
$W
\times
\Proj^1$ given by the equation $w_1U_2=w_2U_1$, where $\bvec {U_1}{U_2}$ are homogeneous
coordinates. We will break it into two pieces:

\noindent
$\bullet$ $E_1$ where $|U_2|\le |U_1|$, which we will parametrize by $(u_2,\ w_1)$, setting
$u_2=U_2/U_1$, and

\noindent
$\bullet$ $E_2$ where $|U_1|\le |U_2|$, which we will parametrize by $(u_1,\ w_2)$, setting
$u_1=U_1/U_2$.
 
The space $\tilde p^{-1}(E_1)$ thus has coordinates $(u_2,\phi_2, \theta_1)$ and 
the space $\tilde p^{-1}(E_2)$ thus has coordinates $(u_1,\phi_1, \theta_2)$.  
Using these coordinates, we will define $h$ on $\tilde p^{-1}(E)$ as follows :
$$\halign{\hskip .5cm \hfil$#$& = $\Bigl($\hfil$#$\hfil, &\hfil$#$\hfil, &\hfil$#$\hfil,
&\hfil$#$\hfil$\Bigr)$&\quad #&\quad$#$\hfil&\hskip 1 cm\hfill(#)\cr 
h(u_2, \phi_2, \theta_1) & \eta_2(|u_2|) e^{i\theta_1}&  0 & \theta_1&
\theta_1+\phi_2 &in  &\tilde p^{-1}(E_1)&\honEone\cr 
h(u_1,\phi_1,\theta_2) & 0& \eta_2(|u_1|) e^{i\theta_2} & \phi_1+\theta_2& \theta_2
&in  &\tilde p^{-1}(E_2),&\honEtwo\cr}
$$
\goodbreak
where $\eta_2 (t) = r-tr$  (see Figure \FigVIIetasfordoublepts).
\bigskip

\centerline{\BoxedEPSF{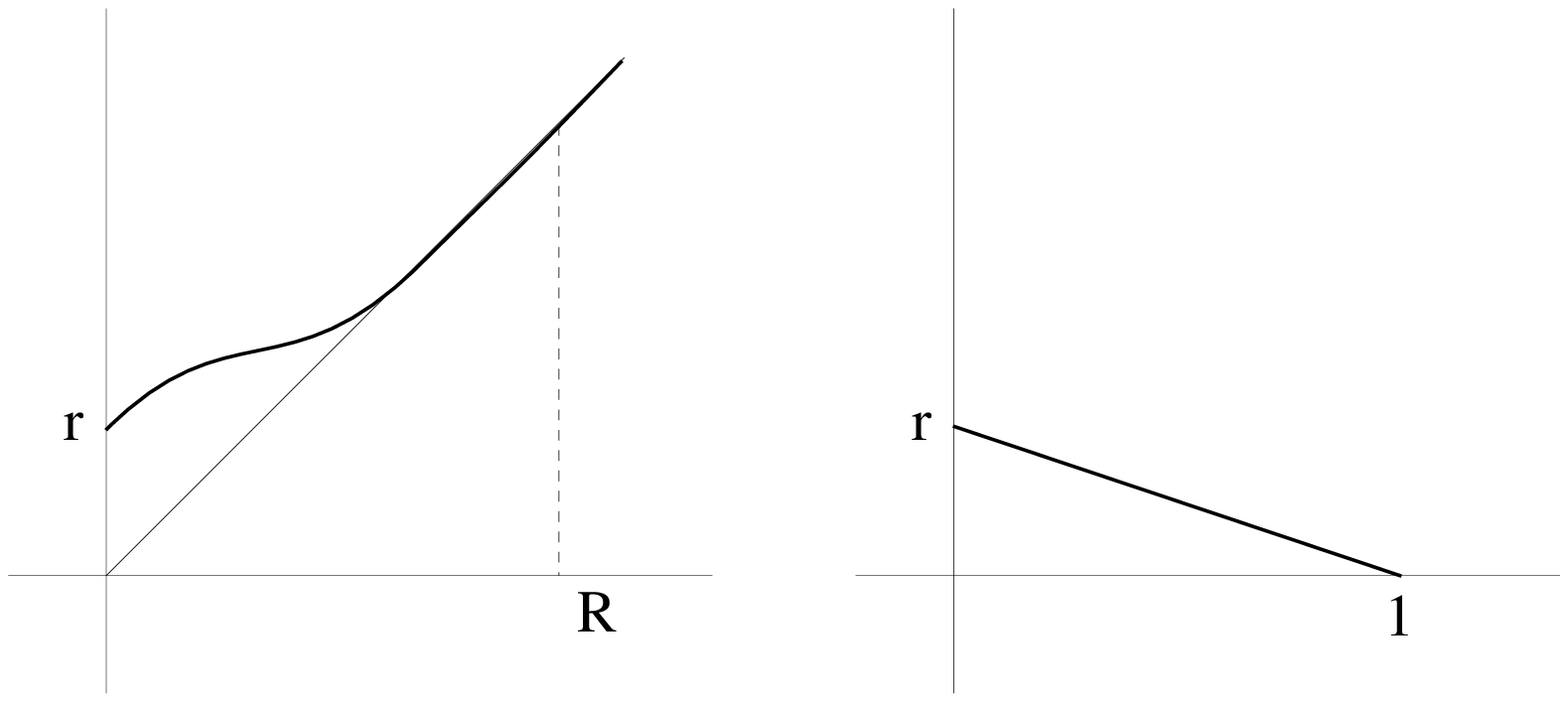 scaled 500}}
\nobreak
\figcap \FigVIIetasfordoublepts {The graphs of $\eta_1$ (left) and of $\eta_2$
(right).}

There are three places where we must check the compatibility of these formulas: above $\tilde
z_1$, above $\tilde z_2$ and above the circle $|u_1|=|u_2|=1$.

Above $\tilde z_1$, the functions $u_2$ and $w_1$ are local coordinates, and since $w_2=u_2w_1$,
we have $\theta_2=\phi_2+\theta_1$.  This means that under the mapping (\honDone), the point 
$(0,0,\theta_1,\theta_2)$ in $\tilde p^{-1}(D_1')$ maps to $(re^{i\theta_1},0,\theta_1,
\theta_2)$, whereas the same point viewed in $\tilde p^{-1}(E_1)$ has coordinates $(0, \phi_2,
\theta_1)$ and maps under (\honEone) to $(re^{i \theta_1}, 0, \theta_1, \theta_1+\phi_2)$, which is
the same point.

The analogous check above $\tilde z_2$ is left to the reader.

Above the circle $u_1=u_2=1$, we must pass from the coordinates $(u_2,w_1)$ to the
coordinates $(u_1, w_2)$.  Under (\honEone), the point $(e^{i\phi_2},\phi_2,\theta_1)$ in $\tilde
p^{-1}(E_1)$ maps to $(0,0,\theta_1, \theta_1+\phi_2)$, whereas under (\honEtwo), the same
point, viewed as an element of $\tilde p^{-1}(E_2)$, has coordinates $(e^{i\phi_1}, \phi_1,
\theta_2) = (e^{-i\phi_2}, -\phi_2,\theta_1+\phi_2)$, and also maps to $(0,0,\theta_1,
\theta_1+\phi_2)$.

Finally, observe that the map $h$ as constructed does transform the circle action on
$\orb(\tilde D)$ into the circle action on $Z$.  In $\tilde p^{-1}(E_1)$, the circle action is 
$$
\Theta*(u_2, \phi_2, \theta_2) = (u_2, \phi_2, \theta_2+\Theta),
$$
which under (\honEone) becomes 
$$
\Theta*(w_1,0,\theta_1, \theta_2) = (w_1,0,\theta_1+\Theta, \theta_2+\Theta).
$$
The computation for $\tilde p^{-1}(E_2)$ is identical.

This ends the construction of $h$ on $\partial Z$; now to extend it to the interior.

We will call $B$  the subset of $W$ in which
$$
\left(|w_1|^2 + |w_2|^2\right)^2 < 4r^2|w_1w_2|.
$$
Again, we can think of $B$ as a subset of $W$, or of $\tilde W_z$, or of $\orb(W, D \cap W)$, or
of $\orb(\tilde W_z, \tilde D \cap \tilde W_z)$, since all the canonical projections are isomorphisms
on the inverse images of $B$. 

\centerline{\BoxedEPSF{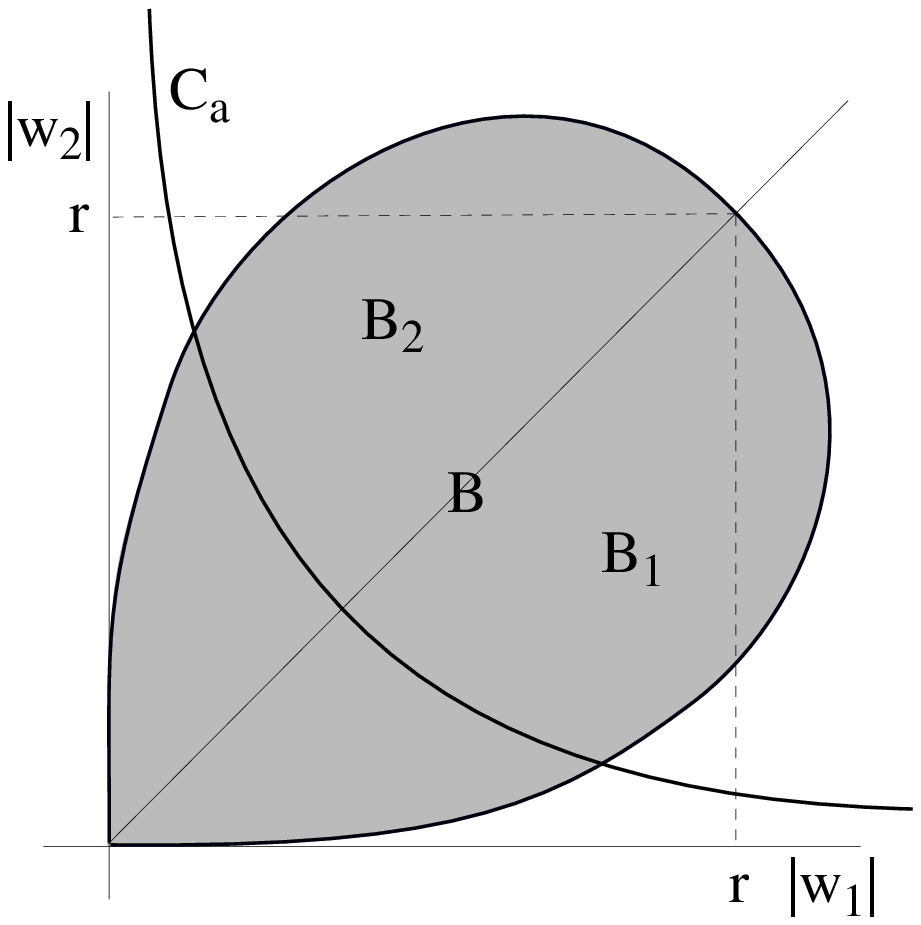 scaled 600}}
\figcap \picofB {The region in $|w_1|,|w_2|$ defining $B$.}

This region contains $E$ in its closure, since every line through the origin has a
neighborhood of the origin in $B$, except for the axes themselves. This is illustrated in Figure 
\picofB. 

Let  $W_1, W_2 \subset W$ be the regions where $|w_1| \ge |w_2|$ and
$|w_1| \le |w_2|$ respecively. Then the space $B$ is the union of the two subsets $B_1= W_1 \cap B$
and
$B_2= W_2 \cap B$; these have $E_1$ and $E_2$ in their closures in $\tilde W_z$ (resp.
$\tilde p^{-1}(E_1)$ and $\tilde p^{-1}(E_1)$ in their closures in $\orb(\tilde W_z, \tilde D\cap
\tilde W_z)$). 

The curve $C_a$ of equation $w_1w_2=a^2$ intersects $B$ in an annulus $B_a$ for $0<|a|<r$,
which shrinks to a circle when $|a|=r$ and is empty for $|a|>r$.  We will extend $h$ so
that it maps each $C_a$ into itself.  Moreover, we will use the local coordinates $u_2, w_1$ in
$B_1$ and $u_1, w_2$ in $B_2$, and keep the coordinates $w_1,w_2$ outside of $B$. Throughout,
we will use the coordinates $w_1,w_2$ in the range.  

In these coordinates, our mapping $h$ is
written  $$
\halign{\hskip 1.5cm \hfil$#$& = $\Bigl($\hfil$#$\hfil&,\ \hfil$#$\hfil$\Bigr)$&\quad #&\quad$#$\hfil&\hskip 1
cm\hfill(#)\cr 
h(w_1,w_2) & \eta_{1,|a|}(|w_1|) \frac {w_1}{|w_1|}&\frac{|a|^2}{\eta_{1,|a|}(|w_1|)} \frac
{w_2}{|w_2|}& in&  W_1-B_1&\hextonWone\cr
\noalign{\medskip}
h(u_2, w_1) & \eta_{2,|a|}(|u_2|) e^{i\theta_1}&\frac{a^2}{\eta_{2,|a|}(|u_2|)}e^{-i \theta_1}
&in  &B_1&\hextonBone\cr
\noalign{\medskip} 
h(u_1,w_2) & \frac{a^2}{\eta_{2,|a|}(|u_1|)}e^{-i \theta_2}&\eta_{2,|a|}(|u_1|) e^{i\theta_2}&
in &B_2&\hextonBtwo\cr 
\noalign{\medskip} 
h(w_1, w_2) & \frac{|a|^2}{\eta_{1,|a|}(|w_2|)}\frac {w_1}{|w_1|}& \eta_{1,|a|}(|w_2|)\frac
{w_2}{|w_2|}& in&  W_2-B_2,&\hextonWtwo\cr}
$$
where $\eta_{2,|a|}(t) = r-(r-|a|)t$ and 
$$
\eta_{1,\rho}:\left[\sqrt{\rho(r+\sqrt{r^2-\rho^2})}, \,\infty\right) \to
\left[\eta_{2,\rho}\Bigl(\frac{(r-\sqrt{r^2-\rho^2}}{\rho}\Bigr),
\,\infty\right) 
$$ 
is a diffeomorphism, to be chosen below.

There are now three compatibilities to check: on the boundary $\partial Z$, on the boundary of $B$
and when
$|u_1|=|u_2|=1$.  The first will be satisfied if we require that  $$
\lim_{|a| \to 0} \eta_{1,|a|} = \eta_1\quad\text{and}\quad \lim_{|a| \to 0} \eta_{2,|a|} = \eta_2.
$$
For the third, our extension of $h$ is precisely the identity on this locus. Indeed, using
formula (\hextonBone), we see
$$
\eta_{2,|a|}(|u_2|) e^{i\theta_1}= \eta_{2,|a|}(1) e^{i\theta_1}= |a|e^{i\theta_1}=w_1
$$
and  
$$
\frac{a^2}{\eta_{2,|a|}(|u_2|)}e^{-i \theta_1}= \frac {a^2}{w_1}= w_2
$$
there; the check for (\hextonBtwo) is identical.
  
For the second, we need to guarantee that
$$
\eta_{1,|a|}\left(\sqrt{|a|\bigl(r+\sqrt{r^2-|a|^2}}\bigr)\right) =
\eta_{2,|a|}\left(\frac{r-\sqrt{r^2-|a|^2}}{|a|}\right).
$$
A look at the curve on which the graph of $\eta_{1,\rho}$ must start (Figure
\graphaofetasfordblepts) shows  that this is easy to accomplish.
\bigskip 

\centerline{\BoxedEPSF{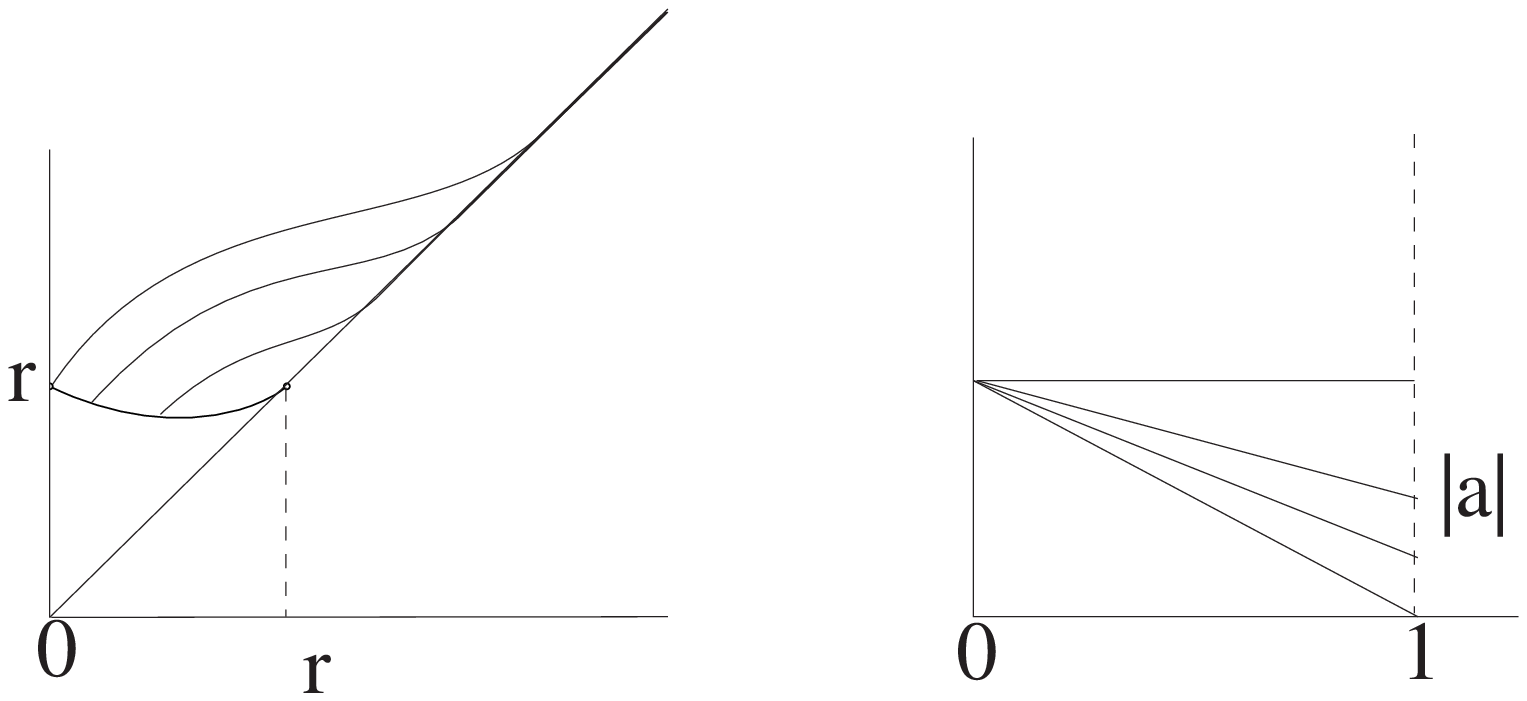 scaled 500}}
\longfigcap  \graphaofetasfordblepts {The graphs of the functions $\eta_{1,\rho}$ (left) and
$\eta_{2,\rho}$ (right) for various values of $\rho$.}
\qed

Again, in order to prove Theorem \topofprojlim\ we will need to control how much $h$ distorts
distances.

\theorem \controlblowupdblept For any Riemann metric $\mu$ on $\orb(D)$ which is
continuous in the coordinates $\rho, \theta, \Theta$, and for any $\epsilon>0$ and
$\alpha<1$, we can choose $h$ so that 
$$
d_\mu(h(x), \tilde \pi(x))<\epsilon
$$
and
$$
d_\mu\bigl( h(x),  h(y)\bigr) \ge \alpha d_\mu\bigl(\tilde \pi(x), \tilde \pi(y)\bigr)
$$
for any $x,y$ in $\orb(\tilde D)$.
\endstat

\proof This is a matter of choosing $r$ and $\eta_1$ correctly in the
constructions above.

In our coordinates, the mapping $\tilde \pi$ is written
\medskip

\halign{\hskip 2.5 cm \hfil$#$&= $#$\hfil\quad if &$#$\hfil\cr
\tilde \pi\circ h^{-1}(w_1,0,\theta_1, \theta_2)& (0,0,\theta_1, \theta_2)  &|w_1|\le
r\cr
\tilde \pi\circ h^{-1}(w_1,0,\theta_1, \theta_2)& (\eta_1(|w_1|)\frac{w_1}{|w_1|},0,\theta_1,
\theta_2) & |w_1|\ge r\cr
\tilde \pi\circ h^{-1}(0,w_2,\theta_1, \theta_2)& (0,0,\theta_1, \theta_2) & |w_2|\le r\cr
\tilde \pi\circ h^{-1}(0,w_2,\theta_1, \theta_2)& (\eta_1(|w_2|)\frac{w_2}{|w_2|},0,\theta_1,
\theta_2) & |w_1|\ge r.\cr}
\noindent This map  collapses the thickened torus $|w_1|,|w_2|\le r$ back onto the torus above
$z$. Clearly
$$
d_\mu( h(x), \tilde \pi (x))<\epsilon 
$$
will be small if $r$ is small and $\eta_1$ is uniformly close to the identity.  Again the
second part is an application of Lemma \natylemmaontubneighborhoods.
\qed

\subheading{Infinitely many blow-ups} Suppose that  we repeat infinitely many
times the following procedure, as in Section IV.
 
Take a surface $X_0$ containing a divisor with normal crossings $D_0\subset X_0$; although it is
not essential, we will assume that $X_0$ is compact. Choose a point $z_0 \in D_0$, blow it up to
create a surface
$X_1=
\tilde{(X_0)}_{z_0}$, with a projection $\pi_1:X_1 \to X_0$; set $D_1 = \pi_1^{-1}(D_0)$.  Now
choose a new point $z_1 \in D_1$, etc.  Denote by 
$$
\tilde \pi_i:\orb(X_i, D_i) \to \orb(X_{i-1}, D_{i-1})
$$
the map induced from $\pi_i$.

Theorems  \blowupofordpoint\ and \blowupdoublepoint\ assert that at each stage the pair
$\orb(D_i) \subset \orb(X_i,D_i)$ is homeomorphic to the pair $\orb(D_0) \subset
\orb(X_0,D_0)$.  It seems reasonable to think that the same will remain true in the limit,
and it is.

\theorem \topofprojlim The projective limit of pairs
$$
\varprojlim (\orb (D_i),\tilde \pi_i) \subset \varprojlim ((\orb(X_i,D_i),\tilde \pi_i)
$$
is homeomorphic to the pair $\orb (D_0) \subset \orb(X_0,D_0)$, and the canonical map 
$$
\varprojlim (\orb(X_i, D_i),\tilde \pi_i) \to 
\orb(X_0, D_0)
$$ 
can be approximated by homeomorphisms.
\endstat

\proof Theorems  \blowupofordpoint\ and \blowupdoublepoint\ guarantee that we can choose
homeomorphisms 
$$
h_i:\orb(X_i,D_i) \to \orb(X_{i-1}, D_{i-1}). 
$$

By induction, we will give $\orb(X_i,D_i)$ the metric which makes $h_i$ an isometry; and
Theorems \controlledblowupordpt\ and \controlblowupdblept\ then show that each $h_i$ can further
be chosen so that for any sequences $\epsilon_i>0,\ 0<\alpha_i<1$,  
\roster
\item the composition $\tilde \pi_i \circ h_i^{-1}$ is $\epsilon_i$ close to the identity map (of
$\orb (X_{i-1},D_{i-1})$);
\smallskip
\item if $x_i,y_i \in \orb(D_i)$ and $x_{i-1} = \tilde \pi_i(x_i),\ y_{i-1} =
\tilde \pi_i(y_i)$, then
$$
d(h_i(x_i),h_i(y_i)) \ge \alpha_i d(x_{i-1},y_{i-1}).
$$
\endroster

Moreover, again by choosing the $r$ sufficiently small, we can assume that for any compact
subset $K \subset X-D$, there exists $m$ such that the mapping $h_i$ coincides with $\tilde
\pi$ on $K$ for $i>m$.  

\remark The open set $X-D$ is naturally embedded as a subset in all $\orb(X_i,D_i)$; in terms
of this natural embedding, the choice of $r$ above guarantees that $h_i$ is the identity on the
compact subset
$K$ for $i>m$.  

Denote by $g_n$ the composition
$$
\varprojlim_i \orb (X_i,D_i) \to \orb(X_n,D_n) \overset {h_1 \circ \dots \circ h_n} \to {\to}
\orb (X_0,D_0),
$$
 where the first map is the canonical projection.  The first condition above guarantees
that the $g_n$ converge uniformly if the $\epsilon_n$ form a convergent series.  Call $g
= \lim_{n \to \infty} g_n$.

The mapping $g$ is obviously surjective. The second condition guarantees that $g$ is injective on
$\varprojlim \orb(D_i)$ if $\prod
\alpha_n >0$. Moreover, it is injective on $X-D$ since the sequence is eventually constant and
injective on every compact subset.  Therefore it is injective.

Finally, since $g$ is bijective, continuous and the domain is compact (a projective limit of
compact sets is compact),  it is a homeomorphism.  
\qed

%% file: Veselov.HenonExample.new
\heading VIII. Real oriented blow-ups and the Hopf fibration \endheading

The Hopf fibration $p:S^3 \to S^2$ is a famous example from topology, where the
circle acts on the 3-sphere, and the quotient is the $2$-sphere. If we think of $S^3$
as the unit sphere in $\C^2$, and $S^2$ as $\Proj^1_\C$, then it can be written as
$$
p:(z_1,z_2) \mapsto [z_1:z_2]. 
$$ 
It can also be thought of as the quotient of the $3$-sphere by the circle action
$$
\R/2\pi \Z \times S^3 \to S^3,\quad  \Theta*(z_1,z_2) =(e^{i\Theta}z_1,e^{i\Theta}z_2)
$$
The Hopf fibration seems
a natural candidate to be the real oriented blow-up of a projective line in a surface,
and it is.  But it comes up in two different settings, which differ in a subtle but
essential way.

More specifically, we will consider a line in $\Proj^2_\C$, for instance the line at
infinity $l_\infty$, and the exceptional divisor $E$ that you get by blowing up the
origin.  The plane $\Proj^2_\C$ has a  real oriented blow-up along both of these
lines; we can  imagine these real oriented blow-ups as the boundaries of tubular
neighborhoods of the two lines respectively.  If we take the tubular neighborhood of
$E$ to be the open unit ball, and the tubular neighborhood of $l_\infty$ to be the
complement of the closed ball, then both have the same boundary $S^3$.  We will see
in Propositions \hopfisom\ and \modifiedHopf\ that there really is a homeomorphism of
the $3$-sphere with
$\orb(E)$ and $\orb(l_\infty)$, and that the first homeomorphism transforms the canonical
circle action on $\orb(E)$ into the Hopf action above and that the other transforms the
standard circle action into
$$
\Theta*(z_1,z_2) =(e^{-i\Theta}z_1,e^{-i\Theta}z_2).
$$

If you think of rotating a line in $\C^2$ in the direction of the boundary of a
disc, or the boundary of the outside of a disc, you will see where the minus sign
comes from.

You might think that this minus sign is the difference between the two cases, but it
isn't: the antipodal map $S^3 \to S^3$ is an orientation-preserving
diffeomorphism which transforms one circle action into the other, so these two
circle actions are essentially the same.  The actual difference between the two
cases is the orientation of $S^3$, as the boundary of the unit ball in the case of
$\orb(E)$, and as the boundary of the outside of the unit ball in the case of
$\orb(l_\infty)$.  This may seem a minor problem, but it turns out to have essential
consequences.

Indeed, if you further blow up a point $z \in l_\infty$, making 
$$
\pi:\widetilde {(\Proj^2_\C)}_z \to \Proj^2_\C,
$$
set $D=\pi^{-1}(l_\infty)$, and still denote by $p:\orb (D) \to D$ the projection,
then there are fibers of $p$ which link with linking number 0, and there are none
which link with linking number $\pm 2$.

On the other hand, if you blow up a point $z \in E$, call $E_1$ the new exceptional
divisor and $E'$ the proper transform of $E$, set $D=E' \cup E_1$, and again
denote by
$p:\orb (D) \to D$ the projection, then there are fibers which link with linking number
2, and there are none which link with linking number 0.

These two situations differ geometrically, not just by a conventional sign.  Let us
now justify these claims.

\subheading {The real blow-up of an exceptional divisor} 
Let $\pi:\widetilde{\C^2}_0\to \C^2$ be the plane $\C^2$ blown up at the origin,
and $E=\pi^{-1}(0,0)$ be the exceptional divisor. Recall the two charts
$$
u_1=x,\ v_1=\frac yx \quad \text{and}\quad u_2=y,\  v_2=\frac xy
$$
which provide coordinates on overlapping open sets which cover $\widetilde{\C^2}_0$.
In each of these charts, $\orb(\widetilde{\C^2}_0,E)$ is the closure of the set 
$$
\set{(u_j,v_j,\theta_j)\in \C^*\times \C \times
\R/2\pi\Z}{\frac{u_j}{|u_j|} = e^{i\theta_j}},\ j=1,2
$$
in $\C\times \C \times \R/2\pi\Z$, and that $p^{-1}(E)$ is given by the equation
$u_j=0$.

Thus  $\orb(E)$ is parametrized by the two charts $(v_1,\theta_1)$ and
$(v_2,\theta_2)$, and on the overlap where $v_1, v_2 \ne 0$, we have
$$
v_2 = \frac 1{v_1},\ \theta_2 = \theta_1 + \arg v_1.
$$

\proposition {\hopfisom} The mapping $S^3 \to \orb(E)$ given by 
$$
(z_1,z_2) \mapsto
\cases
v_1 = \frac{z_2}{z_1},\ \theta_1 = \arg z_1&\quad\text{if $z_1 \ne 0$}\\
v_2 = \frac{z_1}{z_2},\ \theta_2 = \arg z_2&\quad\text{if $z_2 \ne 0$}
\endcases
$$
is a diffeomorphism which takes  the orientation of $S^3$ as the boundary of
the unit ball in $\R^4$ to the chosen orientation of $\orb(E)$, and which transforms the
Hopf circle action 
$$
\Theta *(z_1,z_2)=(e^{i\Theta}z_1,e^{i\Theta}z_2)
$$
into the canonical circle action
$$    
\Theta*(v_j,\theta_j) = (v_j, \theta_j+\Theta)
$$
on $\orb(E)$.
\endstat

\proof The main thing to check is that the mapping is compatible with the
identification $\theta_2 = \theta_1+ \arg v_1$, which becomes 
$\arg z_2 = \arg z_1 + \arg(z_2/z_1)$. 

The map is injective: if we know $v_1, \theta_1$, then from $z_2= v_1z_1$ and the
equation
$|z_1|^2+|z_2|^2=1$ we see
$$
|z_1|^2 = \frac 1{1+|v_1|^2},
$$
and since we also know the argument of $z_1$, we know $z_1$, hence $z_2$.

The surjectivity is also clear from the argument above.

The compatibility with the circle action is 
$$
\frac {e^{i\Theta}z_1}{e^{i\Theta}z_2} = \frac{z_1}{z_2}\quad, 
\quad \arg e^{i\Theta} z_1 = \arg z_1 + \Theta.
$$
The vectors
$$
\bmatrix 0\\1\\0\\0\endbmatrix,\bmatrix 0\\0\\1\\0\endbmatrix,\bmatrix
0\\0\\0\\1\endbmatrix 
$$
form a direct basis of $T_{(1,0)}S^3$ if you orient $S^3$ as the boundary of the ball.
The first of these vectors is tangent to the oriented orbit through $(1,0)$,
whereas the last two project under $p$ to a direct basis of $T_{p(1,0)}E$.
\qed

\subheading {The real oriented blow-up of a line in $\Proj^2_\C$}

We now examine the blow-up of a line in $\Proj_\C^2$, which we will take to be the line
at infinity  $l_\infty$.  This time, local coordinates on a neighborhood of $l_\infty
\in \Proj^2_\C$ are
$$
u_1= \frac 1x,\ v_1 = \frac yx\quad \text{and}\quad u_2 = \frac 1y, \ v_2 = \frac xy.
$$
This leads to the charts $\C \times \R/2\pi \Z \to \orb (l_\infty)$ given by
$$
(v_j, \theta_j),\ j=1,2.
$$
On the overlap $v_1, v_2 \ne 0$ these coordinates are identified by
$$
v_2 = \frac 1{v_1}, \theta_2 = \theta_1- \arg v_1.
$$

This is also a variant of the Hopf fibration.

\proposition {\modifiedHopf} a) The mapping $S^3 \to \orb (l_\infty)$ given by 
$$
(z_1,z_2) \mapsto
\cases
v_1 = \frac{z_2}{z_1},\ \theta_1 = -\arg z_1&\quad\text{if $z_1 \ne 0$}\\
v_2 = \frac{z_1}{z_2},\ \theta_2 = -\arg z_2&\quad\text{if $z_2 \ne 0$}
\endcases
$$
is a diffeomorphism which carries the orientation of $S^3$ as the boundary of the
complement of the 4-ball to the standard orientation $\orb(l_\infty)$. 

b) This diffeomorphism transforms the circle
action 
$$
\Theta *(z_1,z_2)=(e^{-i\Theta}z_1,e^{-i\Theta}z_2)
$$
into the canonical circle action
$$    
\Theta*(v_j,\theta_j) = (v_j, \theta_j+\Theta).
$$
\endstat

The proof is identical to the case above.

\subheading {Stereographic projection}  What really makes these two cases different? 
It is the orientation which $S^3$ acquires.  We invite the reader to check that with
the definition of
\orientrealblowup, $\orb(D)$ is oriented so that a tangent vector to an oriented orbit,
followed by an oriented basis tangent to a section, defines the orientation of the
space
$\orb(D)$.  Thus if we reverse the direction of the circle orbits (and keep the
orientation of the base), we change the orientation of the space.

To better picture the difference of the two cases above, we will identify
$S^3-\{(0,0,0,1)\}$ to $\R^3$ via stereographic projection.

\proposition \computelinkingnumber Set 
$z_1=x_1+iy_1= r_1 e^{i\theta_1},\  z_2= x_2+iy_2= r_2e^{i\theta_2}$.

a) The
stereographic projection $S^3-\{\text{north pole}\} \to \R^3$ given by
$$\eqalign{
X&=\frac {x_1}{1-y_2}\cr
Y&=\frac {y_1}{1-y_2}\cr 
Z&=\frac {x_2}{1-y_2}}\quad \text{resp.}\quad 
\eqalign{
X&=\frac {x_1}{1-y_2}\cr
Y&=\frac {y_1}{1-y_2}\cr 
Z&=-\frac {x_2}{1-y_2}}\peqno \stereographic
$$
maps the orientation of the 3-sphere as the boundary of the inside (resp. outside) of
the unit ball in $\C^2$ to the standard orientation of $\R^3$.

b) The tori $T_r$, $0 \le r \le 1$, with parametric equations 
$$ \eqalign{ 
X&= \left(\frac{\sqrt{1-r^2}}{1-r\sin\alpha}\right)\cos \beta \cr 
Y&= \left(\frac{\sqrt{1-r^2}}{1-r\sin\alpha}\right)\sin \beta \cr
Z&=  \frac{r \cos \alpha}{1-r \sin\alpha}}\quad \text{resp.}\quad 
\eqalign{ 
X&= \left(\frac{\sqrt{1-r^2}}{1-r\sin\alpha}\right)\cos \beta \cr 
Y&= \left(\frac{\sqrt{1-r^2}}{1-r\sin\alpha}\right)\sin \beta \cr
Z&=  -\frac{r \cos \alpha}{1-r \sin\alpha}}   \peqno \stereoinpolar
$$
are the images of the tori $|z_1|=\sqrt{1-r^2}, |z_2|=r$.
They are tori of revolution around the $z$-axis, with two
singular tori: the unit circle in the $(x,y)$-plane (for $r=0$) and the
$z$-axis (for $r=1$). The other tori are obtained by rotating around the $z$-axis the
circle in the $(X,Z)$-plane with center at $X= 1/\sqrt{1-r^2}, Z=0$ and radius
$r/\sqrt{1-r^2}$.

c) The circle actions
$$
\Theta*(z_1,z_2) = (e^{i\Theta}z_1, e^{i\Theta}z_2)\quad \text{resp.}\quad 
\Theta*(z_1,z_2) = (e^{-i\Theta}z_1, e^{-i\Theta}z_2)
$$
become, in these coordinates, 
$$ 
\Theta*(r,\alpha,\beta) = (r, \alpha+\Theta, \beta+\Theta)\quad \text{resp.}\quad
\Theta*(r,\alpha,\beta) = (r, \alpha-\Theta, \beta-\Theta). \peqno \circleaction 
$$
\endstat

\proof We will prove only the parts concerning the blow-up of the line at infinity,
which appear in the right-hand column, as those are the ones we are concerned with.

a) Recall that if $M$ is a smooth k-manifold with boundary and
$m\in \partial M$, then 
$$
v_1, v_2, \dots, v_{k-1} \in T_m \partial M
$$
 define the
standard orientation of $\partial M$, if $v_0, v_1, \dots, v_{k-1}$ define the standard
orientation of $M$, where $v_0$ is a  outward-pointing vector of $T_mM$.  
The vectors
$$
v_0=\bmatrix 0\\0\\0\\1\endbmatrix,\ v_1= \bmatrix 1\\0\\0\\0\endbmatrix,\ 
v_2=\bmatrix 0\\1\\0\\0\endbmatrix,\ v_4=\bmatrix \phantom{-} 0\\
\phantom{-}0\\-1\\ \phantom{-}0\endbmatrix
$$
are a direct basis of $\R^4$ with $v_0$ pointing out of the outside of the unit
sphere at the south pole, so $v_1, v_2, v_3$ define the orientation of $S^3$ as the
boundary of the outside. Under the stereographic projection (\stereographic) these
vectors get mapped to  the vectors 
$$
\bmatrix 1\\0\\0\endbmatrix, \bmatrix 0\\1\\0\endbmatrix, \bmatrix 0\\0\\1\endbmatrix, 
$$ which give the standard orientation of $\R^3$.

Rewriting the mapping (\stereographic) in polar coordinates $$x_1=\sqrt
{1-r^2}\cos\beta,\,\   y_1=\sqrt {1-r^2}\sin\beta,\,\  x_2=r\cos\alpha,
\,\  y_2=r\sin\alpha, \,\  0\le r \le 1,$$ we will immediately  get 
formulas (\stereoinpolar) and statement c) of the proposition.
 
To prove part b) it is enough to check that the curves in the $(X, Z)$ plane given by
the parametric  equations
$$
X= \frac{\sqrt{1-r^2}}{1-r\sin\alpha} \quad\text{and}\quad
Z=  \frac{r \cos \alpha}{1-r \sin\alpha} 
$$
are actually the circles with centers at $(1/\sqrt{1-r^2}, 0)$ and radii 
$r/\sqrt{1-r^2}$. Indeed,
$$
\left(\frac {\sqrt {1-r^2}} {1-r\sin\alpha} - \frac 1{\sqrt {1-r^2}} \right )^2
+\left(\frac{r\cos \alpha}{1-r\sin \alpha}\right)^2=  \frac {r^2(r-\sin\alpha)^2 +
r^2\cos^2\alpha(1-r^2)} {(1-r\sin\alpha)^2(1-r^2)}= \frac {r^2}{1-r^2}
$$
\qed

\subheading {Blowing up one point of a line}  Suppose we now want to blow up one point
of $E$ or of $l_\infty$.  We can symbolically represent the two variants of the Hopf
fibration in the following picture.

\centerline{\BoxedEPSF{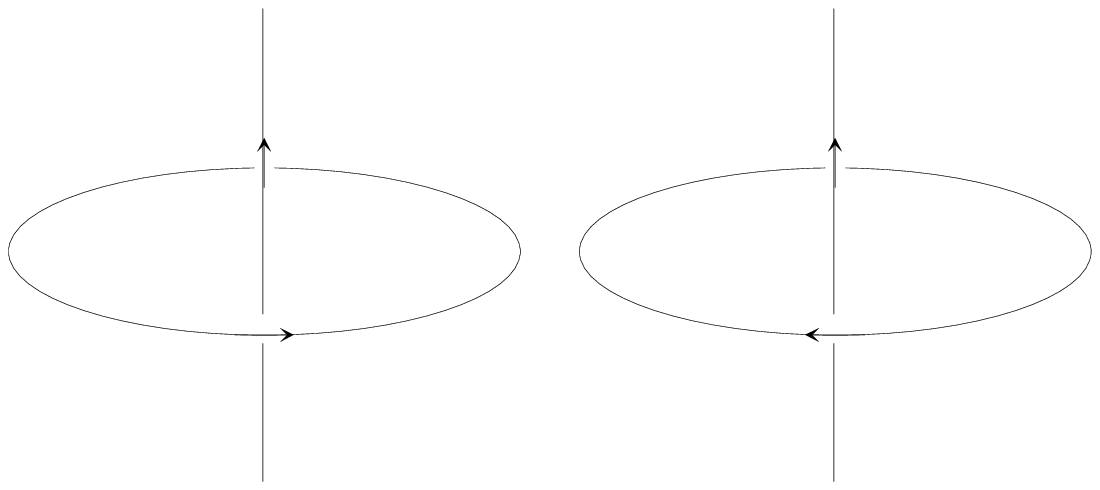 scaled 600}}
\medskip
Figure \twohopf. The real oriented blow-up of an exceptional divisor (left) and of the
line at infinity (right).

In both $\orb(\tilde\Proj^2_0, E)$ and  $\orb(\tilde\Proj^2_2, l_\infty)$,
$S^3-\{\text{north pole}\}$ is identified by stereographic projection with
$\R^3$ with its standard orientation.  In both, the $z$-axis oriented up is one
oriented fiber, and the unit circle in the $x,y$-plane is another, oriented
counterclockwise in the first, clockwise
in the second. Now we will blow up
one point $\bold p$, on $E$ or $l_\infty$ respectively, and see what happens to the real
oriented blow-ups.  Choose the unit circle in the $x,y$-plane to correspond to
$\bold p$, and thicken it as in Theorem \blowupofordpoint, corresponding to the
region $r\le \epsilon$ for some small $\epsilon>0$.  This solid torus $T$ intersects
the  $x,z$-plane in a disc which has a canonical orientation, since it projects
homeomorphically to $E$ or $l_\infty$, and the boundary orientation of the boundary of
the disc in both cases corresponds to $\alpha$ increasing from $0$ to $2\pi$. This
gives one generator $a$ of $H_1(\partial T)$.  Letting $\beta$ increase from $0$ to
$2\pi$ defines a second generator $b$ of $H_1(\partial T)$; this curve also
corresponds to the boundary of a disc on the outside of $T$, which also has a boundary
orientation, and our conventions are consistent.

In these coordinates, an oriented orbit of the circle action on the torus has homology
class
$a +b$ in the case of $\orb(E)$, and homology class $-a-b$ in the
case of $\orb(l_\infty)$.

If we change the circle action inside the torus, as in Theorem \blowupofordpoint,
the homology class of an orbit is $a+2b$ in the case of $\orb(E)$, and
$-b$ in the case of $\orb(l_\infty)$.  This is not a difference in sign
convention: fibers inside the torus  link with linking number $- 2$ in the first case,
and in the second case have linking number $0$  (we are following Milnor's convention
\cite{Mil4} for the sign of the linking number).

\centerline{\BoxedEPSF{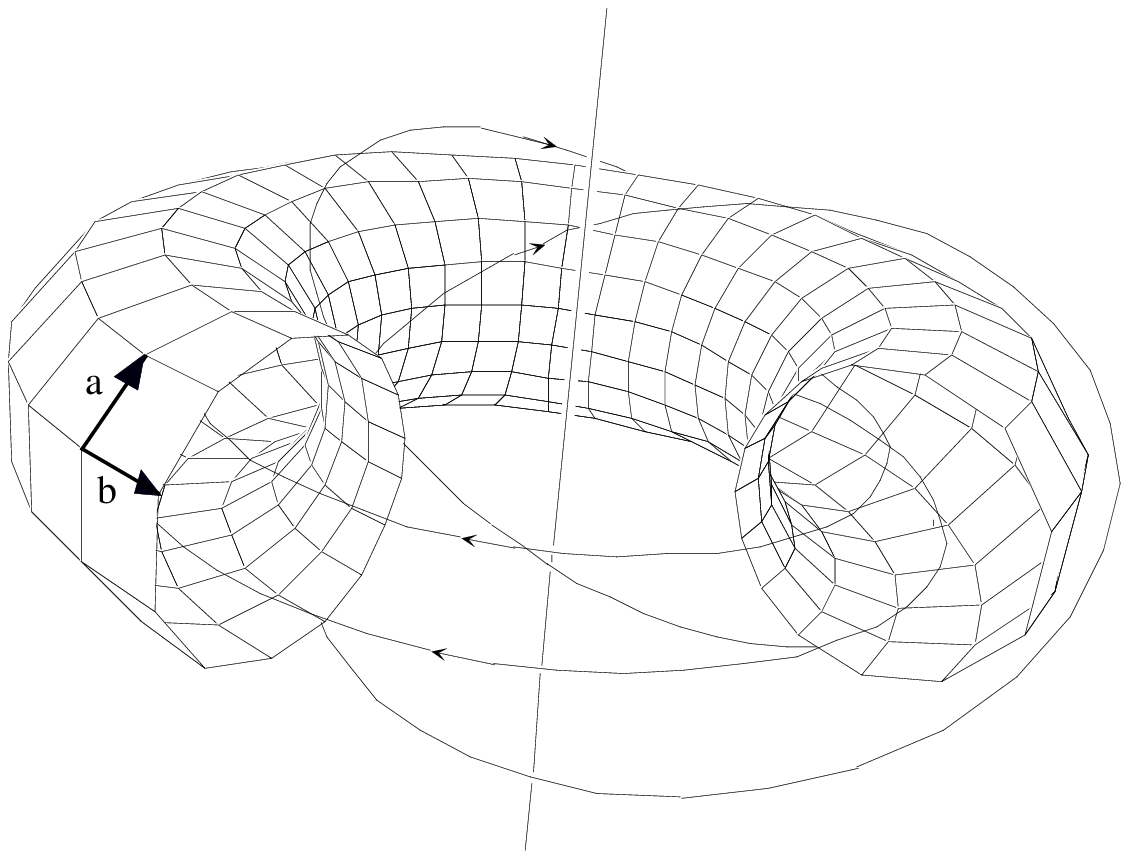 scaled 800}}

\longfigcap \rightsecondblowup {Blow up $\Proj^2$ at a point $z\in l_\infty$, and
then take the real oriented blow-up of the divisor consisting $l_\infty$ and the
exceptional divisor. You obtain a 3-sphere containing a torus, and a circle action on
both components of the complement of the torus.  This picture represents the
stereographic projection of this space, with the inside of the torus corresponding to the
exceptional divisor, where the circle orbits do not link, and the outside
corresponding to the line at infinity, where the circles link with linking number 1.
Curves describing the homology classes $a$ and $b$ are drawn on the torus; note that
the circle orbits on the outside are in the class
$-a-b$ and those inside are in the class $-b$.}

So these two alternatives both exist, and even a trusting reader might wonder whether
we have gotten it right: there were about six places where a wrong sign would lead to
the two possibilities getting interchanged. It would be most reassuring to find some
independent way of choosing the correct picture, and there is.

Indeed, the divisors $D$ under consideration consist of two lines: in one case they
have respectively self-intersection $-2$ and $-1$, and in the other, $0$ and $-1$ (in
both cases, -1 corresponds to the new exceptional divisor). Figure \secblowup\ represents
these two configurations.
\medskip
\centerline{\BoxedEPSF{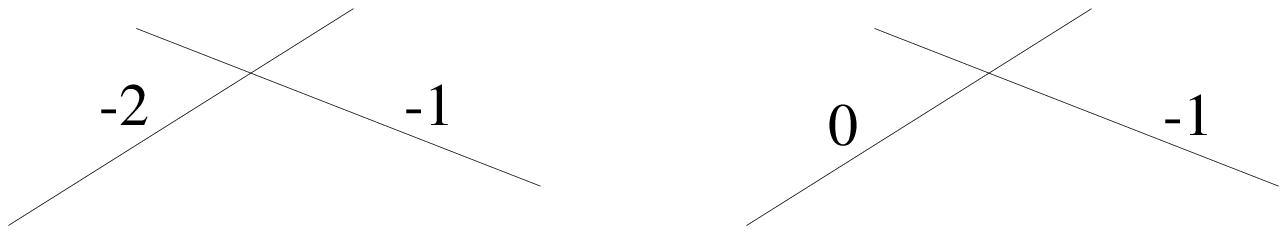 scaled 500}}

\figcap \secblowup {The divisor being blown up.} 

 The self intersection
number $0$ really means that the normal bundle to the proper transform $l_\infty'$ is
trivial. A limit on the torus $p^{-1}(\tilde z)$ of circle orbits coming from the
exceptional divisor is a vector tangent at $z$ to the last exceptional divisor, plus a
normal vector, i.e., a vector tangent to $l_\infty$ which makes a full turn around. 

This can be moved slightly to a circle in $l_\infty'$ turning around $z$, plus
a constant normal vector.  Since the constant normal vector in the small disc
around $z$ can extended to a vector field on
$l_\infty'$ normal to $l_\infty'$, we see the disc on the outside bounded by such an
orbit. So our case, with self intersection of $\l_\infty'$ zero, definitely
corresponds to the proposed picture.
\bigskip
\centerline{\BoxedEPSF{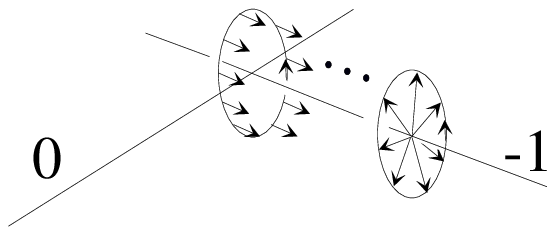 scaled 800}}
\medskip
\longfigcap \whichHopfright {The circle at right represents a normal vector to the
exceptional divisor, turning around the divisor, i.e., an orbit of the circle action. 
This circle can be deformed into the circle at left, on the original line (at infinity),
turning around
$z$, and with a constant normal vector.}


\heading IX. Real oriented blow-ups for complex H\'enon mappings \endheading

Each space $X_{[-N-1,N+1]}$ is obtained from $X_{[-N,N]}$ by a sequence of
blow-ups, first at $\bold q_N$ and $\bold p_{N+1}$ and then at points of the most
recent exceptional divisor; let 
$$
\pi_{N+1,N}:X_{[-N-1,N+1]} \to X_{[-N,N]}
$$
denote the blow-down mapping.  By Propositions \divrightasideal\ and 
\pitildeexistsfordblepoint, the mapping $\pi_{N+1,N}$ induces a projection
$$
\tilde \pi_{N+1,N}:\orb\bigl(X_{[-(N+1),(N+1)]}, D_{[-(N+1),N+1]}\bigr) \to
\orb\bigl(X_{[-N,N]}, D_{[-N,N]}\bigr),
$$
which allows us to consider the projective limit   
$$
\orb(X_\infty, D_\infty) = \varprojlim \Bigl(\orb(X_{[-N,N]}, D_{[-N,N]}); \tilde
\pi_{N+1,N}\Bigr).
$$

There is a canonical inclusion $j_N:\C^2 \to X_{[-N,N]}$ (as in Proposition
\pointsofinfproduct, b), which lifts to $\tilde j_N:\C^2 \to \orb (X_{[-N,N]},
D_{[-N,N]})$  since the real oriented blow-up is taken along
$D_{[-N,N]}=X_{[-N,N]}-\C^2$. These inclusions are compatible with $\tilde
\pi_{N+1,N}$, leading to an inclusion $\tilde j_\infty:\C^2 \to \orb(X_\infty, D_\infty)$.

\theorem \realextensiontheorem a) The mapping $\tilde j_\infty$ is
injective, with dense image, allowing us to think of $\C^2$ as a subset of $\orb(X_\infty,
D_\infty)$.

b) The H\'enon mapping $H:\C^2 \to \C^2$ extends continuously to an automorphism
$$
\orb(H_\infty):  \orb(X_\infty, D_\infty) \to \orb(X_\infty, D_\infty).
$$
\endstat

\proof (a) The injectivity of $\tilde j_\infty$ is clear since all the $\tilde
j_N$ are injective.  Moreover, all the $j_N$ have dense image by Proposition
\pointsofinfproduct, and so do the $\tilde j_N$ since the interior of a manifold
with boundary is dense in the manifold.  The density of the image of $\tilde
j_\infty$ follows immediately, since the topology on the projective limit is
inherited from the topology of the product.  (This also follows from the much
more general Theorem \topofprojlim).

(b) Clearly the H\'enon mapping, i.e., the shift, induces an isomorphism
$$
X_{[-N,N+1]} \to X_{[-(N+1),N]},
$$
hence a homeomorphism 
$$
\orb\bigl(X_{[-N,N+1]},D_{[-N,N+1]}\bigr) \to
\orb\bigl(X_{[-(N+1),N]},D_{[-(N+1),N]}\bigr)
$$
by Proposition \naturalityofrealblowup.
The result will follow since 
$$
\orb(X_\infty,D_\infty) = \varprojlim_{M,N \to \infty}
\orb\bigl(X_{[-N,M]},D_{[-N,M]}\bigr)
$$
when $N$ and $M$ can go to infinity in any way one wants: the pairs $[-N,N]$ are cofinal in
the projective system of pairs $[-N,M]$.
\qed

\remark The existence of $\orb(H_\infty)$ is a bit less trivial than one might expect.  For instance,
there is {\it no well-defined\/} mapping $\orb (\tilde H):\orb(\tilde X, \tilde D)
\to
\orb (\Proj^2,l_\infty)$.  Indeed, although $\tilde D = \tilde H^{-1}(l_\infty)$
set theoretically, it is {\it not true\/} in the sense of ideals.  In the sense of
ideals, $\tilde H^{-1}(l_\infty)$ contains $A',\ B$ and $L_1,\dots, L_d$ with
multiplicity $d-1$, and $L_{d+k}$ with multiplicity $d-k-1$ for $k=0, \dots,
d-3$, and finally $\tilde A$ with multiplicity 1  (for the notation, see Figure
\divatinfty). At the double points where two irreducible components of $\tilde D$
meet and their multiplicities are the same, the map $\tilde H$ extends anyway
because of Example \zerodivs.  But where the two multiplicities are different,
Theorem \naturalityofrealblowup\ does not guarantee that  $\tilde H$ extends, and
in fact it doesn't, as can be verified by explicit computation.
\endremark

The remainder of this section is devoted to understanding the structure of $\orb(X_\infty,D_\infty)$
in detail, as this is equivalent to understanding the dynamics of H\'enon mappings at infinity. 

\theorem \mainthmrealblowup (a) The pair $(\orb(X_\infty, D_\infty), \orb(D_\infty))$ is
homeomorphic to the pair\break $(B_4,S^3)$, the closed 4-ball bounded by the 3-sphere.

b) The mapping $p: \orb(D_\infty) \to D_\infty$ has as its fibers:

\noindent$\bullet$ a circle above ordinary points.

\noindent$\bullet$ a torus above double points,

\noindent$\bullet$ a $d$-adic solenoid $\Sigma^-$ above  $\bold p^\infty$ and a
$d$-adic solenoid $\Sigma^+$ above $\bold q^\infty$.
\endstat

\proof  Part (a) has already been proved in Theorem \topofprojlim.  More
precisely, we have seen that the real oriented blow-up of $\Proj^2$ along the line at
infinity is a 4-ball bounded by a 3-sphere, and so $\varprojlim \orb(X_{[-N,N]},
D_{[-N,N]})$ is also. Indeed, $\varprojlim \orb(X_{[-N,N]},
D_{[-N,N]})$ is obtained by infinitely many times making a blow-up of a surface $X$ at a
point $z$ of a divisor $D$, then taking the real oriented blow-up of the resulting surface
$\tilde X_z$ along the inverse image $\tilde D$ of the divisor $D$.  That is precisely the
situation of Theorem
\topofprojlim.   Moreover, the first two cases of part (b) follow immediately from
Corollary
\descoffibersrealblowups.

The third statement is a bit more delicate. The point $\bold p^\infty$ is represented by
the sequence
$$
\bold p_1 \in X_{[0,0]}, \quad \bold p_2 \in X_{[-1,1]},\quad  \bold p_3 \in
X_{[-2,2]},
\quad \dots
$$
and above this point we see the projective limit of the system of circles
$$
p_0^{-1}(\bold p_1) \leftarrow p_1^{-1}(\bold p_2) \leftarrow p_2^{-1}(\bold p_3)  \dots.
$$
\remark {}We are adding an index to avoid ambiguity,  calling 
$$
p_N:\orb\bigl(X_{[-N,N]}, D_{[-N,N]}\bigr) \to X_{[-N,N]}
$$
 the canonical projection; the added notation is
necessary as the
$\bold p_{N+1}$ can be viewed as points in all $X_{[-M,M]}$ with $M>N$ (or in $X_\infty$), but only in
$X_{[-N,N]}$ are $\bold p_{N+1}$ and $\bold q_N$ simple points of $D_{[-N,N]}$. \endremark

There is a canonical parametrization of $p_N^{-1}(\bold p_{N+1}) \subset \orb(X_{[-N,N]},D_{[-N,N]})$,
obtained from the composition
$$
\underset{\overset \bigcap \to {\orb(X_{[-N,N]},D_{[-N,N]})}} \to {p_N^{-1}(\bold p_{N+1})}
\overset \text{projection to} \to {\underset{\text{$N$th coordinate}}\to {\longrightarrow}}
\underset{\overset \bigcap \to {\orb(\tilde X, \tilde D)}}\to  {p^{-1}(\tilde {\bold p})}
\overset \text{applying $\tilde H$} \to {\longrightarrow}
\underset{\overset \bigcap \to{\orb(\Proj^2, l_\infty)}}\to {p^{-1}(\bold p)} 
\overset {\arg 1/y} \to  {\longrightarrow} 
\R/2\pi \Z.
$$
We will denote this coordinate by $\theta_N$. 

There is also a natural parametrization of $p_N^{-1}(\bold q_{-N}) \subset
\orb(X_{[-N,N]},X_{[-N,N]})$, which is a bit simpler, obtained from the similar composition
$$
\underset{\overset \bigcap \to {\orb(X_{[-N,N]},D_{[-N,N]})}} \to {p_N^{-1}(\bold q_{-N})}
\overset \text{projection to} \to {\underset{\text{$-N$th coordinate}}\to {\longrightarrow}}
\underset{\overset \bigcap \to {\orb(\tilde X, \tilde D)}}\to  {p^{-1}({\bold q'})}
\simeq
\underset{\overset \bigcap \to{\orb(\Proj^2, l_\infty)}}\to {p^{-1}(\bold q)} 
\overset {\arg 1/x} \to  {\longrightarrow} 
\R/2\pi \Z.
$$
We will denote this coordinate by $\phi_N$.

Now part (c) follows from Proposition \computeanglesunderpitilde.

\proposition \computeanglesunderpitilde We have
$$
\tilde\pi_{N+1,N}(\theta_{N+1})= d \theta_{N+1} + \arg a \quad \text{and}\quad
\tilde\pi_{N+1,N}(\phi_{-N-1})= d \phi_{-N-1}.$$
\endstat

\noindent{\bf End of proof of \mainthmrealblowup\ using \computeanglesunderpitilde.} One description
(see \cite {HO1}, Section 3) of the
$d$-adic solenoid $\Sigma_d$ is as
$$
\Sigma_d = \varprojlim (\R/2\pi \Z, \theta \mapsto d\theta).
$$
Clearly the second part of Proposition \computeanglesunderpitilde\ shows that precisely the space
$\Sigma^+$ above $\bold q^\infty$ is canonically the $d$-adic solenoid.  

For the point $\bold p^\infty$,
observe that if we set
$\psi=
\phi +\arg a/(d-1)$, then the mapping $\phi \mapsto d \phi+\arg a$ becomes
$$
\psi \mapsto d\left(\psi-\frac{\arg a}{d-1}\right)+\arg a + \frac {\arg a}{d-1} = d\psi.
$$
Thus the subset $\Sigma^-\subset \orb(X_\infty, D_\infty)$ above $\bold q^\infty$ is also a
$d$-adic solenoid.

\proofof \computeanglesunderpitilde Let us first compute the map 
$$
\underset{\overset \bigcap \to {\orb(\tilde X, \tilde D)}} \to {p^{-1}(\tilde{\bold p})} \to
\underset{\overset \bigcap \to {\orb(\Proj^2,l_\infty)}}\to{p^{-1}(\bold p)}
$$
induced by the blow-down mapping
$\tilde X_H \to \Proj^2$, where both domain and range are identified to $\R/2\pi \Z$:

$\bullet$ ${p^{-1}(\bold p)}$ using $\arg 1/y$;

$\bullet$ ${p^{-1}(\tilde{\bold p})}$ using $(\arg 1/y)\circ \tilde H$.

In Section III, we began by using the coordinates $u= x/y, v=1/y$ near $\bold p$, so we see that $\arg
v$ is our parameter for ${p^{-1}(\bold p)}$. Still in the notation of Section III, 
$\arg X_1$ gives the parametrization of ${p^{-1}(\tilde{\bold p})}$. Formula
(\furthersimpleblowups), for the case $k=d-1$, tells us that the point $\arg X_1=\theta$ of 
${p^{-1}(\tilde{\bold p})}$ is mapped by $\orb(\tilde H)$ to the point where $v$ has argument
$$
\arg \frac {X_1\bigl(X_1^{d-1}Y_{d-1}+a\sum_{j=0}^{d-2}
Q_jX_1^j\bigr)}{X_1^{d-1}Y_{d-1}+a\sum_{j=0}^{d-2} Q_jX_1^j} = \arg X_1.
$$

Thus we must see how the blow-down maps $p^{-1}(\tilde{\bold p})$ to $p^{-1}(\bold p)$, working in
the coordinate $\theta = \arg X_1$ in the domain, and $\theta= \arg v$ in the range, since these
correspond under the H\'enon mapping. 

More precisely, consider the mapping from the circle $p^{-1}(\tilde {\bold p}) \subset
\orb(\tilde X_H, \tilde D)$ to the circle $p^{-1}(\bold p) \subset \orb(\Proj^2,
l_\infty)$. Let $\bold a_0 = \bold p,\ \bold a_1, \dots, \bold a_{2d-2}$
be the successive points at which we performed blow-ups to get from $\Proj^2$ to $\tilde X_H$,
and finally set $\bold a_{2d-1}= \tilde {\bold p}$. The point $\bold a_0$, and the
points $\bold a_d, \dots, \bold a_{2d-1}$ are simple points of the divisor constructed
so far; the points $\bold a_1, \dots, \bold a_{d-1}$ are double points. The inverse images of these
points are parametrized by
$$\eqalign{
\arg v                     \quad&\text {at}\ \bold a_0,\cr
\pvec {\arg X_1}{\arg u}   \quad&\text {at}\ \bold a_1,\cr
\pvec {\arg X_1}{\arg X_k} \quad&\text {at}\ \bold a_k,\quad k=2, \dots, d-1,\cr
\arg X_1                   \quad&\text {at}\ \bold a_k,\quad k=d, \dots, 2d-1.}\peqno \argumentsconseq
$$

In these coordinates, the blow-down mapping $\tilde X_H \to \Proj^2$ induces the
composition 
$$\eqalign{
\underset {\bold a_{2d-1}} \to {\vphantom{\bvec 11}\theta} \mapsto \dots \mapsto &\underset
{\bold a_{d}}
\to {\vphantom{\pvec 11}\theta}\mapsto \underset {\bold a_{d-1}} \to {\bvec
{\theta}{\theta+\arg a}}\mapsto 
\underset {\bold a_{d-2}} \to {\pvec {\theta}{2\theta+\arg a}}\mapsto \dots\cr
&\mapsto     
\underset {\bold a_1} \to {\pvec {\theta}{(d-1)\theta+\arg a}}\mapsto
\underset {\bold a_0} \to {\vphantom{\bvec 11}d\theta +\arg a}.}\peqno \argumentsaddeq
$$
All of these are straightforward applications of Propositions \compsimpblowupontorus\ and
\compmapontorofdblept.  Let us spell out the mapping which takes $\bold a_d$ to $\bold a_{d-1}$.  In
that case we are taking the circle above the point $X_1=0, X_d=a$, parametrized by $\arg X_1$, to the
torus above $X_1=0, X_{d-1}=0$.  The blow-down mapping is 
$$
\pvec {X_1}{X_{d}} \mapsto \pvec {X_1}{X_{d-1}} = \pvec {X_1}{X_1X_d},
$$
and in particular the circle $X_1=\rho e^{i\theta},\ X_d=a$ is mapped to the circle
$X_1=\rho e^{i\theta}, X_{d-1} =\rho a e^{i\theta}$.  If we let $\rho \to 0$ and remember only the
arguments, we get the desired formula. 

This almost proves the first part of Proposition \computeanglesunderpitilde; by definition, the
mapping\break ${p_N^{-1}(\bold p_{N+1})} \to {p_{N-1}^{-1}(\bold p_{N})}$ is precisely the mapping
above, the domain and range being identified to $\R/2\pi \Z$ by $H^{\circ N}$ and
$H^{\circ (N-1)}$ respectively.

There are two ways of approaching the first part of Proposition \computeanglesunderpitilde: to make the
sequence of blow-ups at $\bold q$, repeating the material of section III to make $H^{-1}$ well-defined,
or to make a change of variables to make $H^{-1}$ conjugate to a H\'enon mapping, for a new
polynomial $p_1$ and and a new Jacobian $a_1$, of course.  Remember that we used the fact that
the polynomial $p$ is monic in Section III, so $p_1$ must also be monic.  We invite the
reader to show that in the variables
$x_1, y_1$, where
$\zeta x_1= y$ and $\zeta y_1=x$ with $\zeta^{d-1}=a$, we have
$$
H^{-1}: \pvec {x_1}{y_1} \mapsto \pvec {p_1(x_1)- y_1/a}{x_1}
$$
with $p_1$ monic.  Thus, in these coordinates the blow-down takes 
$$
\arg (1/y_1) \to d\arg (1/y_1) + \arg (1/a).
$$
We invite the reader to check that this means exactly that in the original variables,
 the formula of
 Proposition \computeanglesunderpitilde\ is satisfied.
\qed

\heading{X. The topology of $\orb(X_\infty,D_\infty)$} \endheading
 
We now want to describe the mapping $\orb(H_\infty)$ more precisely.  A first statement
says that appropriate restrictions of $\orb(H_\infty)$ and $\orb(H_\infty)^{-1}$ are
solenoidal (\cite{HO1}, Section 3).  

\remark {} Actually, the definition of solenoidal which appears in \cite{HO1} isn't
quite the right one in our setting, because there the mappings are differentiable, and the notion of
expanding and contracting is with respect to this structure; $\orb(D_\infty)$ doesn't have a
natural differentiable structure, especially on the solenoids $\Sigma^\pm$ which are naturally
accumulations of corners. However, the main result we will be using is Theorem 3.1 of \cite{HO1} which
uses none of this structure (and in the main case of interest, we will construct the mapping
explicitly anyway, in Proposition \solenoidalprop).  The construction of
Theorem \topofprojlim\ provides the necessary contraction in the fiber directions. \endremark

   Let us denote by
$T_{\bold p_i}$ and
$T_{\bold q_i}$ the tori
$\tilde p^{-1}(\bold p_i)$ and $\tilde p^{-1}(\bold q_i)$.  Each of these separates
$\orb(D_\infty)$ into two pieces.  We will denote by $T_{\bold p_i}^+$ the one which
contains the attractive solenoid $\Sigma^+$, and by $T_{\bold p_i}^-$ the one which
contains the repelling solenoid $\Sigma^-$, and similarly for $\bold q_i$.
 
\proposition \solenoidalprop a) The mapping $\orb(H_\infty)$ maps $T_{\bold p_i}^+$ (resp.
$T_{\bold q_i}^+$) into itself, and is solenoidal of degree $d$.

b) The mapping $\orb(H_\infty)^{-1}$ maps $T_{\bold p_i}^-$ (resp.
$T_{\bold q_i}^-$) into itself, and is solenoidal of degree $d$.
\endstat

\proof  We need to examine carefully the sequence of
blow-ups which makes $H$ well-defined, to understand how the tori corresponding to the
double points of $D_\infty$ are embedded in the 3-sphere $\orb(D_\infty)$. Recall that we called
$\bold a_0 = \bold p,\ \bold a_1, \dots, \bold a_{2d-2}$
the successive points at which we performed blow-ups to get from $\Proj^2$ to $\tilde X_H$,
and finally set $\bold a_{2d-1}= \tilde {\bold p}$. 

In Section VIII, we showed how to start, creating a torus $T_{\bold a_0}$ in the 3-sphere,
separating the solid torus corresponding to $l_\infty$ from the solid torus corresponding to the
first exceptional divisor (which will eventually be the irreducible component $B$ of $\tilde D$). See
Figure \rightsecondblowup\ to understand how these solid tori are placed after stereographic projection;
the words ``inner'' and ``outer'' will refer to this picture.

The next $d-1$ blow-ups are fairly easy to understand, now that we have started
right.  We thicken the torus $T_{\bold a_0}$, creating an inner torus $T_{\bold a_1}$ and
an outer torus $T_{\bold p'}=T_{\bold p_0}$ (which we can call by its final name, since it will
not be affected by further blow-ups).  Then thicken the inner torus $T_{\bold a_1}$, creating an
inner torus $T_{\bold a_2}$, and one which corresponds to $L_1\cap L_2$ (between 
$T_{\bold a_2}$ and $T_{\bold p_0}$). Then thicken the inner torus $T_{\bold a_2}$ again,
$d-1$ times in all. The inside of the torus $T_{\bold p_0}$ is
$T^-_{\bold p_0}$. See Figure \firstdblowuppic\ for the case $d=3$.

The circle orbits fibering the regions between the successive tori are contained in
$T^-_{\bold p_0}$ with the innermost torus (corresponding to the component $B$) removed, which is
a space with homology
$\Z^2$, generated by $a$ and $b$ (see Figure \rightsecondblowup). At the first thickening, an
oriented circle orbit between the two tori is in the homology class
$-a-2b$, at the next the new thickened  torus is fibered by curves with homology
class $-a-3b$, etc., ending up with a thickened torus fibered by circles in the
homology class
$-a-db$, and an inner solid torus (corresponding to
$B$) with fibers in the homology class $-b$.

In summary, after the first $d$ blow-ups, we have an inner solid torus with fibers in the
homology class $-b$, then a succession of thickened tori with fibers in the
 homology classes 
$$
-a-db,\quad  -a-(d-1)b,\quad  \dots, \quad -a-2b, \peqno \homologyclassesoffirstblowups
$$ and finally the
region $T^+_{\bold p_0}$ corresponding to $\l_\infty$, fibered (by the old Hopf fibers) in
the homology class $-a-b$. See Figures \firstdblowuppic\ and \toriasrotpics.

We must now make $d-1$ more blow-ups of ordinary points.  The first of these can be
realized by thickening a circle orbit in the region corresponding to $L_{d-1}$,
which is fibered by circles in the homology class $-a-db$.  We will then thicken a circle
inside this torus, which we may take to be the ``core circle'', and repeat this
$d-2$ times. The final torus we create this way is $T_{\bold q_1}$.  All the solid tori
are thickenings of the original circle, and hence in the homology class $-a-db$.

We need to start making the second series of blow-ups, as we don't
yet have the torus $T_{\bold p_1}$.  So thicken a fiber inside $ T_{\bold q_1}^-$,
creating a solid torus still in the homology class $-a-db$, and thicken it again $d-1$
times; the outermost torus of the series just created is $T_{\bold p_1}$. We will not
need to describe the further blow-ups.

\centerline{\BoxedEPSF{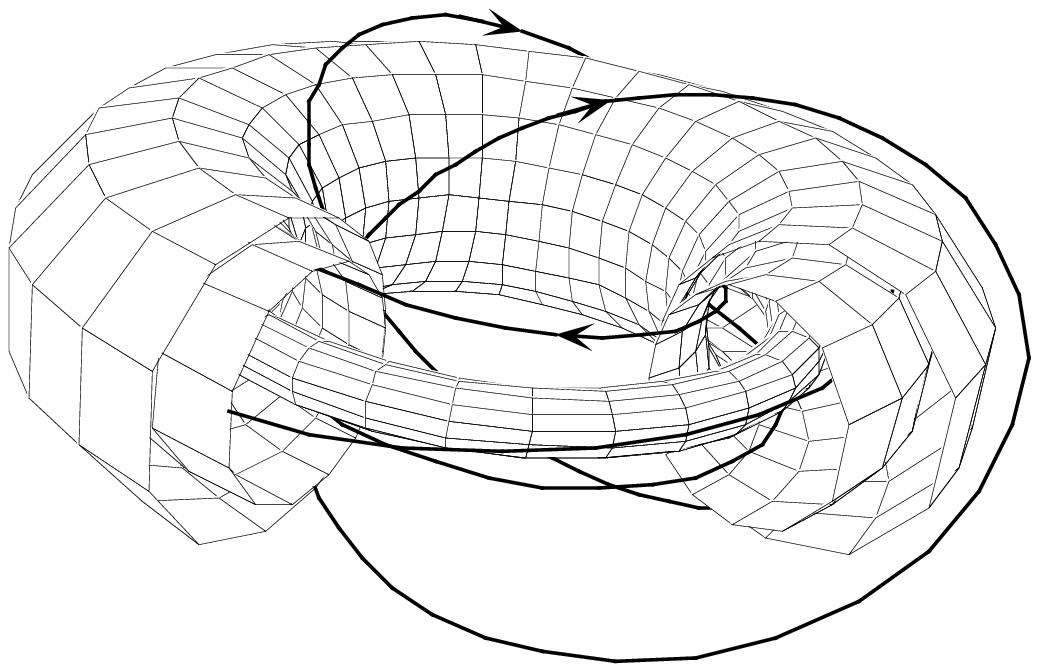 scaled 600}}

\longfigcap \firstdblowuppic {The picture above corresponds to the situation after $d$
blow-ups when $d=3$, and after stereographic projection.  The outside corresponds to the
real oriented blow-up of the line at infinity (now $A'$), with the Hopf circle action,
as shown.  The inner torus corresponds to $B$, with the circle action where the orbits are
not linked. The region between the inner torus and the next corresponds to $L_{d-1}$; that
is where all the further action will take place, thickening a circle orbit (represented
on the drawing, and going around 3 times in one direction as it turns once in the other
direction).  The region between the outer torus and the next corresponds to $L_1$; the
circle action there has orbits which turn twice in one direction as they turn once in
the other; we haven't drawn them to keep the drawing simpler.}

\bigskip

\centerline{\BoxedEPSF{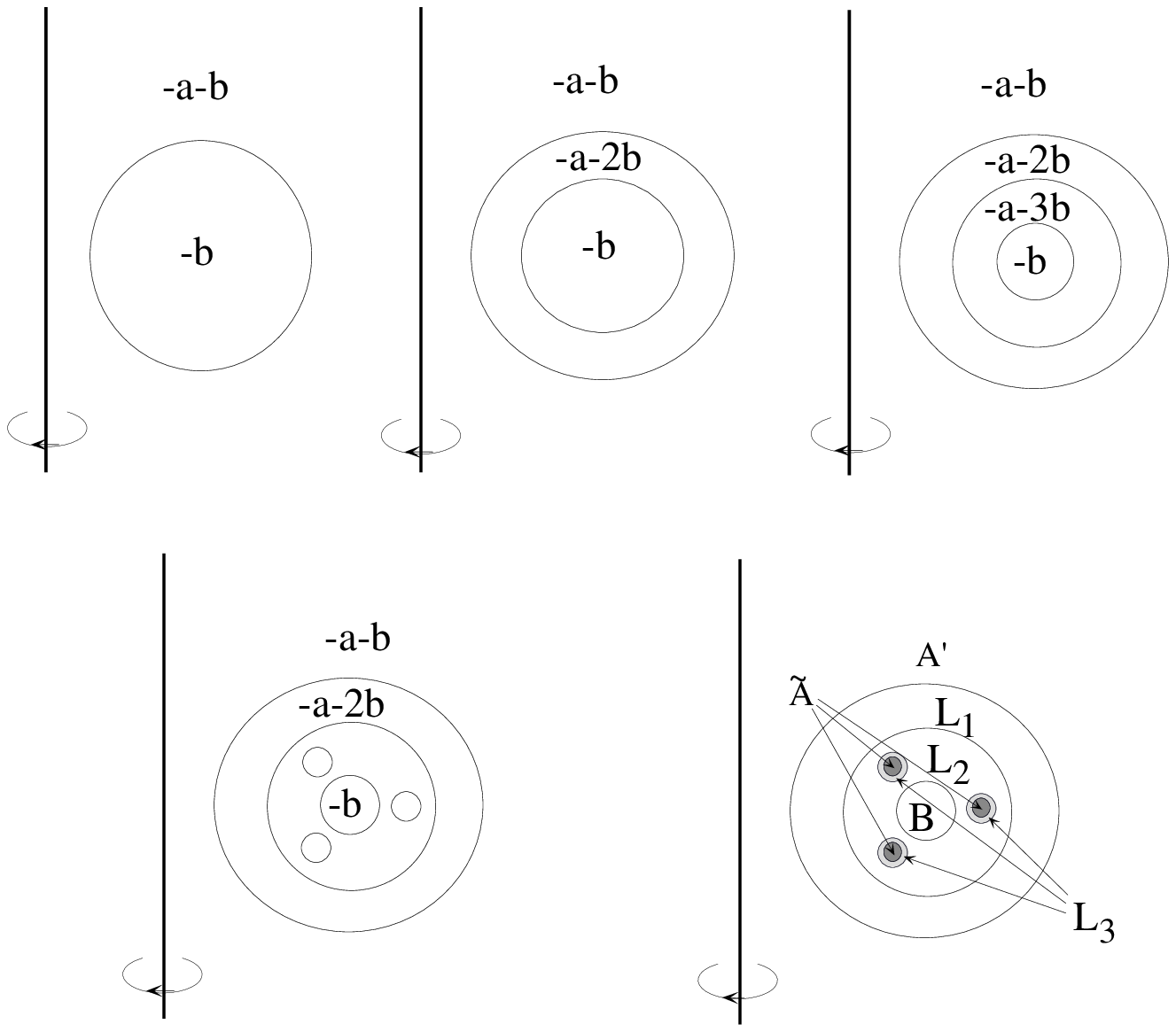 scaled 600}}

\longfigcap \toriasrotpics {You can almost imagine constructing the pattern of tori in $S^3$ by
rotating the figure above around the $z$-axis (shown as a heavy line).  The case represented
corresponds to $d=3$. The ``almost'' is because the small circles are not actually rotated: as they
turn around the $z$-axis they also turn in their annulus, so as to connect up and together form a
single torus, as shown in Figure \firstdblowuppic.}  

Moreover, $\orb(H_\infty)^{-1}$ is a homeomorphism which maps the solid torus $
T^-_{\bold p_0}$  to the solid torus $ T^-_{\bold p_1}$.  We
claim that as a map
$ T^-_{\bold p_0}\to T^-_{\bold p_1}$ it is conjugate to the mapping $\tau_{d,0}$
as defined in
(\cite{HO1}, Section 3).  Recall that  the mapping $\tau_{d,k}:S^1\times D \to S^1
\times D$, where $S^1, D \subset \C$ are respectively the unit circle and the  disc of
radius 2, is given by the formula
$$
\tau_{d,k}(\zeta,z) = \bigl(\zeta^d, \zeta+\epsilon z \zeta^{k+1-d}\bigr).
$$

A ``solenoidal'' map $S^1\times D \to S^1
\times D$ is one which is expanding in the circle direction and
contracting in the disc direction.  Certainly $\orb(H_\infty)$ expands in the circle
direction; in fact it is $\beta \mapsto d\beta$  (this is the coordinate $\beta$ of the stereographic
projection). By choosing our thickenings sufficiently small, the mapping will be contracting on the discs.
\qed

\corollary \solenoidalmapsconjtotau The mappings $\orb(H_\infty):T_{\bold p_i}^+ \to
T_{\bold p_i}^+$ and $\orb(H_\infty)^{-1}:T_{\bold p_i}^- \to
T_{\bold p_i}^-$ are conjugate to $\tau_{d,0}$.
\endstat

\proof
Theorem 3.11 of \cite{HO1} asserts that every unbraided solenoidal mapping from a solid
torus to itself is conjugate to precisely one of the $\tau_{d,k}$, and Propositions 4.1 and
4.6 assert that only $\tau_{d,0}$ extends to the 3-sphere.  Indeed, Proposition 4.6 asserts
that when $k\ne 0$, the forward images of $\tau_{d,k}$ are knotted; but no homeomorphism of
$S^3$ can map an unknotted solid torus to a knotted one.  Since our map
$\orb(H_\infty)$ extends to the 3-sphere $\orb(D_\infty)$, it is enough to
prove that
$$
\orb(H_\infty)^{-1}: T^-_{\bold p_1}\to T^-_{\bold p_1}
$$
is solenoidal and unbraided.  The unbraided part follows from the homotopy class
$-a-db$.
\qed

By Theorem 3.1 of \cite{HO1}, there are maps $\pi_i^+:T_{\bold p_i}^+ \to \R/2\pi \Z$ and
$\pi_i^-:T_{\bold p_i}^- \to \R/2\pi \Z$ such
that the diagrams 
$$
\CD
T^+_{\bold p_i}  @>{\orb(H_\infty)}>> T^+_{\bold p_i} \\ 
@V{\pi_i^+}VV   @V{\pi_i^+}VV     \\ 
\R/2\pi \Z  @>{\theta \mapsto d\theta}>>  \R/2\pi \Z
\endCD
\qquad\qquad
\CD
T^-_{\bold p_i}  @>{\orb(H_\infty)^{-1}}>> T^-_{\bold p_i} \\ 
@V{\pi_i^-}VV   @V{\pi_i^-}VV     \\ 
\R/2\pi \Z  @>{\theta \mapsto d\theta}>>  \R/2\pi \Z
\endCD
$$
commute.

In our case, these functions $\pi_i^+$ and $\pi_i^-$ can be computed explicitly; they are given
by  proposition \computepiplusandminus.  Before stating this proposition, notice
that 

$\bullet$ $T_{\bold q_i}^+$ is the set of $ x\in \orb(X_\infty, D_\infty)$ such
that
$\bigl(p_\infty(x)\bigr)_i=\bold q'$;

$\bullet$ $T_{\bold q_i}^-$ is the set of $ x\in \orb(X_\infty, D_\infty)$ such that
$\bigl(p_\infty(x)\bigr)_{i-1}\in \tilde A$.

$\bullet$ $T_{\bold p_i}^-$ is the set of $ x\in \orb(X_\infty, D_\infty)$ such that
$\bigl(p_\infty(x)\bigr)_{i-1}=\tilde{\bold p}$;

$\bullet$ $T_{\bold p_i}^+$ is the set of $ x\in \orb(X_\infty, D_\infty)$ such that
$\bigl(p_\infty(x)\bigr)_i\in A'$.

There is a natural blow-down mapping  $Q_i:X_\infty \to X_{[i,\infty)}$, which blows $D_\infty$
down onto
$D_{[i,\infty)}$.
 
Since $Q_i$ is a blow-down, it induces a mapping 
$$
\orb(Q_i):\orb(X_\infty,D_\infty) \to \orb\bigl(X_{[i,\infty)},D_{[i,\infty)}\bigr).
$$
More precisely, the blow-downs $X_{[-N,N]} \to X_{[i,N]}, \ i<N$  induce mappings on the
real oriented blow-ups $\orb(X_{[-N,N]},X_{[-N,N]}) \to \orb(X_{[i,N]},D_{[i,N]})$; to construct
$Q_i$, we must pass to the projective limit.

The mapping $Q_i$  maps
$T_{\bold q_i}^+$ to the circle above
$\bold q_i$ in $\orb(D_{[i,\infty]}$.  This circle is canonically parametrized, by $\phi_i$. Let us denote
$$
\Phi_i= \phi_i \circ \orb(Q_i)|_{T_{\bold q_i}^+}: T_{\bold q_i}^+ \to \R/2\pi \Z
$$
 the composition.

Exactly analogously, there is a natural blow-down  $P_i:X_\infty \to X_{(-\infty,i]}$ , which blows
$D_\infty$ down onto $D_{(-\infty,i]}$.
 
Again, since $P_i$ is a projective limit of blow-downs, it induces a mapping 
$$
\orb(P_i):\orb(X_\infty,D_\infty) \to \orb\bigl(X_{(-\infty,i]},D_{(-\infty,i]}\bigr)
$$
which maps $T_{\bold p_i}^-$ to the circle above $\bold p_i$.  This circle is canonically
parametrized, by $\psi_i$ (remember that $\psi_i = \theta_i+\arg(a/(d-1))$.  Let us denote
$$
\Psi_i= \psi_i \circ \orb(P_i)|_{T_{\bold p_i}^-}: T_{\bold p_i}^- \to \R/2\pi \Z
$$
 the composition.

\proposition \computepiplusandminus  (a) We may choose $\pi_i^+= \Phi_i$.

(b) We may choose $\pi_i^-= \Psi_i$.
\endstat

\proof (a) Clearly $\Phi_{i-1}(\orb(H_\infty)(x))= \Phi_i(x)$ when $x \in T_{\bold q_i}^+$, as
the left-hand side is just the right-hand side shifted one to the left. Moreover,
Proposition
\computeanglesunderpitilde\ says that $\Phi_i(y) = d \Phi_{i-1}(y)$ for $y \in T_{\bold
q_{i-1}}^+$. So
$$
\Phi_i\bigl(\orb(H_\infty)(x)\bigr)= d \Phi_{i-1}\bigl(\orb(H_\infty)(x)\bigr) = d
\Phi_i(x).
$$
The argument for part (b) is similar.   
\qed

\remark {} We should discuss what happened to the $d-1$ choices of $\pi^+$ and $\pi^-$.  For $\pi^+$,
our particular choice was given by the coordinate system in $\C^2$, because ultimately, $\psi = \arg
1/x$.  If we conjugate a H\'enon mapping by setting
$x_1=\zeta x, y_1= \zeta y$ where $\zeta^{d-1}=1$, it is easy to show that the H\'enon mapping
remains of the same form  (the polynomial remains monic and the number $a$ is not
changed). For $\pi^-$, we don't actually have a canonical choice. The coordinate $\theta$, which
ultimately comes from $\arg 1/y$, is canonical, but $\psi=\theta+\arg(a/(d-1))$ is exactly ambiguous by
a $d-1$ root of 1, as one would expect.  \endremark

The point of these computations is that since $\pi^+$ and $\pi^-$ are conjugacy
invariants (up to the ambiguity above) of the mapping $\orb(H_\infty)$, we can use them to find a condition
for when the restrictions of $\orb((H_1)_\infty)$ and $\orb((H_2)_\infty)$ to the spheres at infinity are
conjugate, where $H_1$ and $H_2$ are H\'enon mappings with corresponding polynomials $p_1$ and $p_2$,
and Jacobians $a_1$ and $a_2$. In order to pin down our result, we need to know something about
toroidal decompositions of 3-manifolds.

\subheading {Toroidal decompositions}  The 3-manifold $\orb(D_\infty)-(\Sigma^+ \cup \Sigma^-)$ has an
interesting toroidal decomposition.  We will give the definitions and basic properties, largely due to
Jaco, Johannson, Shalen and Waldhausen. Our sources for this material are \cite{Hemp} and especially
\cite{Hat2}.

Let $M$ be an orientable irreducible 3-manifold with boundary.  A properly embedded surface $S \subset M$ is
{\it incompressible\/} if for any closed embedded disc in $D \subset M$ with $\partial D \subset S$,
there is a disc $D' \subset S$ with $\partial D = \partial D'$. The manifold $M$ is atoroidal if each
incompressible torus is isotopic to a boundary component.

 Let $D \subset \C$ be the open unit disc. A {\it Seifert manifold\/} is a 3-dimensional manifold,
foliated by circles, such that each leaf has a neighborhood homeomorphic to the quotient of $D \times
[0,1]$ by the equivalence relation which identifies $(0,z)$ to $(1, e^{2\pi i p/q} z)$ for some rational
number $p/q$, with the foliation induced by the lines $\{z\} \times [0,1]$. The set of leaves is then a
surface with boundary $\Omega$, and the canonical mapping $M \to \Omega$ is referred to as a {\it Seifert
fibration}. This is a locally trivial fibration over the subset $\Omega' \subset \Omega$ corresponding to
the regular leaves; the singular leaves (like the one corresponding to $z=0$ in the model) correspond to
the discrete set $\Omega-\Omega'$.

The key results for us are the following: 

\theorem \Hatcherthreepointthree Let $M$ be a 3-dimensional compact orientable manifold with
boundary. Then there exists a collection of disjoint incompressible tori $T_i \subset M$ such that
each component of $M-\cup_iT_i$ is either atoroidal or a Seifert manifold, and a minimal such
collection is unique up to isotopy.
\endstat

This is exactly Theorem 3.3 of \cite{Hat2}.

\theorem \seifertfibrationunique A Seifert manifold with at least two boundary
components has a unique Seifert fibration up to isomorphism.
\endstat

This follows immediately from Theorem 4.3 of \cite{Hat2}.  Indeed, Hatcher shows that the Seifert
fibration is unique except for a list of exceptions, and all these exceptions have 1 or 0 boundary
components.

\theorem \surfacesinseifert Let $f:M \to \Omega$ be a Seifert fibration, and let
$\Omega' \subset \Omega$ be the complement of the points corresponding to singular
fibers. Suppose that
$M$ is connected and that $\partial M \ne \emptyset$.  Then every incompressible surfaces in $M$
without boundary is isotopic to a surface of the form $f^{-1}(\gamma)$ for some curve $\gamma \subset
\Omega'$, and the isotopy classes of such surfaces correspond exactly to the isotopy classes of such
curves.
\endstat

This follows from Proposition 3.5 of \cite{Hat2}.  Hatcher proves that every incompressible and
boundary-incompressible surface is isotopic to either a vertical or a horizontal surface. 
Horizontal surfaces have non-empty boundary, and surfaces without boundary are vacuously
boundary-incompressible. So our surfaces are isotopic to vertical surfaces, i.e., surfaces of the
form $f^{-1}(\gamma)$.   

We will be interested in applying these notions to the manifold $\orb(D_\infty)-(\Sigma^+ \cup
\Sigma^-)$, which comes with a family of tori $T_{\bold p_i}$, which we will see are
incompressible.

\theorem \toroidaldecompofsphere a) The tori $T_{\bold p_i}\subset \orb(D_\infty)-(\Sigma^+
\cup \Sigma^-)$ are incompressible.

b) Every incompressible torus in
$\orb(D_\infty)-(\Sigma^+
\cup
\Sigma^-)$ is isotopic to exactly one of the $T_{\bold p_i}$.
\endstat

\proof a) The homology $H_1(\orb(D_\infty)-(\Sigma^+ \cup \Sigma^-))$ is isomorphic to $\Z[1/d]
\oplus \Z[1/d]$. This is proved for instance using the Alexander duality theorem (\cite{Spa}, 6.2, Thm. 16),
which asserts that
$$
H_1\left(\orb(D_\infty)-(\Sigma^+ \cup \Sigma^-)\right) = H^1(\Sigma^+)\oplus
H^1(\Sigma^-).
$$
The cohomology above is \v Cech cohomology (isomorphic to Alexander-Spanier cohomology), given by
the inductive limit of the singular cohomology of a basis of neighborhoods of $\Sigma^\pm$. Using the system
of neighborhoods
$T_{\bold p_{-i}}^+$ of
$\Sigma^+$, we see that 
$$
H^1(\Sigma^+) = \varinjlim (\Z,n \mapsto dn) = \Z[1/d].
$$
(This is similar to but simpler than Example \indlimextw.)

An isomorphism is specified by sending the generators $a$, $b$(see Figure \rightsecondblowup)
of $H_1T_{\bold p_0}$ to $(0,1)$ and $(1,0)$ in $\Z[1/d] \oplus
\Z[1/d]$. Under this isomorphism, the corresponding generators of
$H_1(T_{\bold p_i})$ are sent to $(0,d^{-i})$ and $(d^{i},0)$.  In particular, the inclusion is
injective on the homology of such a torus.  If a disc in $\orb(D_\infty)-(\Sigma^+ \cup \Sigma^-)$
bounds a disc in $T_{\bold p_i}$, then the homology class of this boundary is zero in
$H_1(\orb(D_\infty)-(\Sigma^+ \cup \Sigma^-))$, hence also in
$H_1(T_{\bold p_i})$, so the curve bounds a disc in the torus, since any simple closed curve in a
torus  which is trivial in the homology bounds a disc. (This is not true on surfaces
of higher genus, which is why incompressibility is not defined using injectivity of
the inclusion on homology.) This is the definition of an incompressible torus.

b) Let $T'$ be such a torus.  It is contained in the compact manifold $T_{\bold p_i}^- \cap 
T_{\bold p_j}^+$ for $i$ sufficiently small and $j$ sufficiently large, which allows
us to apply Theorem
\Hatcherthreepointthree, where $M$ must be compact.  So it is enough to prove that $T_{\bold p_{i+1}},
\dots, T_{\bold p_{j-1}}$ is a minimal family of incompressible tori in
$T_{\bold p_i}^- \cap  T_{\bold p_j}^+$ such that the components of the boundary are atoroidal or
Seifert manifolds.

First observe that the components of
$$
\left(T_{\bold p_i}^- \cap  T_{\bold p_j}^+\right)-\left( \bigcup_{n=i+1}^{j-1}T_{\bold p_{n}}\right)
$$
are in fact both atoroidal and Seifert manifolds.  Indeed, the region $T_{\bold p_i}^- \cap  T_{\bold
p_{i+1}}^+$ is homeomorphic to the region $M_i$ bounded by the tori corresponding to
$L_{i,d-2} \cap L_{i,d-1}$ and $L_{i,d} \cap L_{i,d-1}$. This region contains the solid torus
corresponding to $B_i$, but that torus can be collapsed onto a circle without changing the
homeomorphism type; call $M_i'$ the resulting manifold.

The manifold $M_i'$ is fibered by the natural circle action, and the circle corresponding to $B_i$
becomes a singular circle of type $(1,d)$.  Thus $T_{\bold p_i}^- \cap  T_{\bold p_{i+1}}^+$ is
a Seifert manifold; let $f_i:T_{\bold p_i}^- \cap  T_{\bold p_{i+1}}^+ \to \Omega_i$ be the
corresponding projection to the set of leaves (the base).  It is also atoroidal, by
Theorem \surfacesinseifert, since
$\Omega_i$ is an annulus with one distinguished point corresponding to the unique singular fiber.  This is
seen as follows. In Figure \toriasrotpics, the manifold $M_i$ corresponds to the annulus between the center
circle and the next, with the $d$ small circles removed.
\bigskip
\centerline{\BoxedEPSF{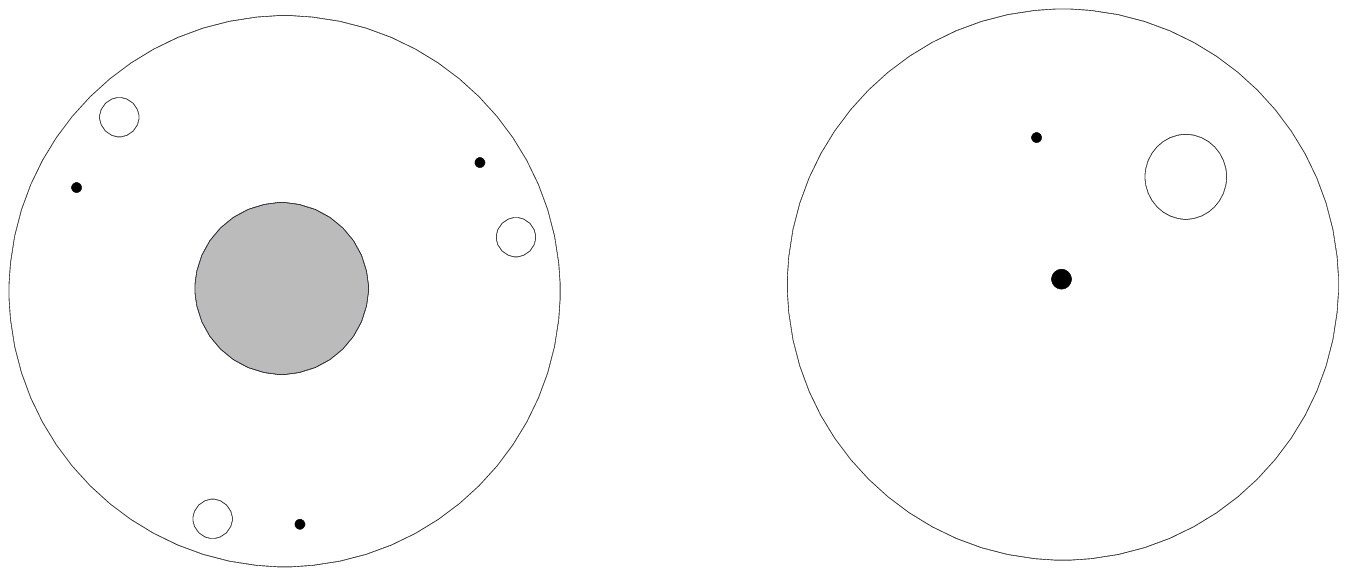 scaled 600}}

\longfigcap \baseSeifertfig {On the right, we have repeated the relevant part of Figure
\toriasrotpics, showing the annulus corresponding to $M_i$ when $d=3$; the three dots
represent the intersection of the plane of the figure with a circle orbit in $M_i$. 
The right-hand side represents the first, after applying $z
\mapsto z^d$ and collapsing the central disc to a point. Every point on the right corresponds to a unique
circle orbit, i.e., the base is an annulus with a single singular fiber (corresponding to the central
point).} 

If you parametrize the disc by $z$, and compose $z \mapsto z^d$ with a collapse of the central disc to a
point (corresponding to collapsing the solid torus corresponding to $B$ to a circle), you manufacture a
space which corresponds exactly to the set of leaves. See Figure \baseSeifertfig.

Any simple closed curve on an annulus with a puncture is homotopic to a point or to a boundary component,
so there are no incompressible tori in $M'_i$ by Theorem \surfacesinseifert. 

Now we need to show that our family $T_{\bold p_{i+1}}, \dots, T_{\bold p_{j-1}}$ is minimal.  It is
clearly enough to show that $T_{\bold p_i}^- \cap  T_{\bold p_{i+2}}^+$ is neither atoroidal nor
Seifert.  It clearly isn't atoroidal, since it contains $T_{\bold p_{i+1}}$, so we must show it isn't
Seifert.  Suppose $f:T_{\bold p_i}^- \cap  T_{\bold p_{i+2}}^+\to \Omega$ is a Seifert fibration,
where
$\Omega$ is some surface. The surface $\Omega$ must have a boundary consisting of two components, but
otherwise we don't know much about it. By Theorem \surfacesinseifert, there is a curve $\gamma \subset
\Omega'$, where $\Omega'$ is the complement of the projections of the singular fibers, such that $T_{\bold
p_i}$ is isotopic to $T_i'= f^{-1}(\gamma)$; in particular, the restriction of $f$ to the components of
$\left(T_{\bold p_i}^- \cap  T_{\bold p_{i+2}}^+\right)-T'_i$ is a Seifert fibration.  

But each of these is already a Seifert fibration, and in fact in a unique way by Theorem
\seifertfibrationunique.  It is then enough to show that the fibers of $f_i$ and $f_{i+1}$ on the
torus $T_{\bold p_i}$ which is the intersection of their domains are not homotopic curves; since
they should both be homotopic to the fibers of $f$, this contradicts the existence of such an $f$.

In the basis $a,b$ for $H_1(T_{\bold p_i})$, we have seen that a fiber of $f_0$  has the homology
class $-a-db$.  We claim that a fiber of $f_1$ has the homology class $-da-b$; knowing this will
end the proof.

We need to repeat the construction of Proposition \solenoidalprop, to understand the sequence of
tori corresponding to the block $B_1$ of $D_\infty$.  Take the 3-sphere, with the sequence of
real oriented blow-ups corresponding to $\orb(\tilde D)$, as represented in Figure \firstdblowuppic.

The fiber above $\bold q$ is a circle orbit outside the torus corresponding to $\bold q'$, which we
may take to be the $z$-axis. Thicken this torus; the fibers inside the thickened torus
will now have the homology class $-a-b+b=-a$ by Theorem \blowupofordpoint.  Indeed, if
you choose a small disc transverse to the $z$-axis, it projects to $\tilde D$ (in
fact, to a neighborhood of $\bold q'$ in the line at infinity $A'$), and the induced
orientation gives its boundary the orientation $+b$.  Now, when we make $d-1$ more
blow-ups, always of double points, the regions between the tori created have circle
orbits in the classes $(-2a-b,\ \dots, -da-b)$.  The region foliated by curves in the
class
$-da-b$ corresponds to the 3-manifold $M_1$.
\qed

\centerline{\BoxedEPSF{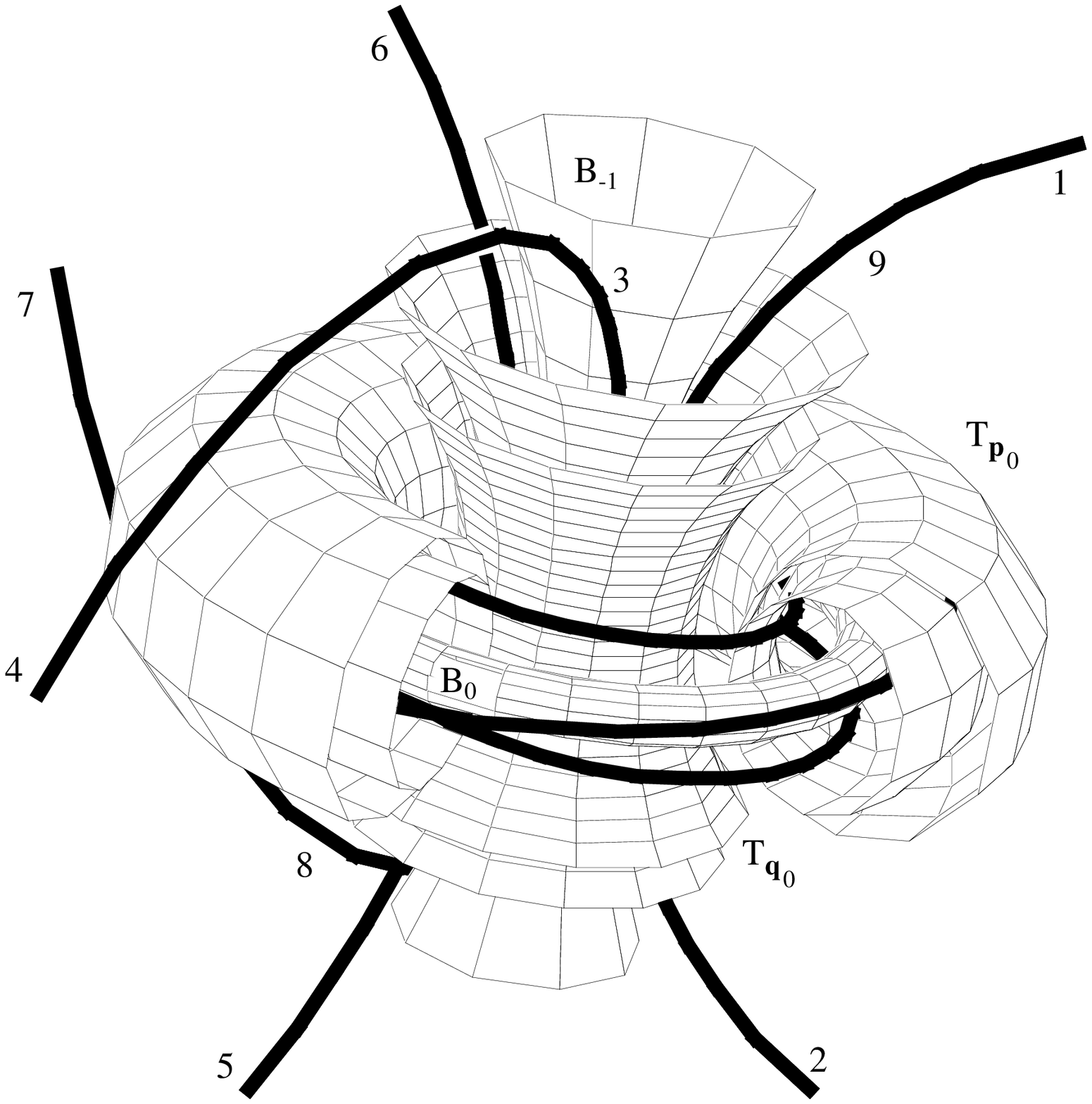 scaled 700}}

\longfigcap \blowupsforandback {The configuration of tori in the 3-sphere at infinity between the torus
corresponding to $B_0$ and $B_{-1}$, in the case $d=3$.  The torus corresponding to $L_{-1,1}\cap
L_{-1,2}$, shown as a heavy curve,  winds three times around the torus corresponding to
$B_0$, in the figure a small thickening of the unit circle in the $x,y$-plane.  The
torus corresponding to $L_{0,2}\cap L_{1,3}$, also shown as a heavy curve, winds three
times around the torus $B_{-1}$, represented in the figure as a thickening of the
$z$-axis.}

\subheading {Conjugacy invariants of $\orb(H_\infty)$}  

Let $H_1$ and $H_2$ be H\'enon mappings. We will give a necessary condition for when
the restrictions of $\orb((H_1)_\infty)$ and $\orb((H_2)_\infty)$ to the spheres at infinity
$\orb(D_{1,\infty})$ and $\orb(D_{2,\infty})$ are conjugate. To (sort of) lighten notation, we will call
these restrictions  $\orb((H_i)'_\infty)$.  

\theorem \conjugacyofexthens In order for $\orb((H_1)'_\infty)$ and $\orb((H_2)'_\infty)$
to be topologically conjugate, it is necessary that they have the same degree, and that
$\arg a_1 \equiv \arg a_2 \mod 2\pi/(d-1)$.
\endstat

\proof  That they must have the same degree is clear by counting fixed points in the solenoids.

We will first
investigate the ``critical locus'' of the map 
$$(\pi_1^+, \pi_0^-):T_{\bold p_0}^- \cap T_{\bold p_1}^+\to (\R/2\pi \Z)^2.
$$
\remark {}The notion of ``critical locus'' isn't quite right: $\orb(D_\infty)$ isn't
naturally a differentiable manifold, it is naturally a {\it manifold with corners},
almost an object of the piecewise-linear category. What we will find is more PL than
$C^1$: it will turn out that the set of points which have  neighborhoods on which
$\pi_1^+$ and
$\pi_0^-$ differ by a constant is non-empty.  In our setting, the critical locus will
be the closure of this open set by definition. This is a much stronger notion of
``critical'' than one would expect: generically for a differentiable mapping from a
3-dimensional manifold to a surface, the critical locus should be a curve.
\endremark

\lemma \comppisonB On $p_\infty^{-1}(B_0)$, we have the identity
$$
\pi_1^+-\pi_0^- = -\arg a +\pi -\frac{\arg a}{d-1}= \pi - \frac d{d-1}\arg a.
$$
\endstat

\proof
The parametrized path $t \mapsto \bvec c{te^{-i\alpha}}, \ t >0$, thought of as a path in
$X_\infty$, approaches a specific point of $B_0$ with coordinate $c$. Thought of as a path in 
$\orb(X_\infty, D_\infty)$, it approaches the point
$x$  above $c$ where $\Psi_0=\alpha+ (\arg a)/(d-1)$. Thus $\pi_0^-(x) = \alpha+(\arg a)/(d-1)$.

To compute $\pi_1^+(x)$, apply $H$ to the path: the path $t \mapsto \bvec {p(c)-ate^{-i\alpha}}{c},\ t
>0$ approaches $\orb(H_\infty)(x) \in p_\infty^{-1}(B_{-1})$. Just as we needed the
argument of $1/y$ to compute $\Psi_0$ (adjusted by $(\arg a)/(d-1)$), we need the argument of $1/x$
(unadjusted) to compute
$\Phi_0(H(x) = \Phi_1(x)$; clearly this argument is $\alpha-\arg a +\pi$. 
Thus $\pi_1^+(x) = \alpha-\arg a +\pi$.
\qed

This means that on this solid torus, the two functions $\pi_1^+$ and $\pi_0^-$ differ by the constant
$\frac d{d-1}\arg a$. 

\lemma \comppisonothertori The solid torus $p_\infty^{-1}(B_0)$ is the critical locus of $\pi_1^+,\pi_0^-$.
\endstat

\proof We will only outline how to do this, for points above $L_1$.  Choose as above a curve in $X_\infty$
tending in 
$\orb(X_\infty, D_\infty)$ to a point above a point $c \in L_i$. For instance, $t \mapsto \bvec
{t^{-i\alpha}}{cte^{-2i\alpha}}, \ t >0$ is a curve approaching a point $x$ above $L_1$. On this point, we
have $\pi_0^-(x) = 2\alpha-\arg c+ (\arg a)/(d-1)$.  If we apply $H$ to this curve and compute $\arg
(1/x)$, we find 
$$
\pi_1^+(x) = \cases d\alpha&\quad \text{if $d>2$};\\
2\alpha-\arg(1-ac)&\quad \text{if $d=2$}.
\endcases
$$
Since $\arg c$ shows up explicitly in the formula for $\pi_1^+-\pi_0^-$, we see that such a point is not
critical.
\qed

\remark {} The point $ac=1$ in the case $d=2$ above corresponds to $B\cap L_1$; a similar point will show
up above $L_{d-1}$ for every $d$.  

We now need to see that the number $\frac d{d-1}\arg a$ from Lemma \comppisonB\ is (almost) a conjugacy
invariant of
$\orb(H_\infty)$. 

Suppose that $F:\orb(D_{1, \infty}) \to \orb(D_{2, \infty})$ conjugates
$\orb(H_{1,\infty}) \to \orb(H_{2,\infty})$, where we have used indices $1$ and $2$ to distinguish the
objects created from $H_1$ and $H_2$. Then $F$ must send $\Sigma_1^\pm$ to $\Sigma_2^\pm$ since these
are the closures of the set of periodic points of  $\orb(H_{1,\infty})$ and $\orb(H_{2,\infty})$
respectively, so it must also send incompressible tori in the complement of the solenoids
$\Sigma_1^\pm$ to incompressible tori in the complement of the solenoids $\Sigma_2^\pm$. Since $T_{\bold
p_i}$ separates the $T_{\bold p_j}$ with  $j>i$ from the $T_{\bold p_j}$ with $j<i$, we
see that the order of the tori must be preserved, and there exists $k$ such that
$$
F(T_{1, \bold p_i}) \quad \text{is isotopic to}\quad T_{2, \bold p_{i+k}}
$$
By composing $F$ with $\orb((H_2)'_\infty)^{\circ k}$, we may assume that $k=0$, and that  
$$
F(T_{1, \bold p_i}) \quad \text{is isotopic to}\quad T_{2, \bold p_{i}}
$$

\lemma \piscoincide The functions  $\pi^+_{1,i}$ and $\pi^+_{1,i}\circ F$ must  differ by a multiple
of $2\pi/(d-1)$ on their common domain of definition.
\endstat 

\proof By composing $\pi^+_{i,1}$ with an appropriate multiple of $2\pi/(d-1)$, we may assume that if $x_1
\in \Sigma_1^+$ is the fixed point of $\orb((H_1)'_\infty)$ with $\pi^+_{i,1}(x_1)=0$, the $F(x_1)=x_2$ is
the fixed point of $\orb((H_2)'_\infty)$ with $\pi^+_{i,1}(x_2)=0$.  After this change, we must show that
$\pi^+_{1,i}$ and $\pi^+_{1,i}\circ F$ coincide on their common domain.

Choose $j$ sufficiently small so that 
$$
T_{\bold p_{2,j}}^+\  \subset\ T_{\bold p_2,i}^+ \cap F(T_{\bold p_1,i}^+).
$$
Now both diagrams
$$
\CD
T^+_{\bold p_{2,j}}  @>{\orb(H_{2,\infty})}>> T^+_{\bold p_{2,j}} \\ 
@V{\pi_{2,i}^+}VV   @V{\pi_{2,i}^+}VV     \\ 
\R/2\pi \Z  @>{\theta \mapsto d\theta}>>  \R/2\pi \Z
\endCD
\qquad\qquad
\CD
T^+_{\bold p_{2,j}}  @>{\orb(H_{2,\infty})}>> T^+_{\bold p_{2,j}} \\ 
@V{\pi_{i,1}^+\circ F^{-1}}VV   @V{\pi_{i,1}^+\circ F^{-1}}VV     \\ 
\R/2\pi \Z  @>{\theta \mapsto d\theta}>>  \R/2\pi \Z
\endCD
$$ 
commute.  The uniqueness statement \cite{HO1}, Theorem 3.1 isn't quite enough to guarantee that the
verical arrows coincide, since they aren't of degree 1.  In that case, the proof guarantees that such maps
differ by a multiple of $1/(d^{i-j}-1)$. This  is enough to guarantee that if
$\pi_{2,i}^+$ and ${\pi_{i,1}^+\circ F^{-1}}$ coincide at a single point, then they coincide everywhere;
indeed, they do coincide at $x_2$.  Now 
$$
\pi_{2,i}^+= \frac 1{d^m} \circ \pi_{2,i}^+ \circ \orb(H_{2,\infty}) \quad \text{and}\quad
\pi_{i,1}^+\circ F^{-1}= \frac 1{d^m} \circ \pi_{2,i}^+ \circ \orb(H_{2,\infty}); 
$$
any difference comes from different branches of $\frac 1{d^m}$.  Since $F(T_{\bold p_{1,i}})$ is isotopic
to $T_{\bold p_{2,i}}$, and we can choose a branch continuously during the isotopy, we see that indeed
$\pi^+_{2,i}$ and $\pi^+_{1,i}\circ F^{-1}$ (the latter adjusted at the beginning of the proof) must agree
on their common domain of definition. 
\qed

Now we need to know something about this common domain of definition.

\lemma \toriintersect The interiors of the torus $F(T_{\bold
p_{1,i}}^+)$ and the interior of the torus $\Cal B_{i-1}$ corresponding to $B_{2,i-1}$ must have non-empty
intersection.
\endstat

\proof An alternative way of saying this is to say that $\Cal B_{i-1}$ cannot be isotoped, in
$\orb(D_{2,\infty})-(\Sigma_2^+ \cup \Sigma_2^-)$, to a torus outside of $F(T_{\bold p_{1,i}}^+)$.  This can
be seen from linking numbers.  The presumed isotopy will take place in the complement of $T_{\bold
p_{2,j}}^+$ for $j$ sufficiently small.  All curves (or unknotted solid tori) outside $T_{\bold p_{2,i}}^+$
have linking number some integer multiple of $d^{i-j}$ with $T_{\bold p_{2,j}}^+$.  But $\Cal B_{i-1}$
has linking number $d^{i-j-1}$ with $T_{\bold p_{2,j}}^+$ (See Figure \linkingnumber).

\centerline{\BoxedEPSF{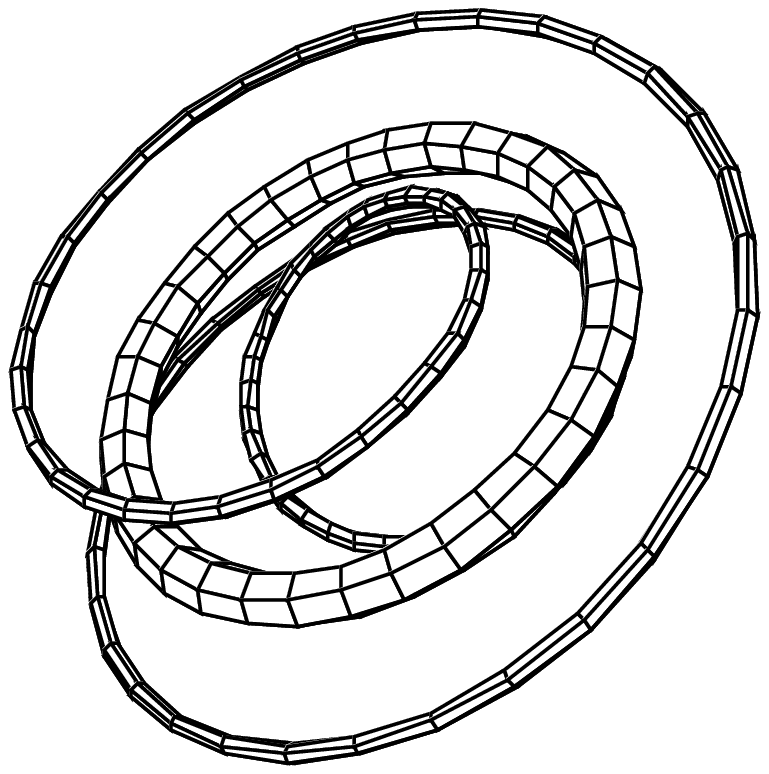 scaled 600}}

\longfigcap \linkingnumber {Represented here is a $3,1$-curve on the boundary of a
solid torus, corresponding to 3 times the generator of the homology of the solid
torus, and the core curve of the solid torus. Clearly they link with linking number 1,
whereas any curve outside the torus links with the
$3,1$-curve with linking number some multiple of 3. Thus this figure represents the case
$i-j=1$ and $d=3$.}

 Since the linking number
must be constant during the isotopy, this is a contradiction.
\qed

Thus $F$ must map some open subset of the torus corresponding to $B_{1,i}$ to some open subset of the torus
corresponding to $B_{2,i}$, in such a way that 
$$\pi_{2,i}^\pm \circ F = \pi_{1,i}
$$
up to a multiple of $1/(d-1)$. This proves Theorem \conjugacyofexthens.
\qed

\remark It seems likely that the condition $\arg a_1 = \arg a_2$ is also sufficient for
conjugacy. We will explore this in a future paper.  It turns out that the
conjugacy properties of mappings like
$\orb((H)_\infty)$ (or the mapping $h_d$ of \cite{HO1}) is quite subtle; and that there is an
infinite-dimensional moduli space, even when the maps are hyperbolic on a neighborhood of
the solenoids.

%% file: Veselov.milnorconj
\heading XI. The compactification of compositions of H\'enon mappings \endheading

A theorem of Friedland and Milnor asserts that any polynomial automorphism of
$\C^2$  is either elementary, in the sense that we can find one variable which
depends only on itself, or conjugate to a composition of H\'enon maps.  Therefore
understanding the appropriate compactification of $\C^2$ to which such a
composition extends is evidently important.

Milnor, in a personal communication, suggested what the 3-sphere at
infinity should look like; we will now state and prove this conjecture.

Let 
$$
H_i:\bvec xy \mapsto \bvec {p_i(x)-a_i y}x ,\quad i=1\cdots k
$$ be $k$ H\'enon mappings, with $a_i \ne 0$ and $p_i$ of degree $d_i\ge 2$. We will
consider $G= H_k\circ \cdots \circ H_1$, which is a polynomial mapping of algebraic
degree $d=d_1\cdots d_k$.

Recall that $\Sigma_d=\underset \leftarrow \to {\lim} (\R/\Z, t \mapsto dt)$ is
the $d$-adic solenoid, and $\sigma_d:\Sigma_d \to \Sigma_d$ the map induced
by $t \mapsto dt$.

We will call the simplest link of two circles with linking number $d$ the one formed by the
circles
$$
\pmatrix \cos t\\ \sin t\\ 0\endpmatrix,\ 0 \le t \le 2\pi \quad \text{and} \quad
\pmatrix (1+ \frac 12 \cos dt) \cos t\\ (1+ \frac 12 \cos dt) \sin t\\ \frac 12 \sin
dt\endpmatrix
$$
\bigskip
\centerline{\BoxedEPSF{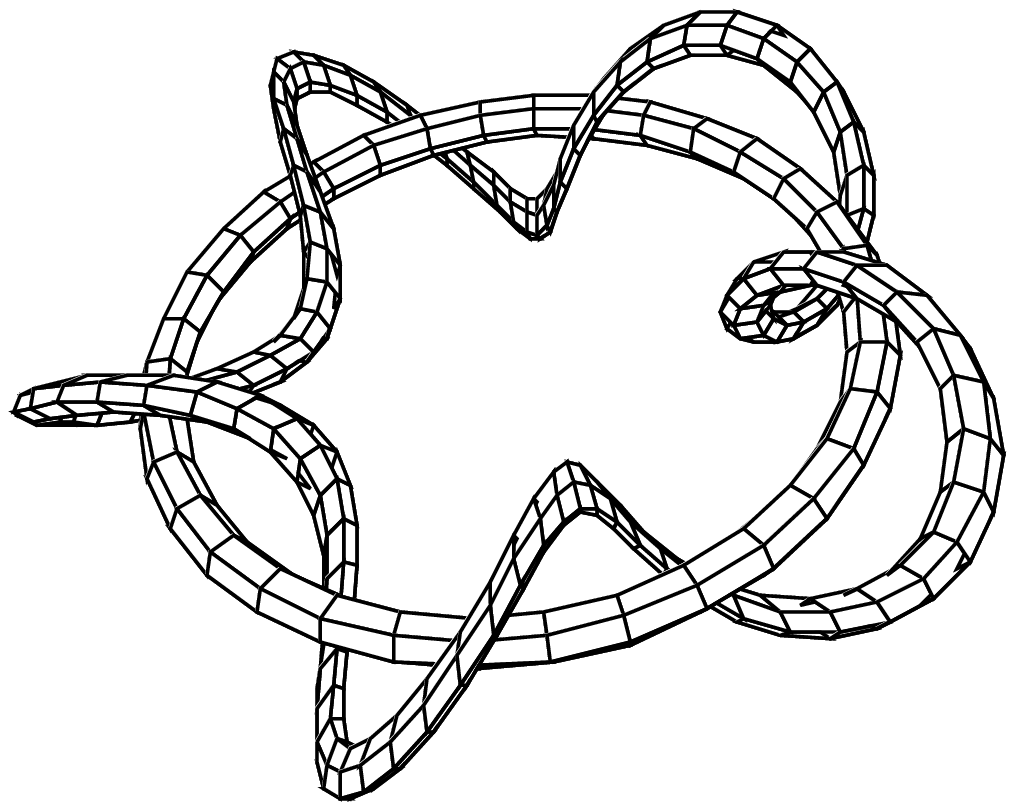 scaled 450}}
\nobreak
\longfigcap \linkedcirclespic {The simplest link of two circles linking with linking number 5.}

\theorem \milnorconj a) There exists a topology on $\C^2 \sqcup S^3$ homeomorphic
to the 4-ball, with $S^3$ corresponding to the boundary, such that $G$ extends continuously to
a homeomorphism of
$g:S^3\to S^3$.

b) The homeomorphism $g$ has two invariant solenoids $\Sigma^+, \Sigma^-$,
one attracting and one repelling, and both homeomorphic to $\Sigma_d$,
and the homeomorphisms can be chosen to be conjugacies between the
restriction $g|_{\Sigma^+}$ and $\sigma_d$, and between $g|_{\Sigma^-}$ and
$\sigma_d^{-1}$.

c) The complement $M=S^3-(\Sigma^+ \cup \Sigma^-)$ has a decomposition
by incompressible tori $(T_i)_{i \in \Z}$, unique up to isotopy, into pieces $M_i$ bounded
by $T_i$ and $T_{i+1}$, homeomorphic to the complement of the simplest link of two circles
with linking number $d_{i\mod k}$. Moreover, $M_i \cap M_j =
\emptyset$ unless $|i-j|\le 1$. The tori can be chosen so that
$g(T_i) = T_{i+k}$.
\endstat

In particular, the topology of the sphere at infinity is different for a composition of
H\'enon maps with total degree $d$ and for a single such mapping: the solenoids are the same
but they are embedded differently in the 3-sphere. 

\proof Let $\tilde X_G$ be the minimal blow-up  of $\Proj^2$ on which
$\tilde G:\tilde X_G \to \Proj^2$ is well-defined.

It can be constructed as follows: set $G_m=H_m \circ \dots \circ H_1$, so
that $G_1=H_1$ and $G_k = G$, and define $\tilde X_{G_m}$ to be the minimal
blow-up of $\Proj^2$ on which $\tilde G_m:\tilde X_{G_m} \to \Proj^2$ is
well-defined. Further denote by $\pi_{G_m}:\tilde X_{G_m} \to \Proj^2$ the
canonical projection, and $\tilde D_{G_m}= \pi_{G_M}^{-1} (l_\infty)$ the divisor at
infinity of $\tilde X_{G_m}$.

 We will construct $\tilde X_{G_m}$ by induction.
Clearly $\tilde X_{G_1} = \tilde X_{H_1}$ is the space constructed in Section III.  Suppose we
have constructed $\tilde X_{G_{m-1}}$, together with $\tilde G_{m-1}$ and $\pi_{G_{m-1}}$.

Set $\tilde X_{G_m}$ to be such that the upper left-hand square of the
diagram 
$$
\CD
\tilde X_{G_m} @>>> \tilde X_{H_m} @>{\tilde H_m}>> \Proj^2\\
@VVV  @V\pi_{H_m}VV\\
\tilde X_{G_{m-1}} @>{\tilde G_{m-1}}>> \Proj^2\\
@V{\pi_{G_{m-1}}}VV\\
\Proj^2
\endCD
$$
is a fiber product in the category of analytic spaces.  Then the top line is a mapping $\tilde
G_m:\tilde X_{G_m} \to \Proj^2$, whereas the left-hand column represents $\tilde
X_{G_m}$ as a modification of $\Proj^2$ at $\bold p$. Thus $\C^2$ is dense
in $\tilde X_{G_m}$, and it is clear by induction that $\tilde G_m$
extends $G_m:\C^2 \to \C^2$. That it
is the minimal modification of $\Proj^2$ to which $G_m$ extends follows from the fact that the
construction of Section III is the minimal modification to which the individual
H\'enon maps extend. Thus $\tilde X_{G_m}$ is the required minimal blow-up.

The divisor above infinity 
$$
\tilde D_G = \tilde X_G - \C^2 = \pi_G^{-1}(l_\infty)
$$
looks as follows. 
\bigskip
\centerline{\BoxedEPSF{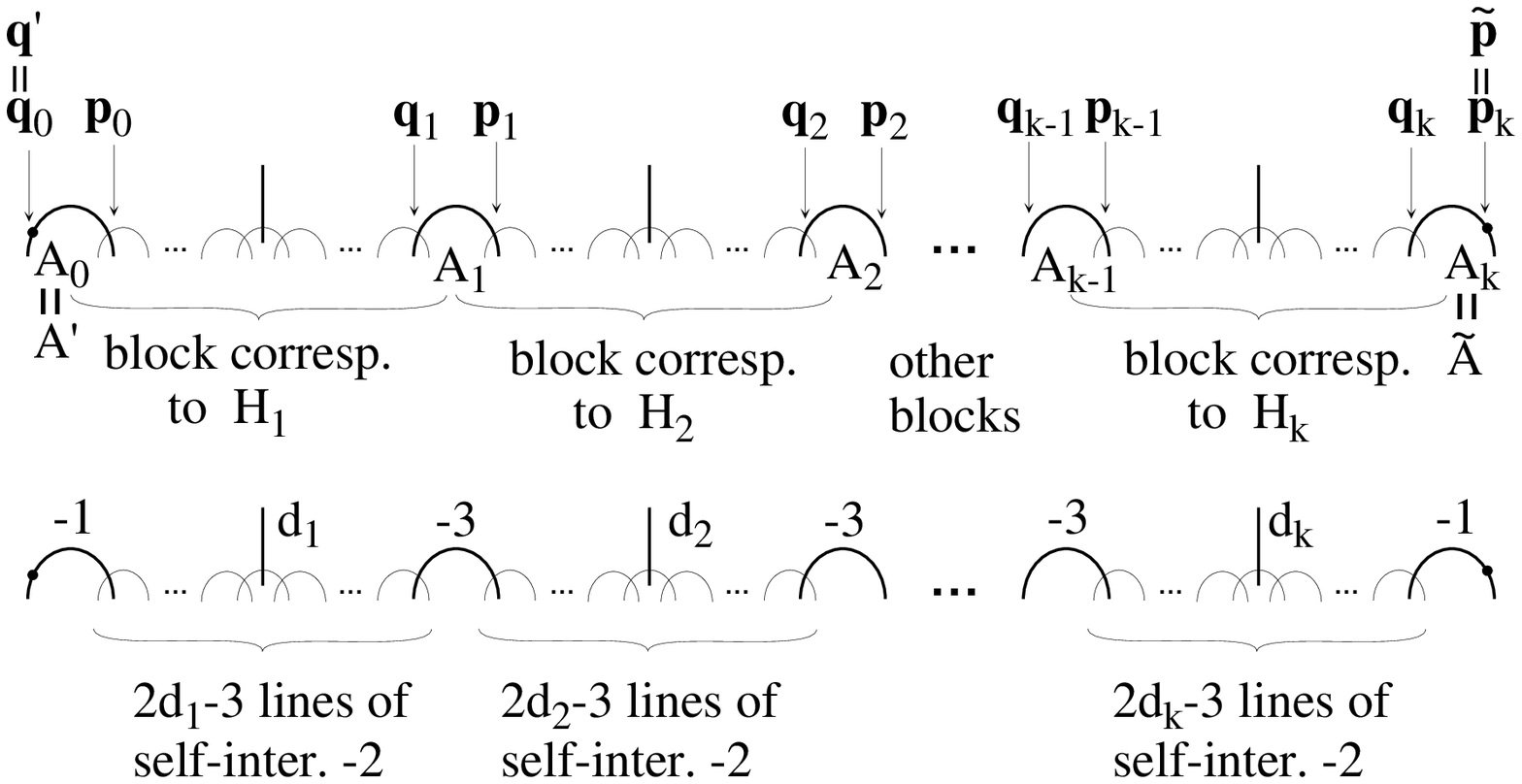 scaled 600}}   
\bigskip
\longfigcap \compHdiv {The divisor $\tilde D_G$. The top figure gives the labels of all the
components and points to which we will need to refer, the second gives the
self-intersections of all these components.}

As before, we will avoid an infinite sequence of blow-ups by a considering a sequence space.
Consider the rational mapping $G^\sharp:\tilde X_G \ratto \tilde X_G$, which is $\tilde G$
wherever it is defined, and define
$$
 \overline \Gamma_G \subset\tilde X_G\times \tilde X_G
$$ 
to be the closure of the graph $\Gamma_{G^\sharp} \subset \tilde X_G \times
\tilde X_G$ of
$G$.  

\lemma {\closureofgraph} A pair $(x,y)$ belongs to $\overline \Gamma_G$ if and only
if either 

\noindent$\bullet$ it is in $\Gamma_{G^\sharp}$, or 

\noindent$\bullet$  $x=\tilde {\bold p} , y \in (\tilde  D-A')\cup\{\bold p'\}$.
\endstat

The proof is analogous to that of Theorem \tildeHasratmap.

Now define the natural compactification of the composition of H\'enon
mappings $G$ as
$$
X_{\infty}(G) = \Bigl\{(\dots, x_{-2}, x_{-1},x_0,x_1,x_2, \dots) \in 
\left(\tilde X_G\right)^\Z | (x_n,x_{n+1}) \in \overline \Gamma_G
\quad \text{for all $n \in \Z$}\Bigr\}.
$$ 
Using Lemma \closureofgraph, this space is not so difficult to
understand.

\proposition \structcompactificationforcomp The space $X_{\infty}(G)$ is compact.  The
complement of two points $\bold q^\infty= (\dots, \bold q', \bold q', \bold q', \dots)$ and 
$\bold p^\infty= (\dots, \tilde{\bold p}, \tilde{\bold p}, \tilde{\bold p}, \dots)$ is an
algebraic manifold.
\endstat

\proof The proof is the same as that of Proposition
\infproductgoodpointssmooth.

Again, to understand the structure of the bad points, we will pass to the real
oriented blow-ups.  

We can define spaces $X_{[-N,M]}(G)$ analogously to the construction 
in Section V, with the
divisors $D_{[-N,M]}(G)$; next we construct  the real oriented blow-ups
$$
\orb(X_{[-N,M]}(G),D_{[-N,M]}(G))
$$
and take their projective limit
$$
\orb(X_\infty(G), D_\infty(G)) = \varprojlim_{N \to \infty} (\orb(X_{[-N,N]}(G),D_{[-N,N]}(G)).
$$
Note that again we are using Propositions \divrightasideal\ and \pitildeexistsfordblepoint\ to
construct the mappings implicit in the projective limit.

The pair $\orb(D_\infty(G)) \subset\orb(X_\infty(G), D_\infty(G))$ is the
compactification of $\C^2$ promised in Theorem \milnorconj. By Theorem \topofprojlim, the pair
is homeomorphic to $(B^4,S^3)$, exactly as in the proof of Theorem
\mainthmrealblowup.  Moreover, the inclusion of $\C^2\subset \orb(X_\infty(G), D_\infty(G))$
and the extension of $G$ to the real oriented blow-up is constructed exactly as in 
\realextensiontheorem.

The proof of (b) is closely analogous to Proposition \computeanglesunderpitilde, but requires
a bit of terminology. First, label components and points of $\tilde D_\infty$ as follows: let
$A_0=l_\infty
\subset
\Proj^2$, and define recursively
$A_i\subset \tilde X_{G_i}$ to be the component of
$H_i^{-1} (A_{i-1})$ which is sent by $H_i$ isomorphically onto $A_{i-1}$. Finally, Let $\bold
p=\bold p_0$ and $\bold q=\bold q_0$; by induction each
$A_i$ contains the two points $\bold p_i, \bold q_i$ which under the isomorphism $G_i
|_{A_i}$ map to
$\bold p_{i-1}$ and $\bold q_{i-1}$, as in Figure \compHdiv.

Next, we will label $\bold p_{m,i}$ and $\bold q_{m,i}$ the points of $D_\infty(G)$
whose $m$th entry is $\bold p_i$, this requires a bit of care when $i=0$ and $i=k$,
which we will leave to the reader. 

Proposition
\computeanglesunderpitilde\ tells us that there are natural angles $\theta_j$ parametrizing
$p^{-1}(\bold p_j)\subset \tilde X_{G_j}$, and that these angles correspond under the H\'enon
mappings (see Equation (\argumentsconseq), where this angle appears as the argument of $v$ and
the argument of $X_1$). Note that we are considering these fibers at the moment when they are
created by the blow-up, so that each lies above a simple point of the divisor
defined so far.  Moreover, the same proposition (specifically, see Equation 
(\argumentsaddeq)) says that the composition of the H\'enon mappings takes
angles $\theta_j$ to angles $\theta_{j-1}$ as indicated in the following diagram:
$$\eqalign{
\theta_k\  \overset {\text {blow-down}}\to {\underset{\tilde X_{G_k} \to \tilde
X_{G_{k-1}}}\to \mapsto}\  &d_k\theta_{k} +
\arg a_k
\overset {\text {blow-down}}\to {\underset{\tilde X_{G_{k-1}} \to \tilde
X_{G_{k-2}}}\to \mapsto}\  d_kd_{k-1}\theta_{k} + d_{k-1}\arg a_k+ \arg a_{k-1}\cr
\noalign{\smallskip}
&\overset {\text {blow-down}}\to {\underset{\tilde X_{G_{k-2}} \to \tilde
X_{G_{k-3}}}\to \mapsto} \quad \dots\quad 
\overset {\text {blow-down}}\to {\underset{\tilde X_{G_1} \to \tilde
X_{G_0}}\to \mapsto} d \theta_k+ \beta,}
$$
where $d=d_k\dots d_1$ and  
$$\beta= d_{k-1} d_{k-2} \dots d_1 \arg a_k + d_{k-2} d_{k-3} \dots d_1 \arg a_{k-1} +
\dots + d_1 \arg a_2 + \arg a_1.
$$

An analogous argument, using appropriate conjugates of the inverses of the $H_j$, will show
that the similar parameter $\phi_j$ of $p^{-1}(\bold q_j)$ is simply multiplied by $d$. Now
the proof ends in the same way as the proof of \computeanglesunderpitilde, showing that the
fibers above $\bold p^\infty$ and $\bold q^\infty$ are both $d$-adic solenoids, in one case
using the angles $\phi_n$, and in the other $\psi_n=\theta_n+\beta/(d-1)$.

Part (c) has substantially already been proved: The tori corresponding to the $\bold p_{n,j}$
do form a sequence of incompressible tori in $\orb(D_\infty(G))-(\Sigma^+ \cup \Sigma^-)$, and
the components of the complements are homeomorphic to the simplest link of two circles which
link with linking number $d_j$.  Moreover, the proof we have given in Theorem
\toroidaldecompofsphere\ that this is the unique toroidal decomposition of 
$\orb(D_\infty(G))-(\Sigma^+ \cup \Sigma^-)$ goes through without change. \qed

%% file: bibliography
\heading Bibliography \endheading

\bigskip
\noindent[B] E. Bedford, {\it Iteration of polynomial automorphisms of $\C^2$}, Proc. Int. Cong. Math.
Kyoto, (1990), 847-858.

\noindent[BS1] E. Bedford and J. Smillie, {\it Polynomial diffeomorphisms of $\C^2$ I: Currents,
equlibrium measure and hyperbolicity\/}, Invent. Math.,
87 (1990), 69-99. 
 
\noindent[BS2] E. Bedford and J. Smillie, {\it Polynomial diffeomorphisms of $\C^2$ II: Stable manifolds
and recurrence\/}, J Amer. Math. Soc., 4 (1991), 657-679. 
 
\noindent[BS3] E. Bedford and J. Smillie, {\it Polynomial diffeomorphisms of $\C^2$ III: Ergodicity, 
exponents and entropy of the equilibrium measure\/}, Math. Ann., 294 (1992), 395-420. 
 
\noindent[BS5] E. Bedford and J. Smillie, {\it Polynomial automorphisms of $\C^2$ V: Critical points
and Lyapunov exponents\/}. (to appear in the Jornal of Geometric Analysis) 
 
\noindent[BS6] E. Bedford and J. Smillie, {\it Polynomial diffeomorphisms of $\C^2$ VI: Connectivity of $J$}.
(submitted for publication)

\noindent[BS7] E. Bedford and J. Smillie, {\it Polynomial diffeomorphisms of $\C^2$ VI: Hyperbolicity and 
external rays}. (submitted for publication)

\noindent[BLS1] E. Bedford, M. Lyubitch and J. Smillie, {\it Distribution of periodic points of polynomial
diffeomorphisms of $\C^2$\/}, Inv. Math 114 (1993), 2, 277-288.
 
\noindent[BLS2] E. Bedford, M. Lyubitch and J. Smillie, {\it Polynomial diffeomorphisms of $\C^2$ IV: The
measure of the maximal entropy and laminar currents\/}, Invent. Math. 112 (1993), 1, 77-125.
 
\noindent[Dil] J. Diller, {\it Dynamics of birational maps of $\Proj^2$\/}, Ind. J.
Math. 45.3 (1996), 721-772.
 
\noindent[Dl] G. Dloussky, {\it Une construction \'el\'ementaire des surfaces d'Inoue-Hirzebruch\/},
Math. Annalen 280 (1988), 663-682.

\noindent[DO] G. Dlousky and K. Oeljeklaus, {\it Foliations on compact surfaces containing global
spherical shells\/}, preprint of CMI, Universit\'e de Provence, 1997.

\noindent[Dou] R. and A. Douady, {\it Alg\`ebre et Th\'eorie Galoisiennes, I and II}, Cedic/Fernand
Nathan, Paris 1977, 1979.

\noindent[DouN] A. Douady and L. H\'erault, {\it Arrondissement des vari\'et\'es \`a coins\/},
Comment. Math.  Helvet. 48 (1973), 484-491. (Appendice to A. Borel and J-P. Serre {\it Corners and
Arithmetics Groups\/}, 436-483.)
 
\noindent[FS1] J. Fornaess and N. Sibony, {\it Complex H\'enon mappings in $\C^2$ and Fatou-Bieberbach
domains\/}, Duke Math. Jour. 65 (1992), 345-380.

\noindent[FS2] J. Fornaess and N. Sibony, {\it Complex Dynamics in Higher Dimension I\/}, Asterisque
222 (1994), 201-231.

\noindent[FS3] J. Fornaess and N. Sibony, {\it Complex Dynamics in Higher Dimension II\/}, Ann.
Math. Studies 137 (1995), 134-182.

\noindent[FS4] J. Fornaess and N. Sibony, {\it Complex Dynamics in Higher Dimension\/}, in Complex
Potential Theory (1995) Kluwer Academic Publishes, 131-186.

\noindent[Fr1]  S. Friedland, {\it Oral communication\/}, 1995.

\noindent[Fr2]  S. Friedland, {\it Entropy of polynomial and rational maps\/}, Ann. Math, 133 (1991), 359-368.

\noindent[FM] S. Friedland and J. Milnor, {\it Dynamical properties of polynomial automorphisms\/},
Ergod. Th. and Dyn. Sys., 9 (1989), 67-99.

\noindent[G] H. Grauert, {\it Uber Modifikationen und exzeptionelle analytische Mengen\/}, Math.
Analen 146,(1962), 331-368. 

\noindent[Gri] P. Griffiths, {\it Infinite Abelian group theory\/}, Chicago
Lectures in Mathematics, U. of Chicago Press, 1970.

\noindent[GH] P. Griffiths and J. Harris, {\it Principles of Algebraic Geometry\/}, Wiley, N.Y., 1978. 

\noindent[Gu] R. Gunning,{\it Introduction to Holomorphic Functions in Several Variables, III:
Homological theory\/}, Wadsworth \& Brooks/Cole, Belmont, CA (1990).

\noindent[Har] R. Hartshorne, {Algebraic Geometry\/}, Springer Verlag

\noindent[Hat1] A. Hatcher, {\it Topology}, to appear.

\noindent[Hat2] A. Hatcher, {\it Notes on 3-dimensional topology}, manuscript in
preparation.

\noindent[Hemp] J.  Hempel, {\it 3-Manifolds\/}, Annals of Math. Stud. 86, Princeton Univ.
Press,  1976.

\noindent[He] M. H\'enon, {A two-dimensional mapping with a strange attractor\/}, Comm. Math. Phys., 50
(1976), 69-77.

\noindent[Hi1] H. Hironaka, {\it Desingularization of complex-analytic spaces\/}, Congr\`es Inter.
Math., II (1970), 627-631.

\noindent[Hi2] H. Hironaka, {\it Triangulation of algebraic sets\/}, in {\it Algebraic Geometry\/},
Arcata, 1974, Proc Symp. Pure Math. 29, AMS, Providence, R.I., 1975, 165-185..

\noindent[Hirz] F. Hirzebruch, {\it Hilbert Modular Surfaces}, l'Ens. Math. 19 (1973), 183-281.

\noindent[HNK] F. Hirzebruch, W.D. Neumann and S.S. Koh, {\it Differentiable manifolds and quadratic
forms\/}, Marcel Dekker, N.Y. 1971.

\noindent[Ho1] H. Hopf, {Uber die Abbildungen der 3-Sph\"are auf die Kugelfl\"ache\/}, Math. Annalen
104 (1931), 637-665.   

\noindent[Ho2] H. Hopf, {Zur Topologie der komplexen Manigfaltigkeiten \/}, Studies and Essays
presented to R. Courant on his 60th birthday, Interscience, N.Y. (1948).   
 
\noindent[Hu] J.H. Hubbard, {\it The H\'enon mapping in the complex domain\/}, in {\it Chaotic
Dynamics and Fractals\/}, Barnsley and Demko, eds, Academic Press, New York (1986), 101-111.  

\noindent[HO1] J.H. Hubbard and R.W. Oberste-Vorth, {\it H\'enon mappings in the complex domain I: The
global topology of dynamical space}, Pub. Math. IHES 79 (1994), 5-46. 

\noindent[HO2] J.H. Hubbard and R.W. Oberste-Vorth, {\it H\'enon mappings in the complex domain II:
Projective and inductive limits of polynomials\/}, in {\it Real and Complex Dynamical
Systems\/}, Branner and Hjorth, eds, Kluwer Academic Publishers (1995), 89-132. 

\noindent[HO3] J.H. Hubbard and R.W. Oberste-Vorth, {\it H\'enon mappings in the complex domain III:
the case of a single attracting fixed point\/}, in preparation.

\noindent[JP] J.H. Hubbard and P. Papadopol, {\it Superattractive fixed points in $\C^n$\/}, Ind. J.
Math. 43.1 (1994), 321-365.

\noindent[L] S. Lang, {\it Algebra\/}, Addison-Wesley, 1965.

\noindent[Mi1] J. Milnor,  {\it Singular points of complex hypersurfaces}, Ann. Math. Stud. 61, Princeton
University Press (1968).

\noindent[Mi2] J. Milnor,  {\it Personal communication\/}, 1995.

\noindent[Mi3] J. Milnor,  {\it Non-expansive H\'enon maps}, Adv. Math., 69 (1988), 109-114.

\noindent[Mi4] J. Milnor,  {\it Topology from the differentiable viewpoint}, The University
Press of Virginia, 1965.

\noindent[S1] I.R. Shafarevic, {\it Algebraic surfaces\/}, Proc. Stek. Inst. Math., 75 (1965), AMS
Providence (1967).

\noindent[S2] I.R. Shafarevic, {\it Basic Algebraic Geometry, I\&II\/}, Springer Verlag, NY 1994.

\noindent[S] J. Smillie, {\it The entropy of polynomial diffeomorphisms of  $\C^2$\/}, Ergod. Th. Dyn. Sys. 10
(1990), 823-827.

\noindent[Spa] E. Spanier, {\it Algebraic topology\/}, McGraw Hill, N.Y. 1966.

\noindent[Th] W.P. Thurston, {\it Three-Dimentional Geometery and Topology\/}, Vol. I, Princeton,
1997.

\noindent[Or] P. Orlik, {\it Seifert manifolds\/}, LNM, Nr. 291, Springer, 1972.

\noindent[Ves] V. Veselov, {\it A compactification of the Fatou mapping as a dynamical system\/}, Thesis,
University of South Florida, 1996.